\documentclass[aos,preprint]{imsart}

\usepackage[]{graphicx} % more modern %take out [draft] to restore figures

\usepackage{amsfonts,float}
\usepackage{amssymb}

\DeclareSymbolFont{sfletters}{OML}{cmbrm}{m}{it}

\DeclareMathSymbol{\salpha}{\mathord}{sfletters}{"0B}
\usepackage{mathrsfs}
\RequirePackage{natbib}
\RequirePackage[colorlinks,citecolor=blue,urlcolor=blue]{hyperref}
\usepackage{xr}
\usepackage{xcite}
\usepackage{overpic}
\makeatletter
\newcommand{\oset}[3][1ex]{%
  \mathrel{\mathop{#3}\limits^{
    \vbox to#1{\kern-2\ex@
    \hbox{$\scriptstyle#2$}\vss}}}}
\makeatother

\usepackage{mathtools}
\RequirePackage{natbib}
\usepackage{comment}
\usepackage{enumerate}
\usepackage{eufrak}

\usepackage{array,varwidth}

\usepackage{overpic}
\usepackage[framemethod=tikz,innertopmargin=7pt]{mdframed}
\usepackage{graphicx} % more modern
\usepackage{subcaption}
\usepackage{MnSymbol}
\usepackage{placeins} %needed for floatbarrier to control mix of figs and biblio at end

\makeatletter

\def \@cflt{%
    \let \@elt \@comflelt
    \setbox\@tempboxa \vbox{}%
    \@toplist
\setbox\z@\vsplit\@outputbox to 0.5\ht\@outputbox
    \setbox\@outputbox \vbox{%
                             \boxmaxdepth \maxdepth
                             \unvbox\z@
                             \vskip .5\textfloatsep
                             \unvbox\@tempboxa
                             \vskip -\floatsep
                             \topfigrule
                             \vskip .5\textfloatsep
                             \unvbox\@outputbox
                             }%
    \let\@elt\relax
    \xdef\@freelist{\@freelist\@toplist}%
    \global\let\@toplist\@empty
}

\makeatother

\usepackage{booktabs,longtable,tabu} 
\usepackage{multirow} 
\usepackage{float}
\usepackage{makecell}
\usepackage{dcolumn}

\newcommand*\widebar[1]{%
   \vbox{%
     \hrule height 0.5pt%                  % Line above with certain width
     \kern0.25ex%                          % Distance between line and content
     \hbox{%
       \kern-0.15em%                        % Distance between content and left side of box, negative values for lines shorter than content
       \ifmmode#1\else\ensuremath{#1}\fi%  % The content, typeset in dependence of mode
       \kern-0.15em%                        % Distance between content and left side of box, negative values for lines shorter than content
     }% end of hbox
   }% end of vbox
}

\RequirePackage[OT1]{fontenc}
\RequirePackage{amsthm,amsmath}

\startlocaldefs
\numberwithin{equation}{section}
\theoremstyle{plain}
\newtheorem{theorem}{Theorem}[section]
\newtheorem{lemma}{Lemma}[section]
\newtheorem{proposition}{Proposition}[section]
\newtheorem{assumption}{Assumption}[section]

\newtheorem{alg}{Algorithm}
\usepackage{algorithm}
\usepackage{algorithmic}

\theoremstyle{remark}
\newtheorem{remark}{Remark}[section]

\newcommand{\R}{\mathbb{R}}

\renewcommand{\hat}[1]{\widehat{#1}}

\newcommand{\bunderline}[1]{\underline{#1\mkern-4mu}\mkern4mu }
\newcommand{\vunderline}[1]{\underline{#1\mkern-2mu}\mkern2mu }

\newcommand{\mnorm}[1]{{\left\vert\kern-0.25ex\left\vert\kern-0.25ex\left\vert #1 
    \right\vert\kern-0.25ex\right\vert\kern-0.25ex\right\vert}}

\newcommand{\I}{\textbf{\textsc{I}}}
\newcommand{\II}{\textbf{\textsc{II}}}
\newcommand{\IIIo}{\textbf{\textsc{III}}}
\newcommand{\III}{\widetilde{\textbf{\textsc{III}}}}
\newcommand{\IV}{\widetilde{\textbf{\textsc{II}}}}
\newcommand{\V}{\tilde{\textbf{\textsc{I}}}}

\newcommand{\dK}{d_{\textup{K}}}

\newcommand{\A}{A}

\newcommand{\G}{\mathbb{G}}

\renewcommand{\P}{\mathbb{P}}
\newcommand{\E}{\mathbb{E}}

\newcommand{\LL}{\mathcal{L}}

\newcommand{\e}{\epsilon}
\newcommand{\ve}{\varepsilon}

\newcommand{\cov}{\text{cov}}
\newcommand{\tr}{\operatorname{tr}}

\newcommand{\var}{\operatorname{var}}

\newcommand{\ttop}{^{\top}}

\newcommand{\op}{_{\textup{op}}}

\newcommand{\ts}{\textstyle}

\endlocaldefs

\usepackage[T1]{fontenc}
\begin{document}

\begin{frontmatter}
\title{ Bootstrapping the Operator Norm in High Dimensions: Error Estimation for Covariance Matrices and Sketching}
\runtitle{Bootstrapping the Operator Norm}

%\title{ Error Estimation for the Operator Norm in High Dimensions: Bootstrap, 
% Covariance Matrices, and Sketching}

\begin{aug}
\author{\fnms{Miles E.} \snm{Lopes}\thanksref{t1}\ead[label=e1]{melopes@ucdavis.edu}},
\author{\fnms{N. Benjamin} \snm{Erichson}\thanksref{t2}\ead[label=e2]{???}}\\
\and
\author{\fnms{Michael W.} \snm{Mahoney}\thanksref{t2}
\ead[label=e3]{??}}
%\ead[label=u1,url]{http://www.foo.com}}
\affiliation{University of California, Davis\thanksref{t1}, International Computer Science Institute\thanksref{t2}, and University of California, Berkeley\thanksref{t2}}

\thankstext{t1}{Supported in part by NSF grants DMS-1613218 and DMS-1915786}
\thankstext{t2}{Supported in part by ARO, DARPA (FA8750-17-2-0122), NSF, and ONR.}
\runauthor{M.~E. Lopes et al.}

%\thankstext{t3}{???}
%\runauthor{run author.}

\end{aug}

\begin{aug}
%\address{\normalsize University of California, Davis}
%\\
%Usually a few lines long\\
%\printead{e1}\\
%\phantom{E-mail:\ }\printead*{e2}}
%
%\address{Address of the Third author\\
%Usually a few lines long\\
%Usually a few lines long\\
%\printead{e3}\\
%\printead{u1}}
\end{aug}

\begin{abstract}
Although the operator (spectral) norm is one of the most widely used metrics for covariance estimation, comparatively little is known about the fluctuations of error in this norm.
To be specific, let $\hat\Sigma$ denote the sample covariance matrix of $n$ observations in $\R^p$ that arise from a population matrix $\Sigma$, and let $T_n=\sqrt{n}\|\hat\Sigma-\Sigma\|\op$. In the setting where the eigenvalues of $\Sigma$ have a decay profile of the form $\lambda_j(\Sigma)\asymp j^{-2\beta}$, we analyze how well the bootstrap can approximate the distribution of $T_n$. 
Our main result shows that up to factors of $\log(n)$, the bootstrap can approximate the distribution of $T_n$ at the \emph{dimension-free} rate of 
 $ n^{-\frac{\beta-1/2}{6\beta+4}}$, 
with respect to the Kolmogorov metric. Perhaps surprisingly, a result of this type appears to be new even in settings where $p< n$.
More generally, we discuss the consequences of this result beyond covariance matrices, and show how the bootstrap can be used to estimate the errors of sketching algorithms in randomized numerical linear algebra (RandNLA). An illustration of these ideas is also provided with a climate data example.
 \end{abstract}


\end{frontmatter}

\section{Introduction}
Within the areas of covariance estimation and principal components analysis, 
it is of central importance to understand
how well a sample covariance matrix 
%\begin{equation*}
$\hat\Sigma=\ts\frac{1}{n}\sum_{i=1}^n X_iX_i\ttop$
%\end{equation*} 
approximates its population version $\Sigma=\E[X_1X_1\ttop]$, where \smash{$X_1,\dots,X_n\in\R^p$} are centered i.i.d.~observations. 
In particular, a major line of research in high-dimensional statistics has focused on the problem of deriving non-asymptotic bounds for the operator (spectral) norm error
$$T_n=\sqrt{n}\|\hat\Sigma-\Sigma\|\op,$$
where the norm is defined as $\|A\|\op=\sup_{\|u\|_2=1}\|Au\|_2$. A partial overview of work on this problem, as well as some of its extensions, may be found in the papers
(\cite{Rudelson:1999},~\cite{Bickel2:2008},~\cite{Cai:2010}
%Bickel1:2008,
~\cite{Adamczak:2011},~\cite{Tropp:2012},~\cite{Lounici:2014},~\cite{Bunea:Xiao:2015},~\cite{Koltchinskii:Bernoulli},~\cite{Minsker:2017}, among numerous others).

As a whole, this line of work 
offers many conceptual insights into the ways that error is influenced by model assumptions. 
However, the literature is less complete with regard to inference, and there are not many guarantees for the problem of constructing confidence intervals for $T_n$, which is equivalent to constructing numerical bounds on the error of $\hat\Sigma$, or confidence regions for $\Sigma$. Accordingly, the challenges of inference on high-dimensional covariance matrices have stimulated much recent activity, and there has been a particular interest to understand the limits of the bootstrap in this context~(\cite{Johnstone:2018}~\textsection X.C,~\cite{Zhou:2018};\,~\cite{Karoui:2019};~\cite{Lopes:Biometrika};~\cite{Naumov:2019}).

Simultaneously with these developments, the burgeoning field of randomized numerical linear algebra (RandNLA) has generated many other error estimation problems of a similar nature~\citep{mahoney2011randomized,Halko:2011,woodruff2014sketching,Kannan:2017,RandNLA_PCMIchapter_chapter}. A prototypical example deals with computing a fast randomized approximation of the product $A\ttop A$, where $A$ is a very large matrix. Most commonly, the matrix $A$ is randomly ``sketched'' into a much shorter matrix $\tilde A$, which can then be used to quickly compute $\tilde A\ttop \tilde A$ as an approximation to $A\ttop A$.
In turn, it is necessary to assess the unknown error $\|\tilde A\ttop \tilde A-A\ttop A\|\op$,  which leads to a notable parallel with the statistical literature: There are many existing theoretical error bounds, but very few tools for numerical error estimation (cf.~Sections~\ref{sec:related} and~\ref{sec:sketching}). Furthermore, the operator norm is of special importance, because it governs the accuracy of numerous matrix computations, and it frequently appears in numerical analysis~\citep{Golub}.

 Motivated by the challenges above, this paper aims to quantify how well the bootstrap can approximate the error distribution $\mathcal{L}(T_n)$ for sample covariance matrices, and likewise in the context of RandNLA. Specifically, we consider a setup where $\Sigma$ has low ``effective rank'' and its ordered eigenvalues satisfy a decay profile of the form
\begin{equation}\label{eqn:alphacond}
\lambda_j(\Sigma) \asymp j^{-2\beta},
\end{equation}
for some parameter $\beta>1/2$. Variations of this setting have drawn considerable attention in recent years, especially in connection with principal components analysis~(e.g.,~\cite{Lounici:2014},~\cite{Bunea:Xiao:2015}, \cite{Reiss:2016},~\cite{Koltchinskii:Bernoulli,KoltchinskiiAOS:2017}, \cite{Nickl:2017}, \cite{Naumov:2019}, and~\cite{Jung:2018}, among others). Moreover, the condition~\eqref{eqn:alphacond} corresponds to problems where sketching algorithms can be highly effective.

\subsection{Contributions}\label{sec:contributions}
To briefly outline our main result, let the Kolmogorov metric be denoted as $\dK(\mathcal{L}(U),\mathcal{L}(V))=\sup_{t\in\R}|\P(U\leq t)-\P(V\leq t)|$ for two generic random variables $U$ and $V$, and let $T_n^*$ denote the bootstrap version of $T_n$, obtained by sampling with replacement from $(X_1,\dots,X_n)$. Then, as long as~\eqref{eqn:alphacond} is satisfied and the observations have suitable tail behavior, it follows that the bound
\begin{equation}\label{eqn:introbound}
 \dK\Big(\mathcal{L}(T_n) \, , \, \mathcal{L}(T_n^*|X)\Big) \ \leq \ c\, n^{-\frac{\beta-1/2}{6\beta+4}}\log(n)^c
\end{equation}
holds with probability at least $1-\frac{c}{n}$, where $\mathcal{L}(T_n^*|X)$ is the conditional distribution of $T_n^*$ given the observations. (Going forward, we use $c$ to denote a positive constant not depending on $n$ whose value may change at each occurrence.) Most importantly, this bound explicitly relates the structural parameter $\beta$ to the rate of approximation in a way that is both \emph{non-asymptotic} and~\emph{dimension-free}.

From the standpoint of methodology, our work illustrates new possibilities for applying the bootstrap in the domains of computer science and applied mathematics. At this interface, the bootstrap has a largely untapped potential to make an impact, because error estimation allows randomized computations to be done \emph{adaptively}, so that ``just enough'' work is done. More specifically, the estimated error of a rough initial solution can be used to predict how much extra computation is needed to reach a high-quality solution --- and this will be demonstrated numerically in Section~\ref{sec:sketching}.
Lastly, to put this type of application into historical perspective, it is worth noting that the bootstrap has been traditionally labeled as ``computationally intensive'', and so in this respect, it is relatively novel to use the bootstrap \emph{in the service of computation}. 

With regard to theoretical considerations, our work contributes to recent developments on bootstrap methods, as well as covariance estimation. For the bootstrap, we expand upon the 
progress achieved in the series of papers~\citep{CCK:AOS,CCK:suprema,CCK:SPA,CCK:AOP}, which address bootstrap approximations for ``max statistics'' of the form
%\begin{equation}\label{eqn:suprep}
$M_n\, = \ \sup_{f\in\mathscr{F}}\,\G_n(f)$,
%\end{equation}
where  $\mathscr{F}$ is a class of functions, and  $\G_n(f)=\ts\frac{1}{\sqrt n}\sum_{i=1}^nf(X_i)-\E[f(X_i)]$.
The basic similarity between $M_n$ and $T_n$ is that they can be represented in a common form, due to the variational representation of $\|\cdot\|\op$. Nevertheless, the statistic $T_n$ seems to present certain technical obstructions with regard to previous results.
First, in order to handle the metric $\dK$, the mentioned works essentially require a ``minimum variance condition'' such as 
$$\inf_{f\in\mathscr{F}} \var(\G_n(f)) \, \geq \, c,$$
which poses a difficulty in our setting, because the minimum variance may decrease rapidly with $n$.  As a result, a challenge arises in showing that our statistic is well approximated (in $\dK$) by the supremum $\sup_{f\in\mathscr{F}_n'}\G_n(f)$, where $\mathscr{F}_n'\subset\mathscr{F}$ is a ``nice'' subset for which $\inf_{f\in\mathscr{F}'_n} \var(\G_n(f))$ decreases slowly with $n$. 
%Second, the use of discretization to replace an
%
%
Second, further challenges are encountered when controlling the discretization error that comes from replacing $\mathscr{F}$ with a discrete $\e$-net. More specifically, this error is significant in our analysis because the relevant class $\mathscr{F}$ is exponentially larger than VC-type --- in the sense that $\e$-covering numbers grow exponentially in $1/\e$, rather than polynomially. 
By contrast, previous applications of bootstrap approximation results for max statistics have predominantly been concerned with VC-type classes, which allow for strong control of the discretization error.

Another technical aspect of our work deals with dimension-free bounds for $\|\hat\Sigma-\Sigma\|\op$, as studied in~\citep{Rudelson:Vershynin:2007,Oliveira:2010,Hsu:Matrix:2012,Koltchinskii:Bernoulli, Minsker:2017}.  In the setting of~\eqref{eqn:alphacond}, this line of work shows that if the observations satisfy $\|X_i\|_2\leq c$ almost surely, or $\|\langle u,X_i\rangle\|_{\psi_2}\asymp \|\langle u,X_i\rangle\|_{2}$ for all $\|u\|_2=1$,
 then the operator norm error can be bounded as
$\|\hat\Sigma-\Sigma\|\op \ \leq \ c\,n^{-1/2}\log(n)^c$
with high probability. However, the $\ell_2$-boundedness condition is often restrictive, while the $\psi_2$-$L_2$ equivalence condition is not well-suited to the discrete distributions that arise from resampling~\citep[][\textsection \,3.4.2]{Vershynin:2018}. Consequently, as a way to streamline our analysis of both $(X_1,\dots,X_n)$ and the bootstrap samples  $(X_1^*,\dots,X_n^*)$, it is of interest to develop a dimension-free bound that can be applied in a more general-purpose way.
Indeed, an extension of this type is also suggested briefly in the paper~\citep{Rudelson:Vershynin:2007}, but to the best of our knowledge, such a result has not been available in the literature. 
Accordingly, one of our secondary results (Proposition~\ref{prop:newcontraction}) serves this purpose by showing that the approach of~\cite{Rudelson:Vershynin:2007} based on non-commutative Khintchine inequalities can be used to weaken the distributional constraints in a flexible manner.  

\subsection{Related work}\label{sec:related}
The most closely related work to ours is the recent paper~\citep{Zhou:2018}, which studies bootstrap approximations for certain variants of $T_n$. To explain the connection, first recall that $T_n$ may be written in terms of a supremum over the unit sphere $\mathbb{S}^{p-1}\subset\R^p$, namely
$T_n=\sup_{u\in\mathbb{S}^{p-1}}\sqrt n\, |u\ttop (\hat\Sigma-\Sigma)u|$.
As an alternative to this, the paper~\citep{Zhou:2018} analyzes sparse versions of $T_n$ obtained by replacing $\mathbb{S}^{p-1}$ with a subset of vectors that are at most $s$-sparse, $\{u\in\mathbb{S}^{p-1}|\, \|u\|_0\leq s\}$, where $1\leq s\leq p$. For these sparse versions of $T_n$, bootstrap approximation results are obtained in the Kolmogorov metric with rates of the form $s^{9/8}/n^{1/8}$, up to logarithmic factors.
As this relates to our work, it should be emphasized that the setting in~\citep{Zhou:2018} is quite different, since  the eigenvalues of $\Sigma$ are not assumed to decay.
Specifically, the difference becomes most apparent when
 $s=p$, so that the sphere $\mathbb{S}^{p-1}$ coincides with the set $\{u\in\mathbb{S}^{p-1}|\, \|u\|_0\leq s\}$. In this case, the analysis without spectral decay requires $p\ll n^{1/9}$ for bootstrap consistency, whereas our setting places no constraints on $p$.

Next, the recent paper~\citep{Karoui:2019} looks at both positive and negative results for bootstrapping sample eigenvalues.  For the positive results, this work assumes that $\Sigma$ is nearly low-rank, and that the dimension satisfies $p\lesssim n$. The main result shows that the bootstrap consistently approximates the joint distribution of  
$L_n=\sqrt{n}(\lambda_j(\hat\Sigma)-\lambda_j(\Sigma))_{1\leq j\leq j_0}$, where $j_0$ is held fixed as $(n,p)\to\infty$. To mention some further  points of contrast, the distribution of $L_n$ is analyzed in an asymptotic manner by adapting fixed-$p$ results~\citep[e.g.][]{Beran:1985,Eaton:1991}, whereas our approach is non-asymptotic.
 Furthermore, the analysis of $L_n$ relies upon the eigenvalues $(\lambda_j(\Sigma))_{1\leq j\leq j_0}$ having multiplicity 1, whereas the analysis of $T_n$ does not. 
Indeed, this illustrates a key difference between $T_n$ and $L_n$, because in the latter case, it is well-known that repeated eigenvalues are a source of difficulty for bootstrap methods. Also, our numerical results in Section~\ref{sec:simult} will confirm that bootstrapping $T_n$ is robust against high multiplicity. For additional background on this topic, we refer to~\citep{Hall:2009} and references therein.

Two more papers on bootstrap methods for high-dimensional sample covariance matrices are~\citep{Naumov:2019} and~\citep{Lopes:Biometrika}. The first of these deals with bootstrapping the Frobenius norm error of spectral projectors, $\|\hat v_j\hat v_j\ttop -v_jv_j\ttop\|_F^2$, where $\hat v_j$ and $v_j$ are respective $j$th eigenvectors of $\hat\Sigma$ and $\Sigma$. Although the statistic $\|\hat v_j\hat v_j\ttop -v_jv_j\ttop\|_F^2$ is qualitatively different from $T_n$, the paper~\citep{Naumov:2019} shares our interest in settings where $\Sigma$ has low effective rank. Also see~\citep{KoltchinskiiAOS:2017,Koltchinskii:Sankhya}.
 In a different direction, the paper~\citep{Lopes:Biometrika} generalizes the parametric bootstrap for high-dimensional models without spectral decay, and it establishes consistency for linear spectral statistics.

Finally, to conclude this section, we describe related work on the estimation of algorithmic error. Here, it is important to note that error estimation has a long history for \emph{deterministic} algorithms, such as those in numerical partial differential equations and finite-element methods, where it is called \emph{a posteriori error estimation}~\citep[][among many others]{
Babuska1:1978,Babuska2:1978,Verfurth:1994,Becker:2001,Jiranek:2010,Ainsworth:2011,Cangiani:2017}. However, in the literature on randomized algorithms, error estimation has received much less attention, and for certain types of computations there are only a few papers addressing error estimation: ~\citep[low-rank approximation:][]{liberty2007,woolfe2008,Halko:2011}, \citep[least-squares:][]{Lopes:ICML},~\citep[classification:][]{Lopes:AOSRF},~\citep[matrix multiplication:][]{Blum,Sarlos:2006,Lopes:JMLR}. 
Among these works, the only ones to address error estimation for the operator norm are~\citep{liberty2007,woolfe2008,Halko:2011}, but this is done specifically for low-rank approximation, which is complementary to our applications. Also, the approach in these works is quite different from bootstrapping, and is based on the idea of bounding error in terms of random ``test vectors'', which is rooted in the classical works~\citep{Freivalds:1979,Dixon:1983}.
In essence, the main difference between the test-vector approach and bootstrapping is that the former is inherently conservative, whereas the latter can be used to \emph{directly estimate} the error distribution.

\paragraph{Outline}Section~\ref{sec:main} presents the problem setup and main result. Section~\ref{sec:cov} describes numerical results for inference tasks related to covariance matrices, including the construction of simultaneous confidence intervals for population eigenvalues. Section~\ref{sec:sketching} introduces the setting of sketching algorithms, and demonstrates the performance of the bootstrap in synthetic problems, as well as in a climate data example. Lastly, all proofs are given in the supplementary material.

\paragraph{Notation and  conventions}
For a vector $v\in\R^m$, and a number $q\geq 1$, the $\ell_q$-norm is $\|v\|_q=(\sum_{j=1}^m|v_j|^q)^{1/q}$. 
For a real matrix $M$, its Frobenius norm is $\|M\|_F=\sqrt{\tr(M\ttop M)}$, and its Schatten-$q$ norm is $\|M\|_{S_q}=\tr((M\ttop M)^{q/2})^{1/q}$. The identity matrix of size $m\times m$ is $I_m$, and the standard basis vectors in $\R^m$ are $\{e_1,\dots,e_m\}$. The sorted singular values of a real matrix $M$ are written as $\sigma_j(M)\geq \sigma_{j+1}(M)$, and similarly, if $M$ is symmetric, then the sorted eigenvalues are written as $\lambda_{j}(M)\geq \lambda_{j+1}(M)$. For a random variable $\xi$, the $L_q$ norm is $\|\xi\|_q =(\E[|\xi|^q])^{1/q}$. Also, if $\psi_q(x)= \exp(x^q)-1$, then the $\psi_q$-Orlicz norm is given by $\|\xi\|_{\psi_q}=\inf\{r>0\, |\, \E[\psi_q(|\xi|/r)]\leq 1\}$.  If $\zeta$ is another random variable, then the conditional distribution of $\zeta$ given $\xi$ is denoted as $\mathcal{L}(\zeta|\xi)$. 
 If $a_n$ and $b_n$ are sequences of non-negative real numbers, we write $a_n\lesssim b_n$ if there is a constant $c>0$ not depending on $n$, and integer $n_0\geq 1$ such that $a_n\leq c b_n$ for all $n\geq n_0$. In addition, we write $a_n\vee b_n=\max\{a_n,b_n\}$, and $a_n\asymp b_n$ if $a_n\lesssim b_n$ and $b_n\lesssim a_n$.

\section{Main result}\label{sec:main}

Our setup is based on a sequence of models indexed by $n$, where all parameters may depend on $n$, unless stated otherwise. In particular, the dimensions $p=p(n)$ and $d=d(n)$ below may vary with $n$. Lastly, if a parameter does not depend on $n$, then it is understood not to depend on $p$ or $d$ either.

\noindent \begin{assumption}[Data-generating model]\label{A:model}
~\\[-0.3cm]
\begin{enumerate}[(i).]
\item There is a deterministic matrix $A\in\R^{d\times p}$ with $d\geq p$, and i.i.d.~random vectors $Z_1,\dots,Z_n\in\R^d$, such that for each $i\in\{1,\dots,n\}$, the observation $X_i\in\R^p$ is generated as
\begin{equation}\label{eqn:Xidef}
X_i=A\ttop Z_i.
\end{equation}
\item The  random vector $Z_1$ has independent entries that satisfy $\E[Z_{1j}]=0$, $\E[Z_{1j}^2]=1$, and $\kappa:=\E[Z_{1j}^4]>1$ for all $j\in\{1,\dots,d\}$, where $\kappa$ does not depend on $n$. In addition, there is a constant $c_0$ not depending on $n$ such that $\max_{1\leq j\leq d}\|Z_{1j}\|_{\psi_2}\leq c_0$.\\[-0.2cm]
\item There are constants $\beta>1/2$ and $c_1,c_2>0$, not depending on $n$, such that for each $j\in\{1,\dots,p\}$, the singular value $\sigma_j(A)$ satisfies
\begin{equation*}
 c_1 j^{-\beta} \, \leq \, \sigma_j(A) \, \leq \, c_2j^{-\beta}.
\end{equation*}
\end{enumerate}
\end{assumption}

\paragraph{Remarks} In statistical applications, the matrix $A$ is typically taken to be the square root $\Sigma^{1/2}$, with $p=d$. However, the extra generality of a rectangular matrix is needed for the application of our work to sketching algorithms in Section~\ref{sec:sketching}. To comment on two other aspects of Assumption~\ref{A:model}, observe that it places no constraints on the relationship between $n$ and $p$, and it allows for many eigenvalues of $\Sigma$ to be repeated.

In order to state our main result, we need to precisely define the statistic $T_n^*$ that arises from bootstrap sampling. Let $(X_1^*,\dots,X_n^*)$ be drawn with replacement from $(X_1,\dots,X_n)$, and define the matrix
$$\hat\Sigma^*=\frac 1n\sum_{i=1}^n X_i^*(X_i^*)\ttop.$$
Then, the bootstrapped counterpart of $T_n$ is defined as
$$T_n^* \ = \ \sqrt n\|\widehat\Sigma^* -\widehat\Sigma\|\op.$$
The following is our main result.

\begin{theorem}\label{THM:MAIN}
Suppose that Assumption~\ref{A:model} holds. Then, there is a constant $c>0$ not depending on $n$ such that the event
\begin{equation*}
d_{\textup{K}}\big( \mathcal{L}(T_n)\, ,\, \mathcal{L}(T_n^*|X)\big)  \ \leq \  c\,n^{-\frac{\,\beta-1/2}{6\beta+4}}\,\log(n)^c
\end{equation*}
occurs with probability at least $1-\frac cn$.
\end{theorem}

\paragraph{Remarks} To explain how the difference $\beta-1/2$ arises in the rate of bootstrap approximation, we offer some informal discussion.  As preparatory notation, define the ellipsoidal boundary set 
$\mathcal{E} = \{Au \, | \, u\in\mathbb{S}^{p-1}\},$
as well as its signed version $\Theta=\mathcal{E}\times \{\pm 1\}$, whose generic element is denoted by $\theta=(v,s)$. 
With these items in place, we will consider the following empirical process indexed by $\Theta$,
\begin{equation*}
 \G_n(\theta)= \ts\frac{s}{\sqrt n} \displaystyle\sum_{i=1}^n \langle v,Z_i\rangle^2-\E[\langle v,Z_i\rangle^2],
\end{equation*}
which allows $T_n$ to be represented as
\begin{equation*}
T_n \ = \ \sup_{\theta\in\Theta}\,\G_n(\theta).
\end{equation*}
Given that the set $\Theta$ is uncountable, a standard reduction is to approximate $T_n$ with the supremum of $\G_n$ over a discrete $\e$-net for $\Theta$, where the metric is taken to be \smash{$\rho(\theta,\tilde\theta)=\|v-\tilde v\|_2+|s-\tilde s|$.} In turn, this requires us to control the discretization error, which leads to bounding the supremum of increments, denoted
$$\Delta_n(\e) \ = \ \sup_{\rho(\theta,\tilde\theta)\leq \e}|\G_n(\theta)-\G_n(\tilde\theta)|.$$
In order for the discrete approximation to succeed, the quantity $\E[\Delta_n(\e)]$ should vanish as $\e\to 0$. However, the demonstration of this property depends on the complexity of $\Theta$ through the parameter $\beta$.

We can gain some intuition for the role of $\beta$ by looking at how it affects $\E[\Delta_n(\e)]$ in a much  simpler case --- where $\G_n$ is replaced by a linear Gaussian process indexed by $\Theta$. Namely, consider the process
%\begin{equation}
$\tilde \G_n(\theta)= n^{-1/2}\sum_{i=1}^n s\langle v,\zeta_i\rangle$,
%\end{equation}
where $\zeta_1,\dots,\zeta_n$ are independent standard Gaussian vectors. In this case, if $\tilde\Delta_n(\e)$ denotes the analogue of $\Delta_n(\e)$ for $\tilde \G_n$, then the following lower bound can be shown using classical facts about Gaussian processes,
\begin{equation}\label{eqn:lowertext}
\E\big[\tilde\Delta_n(\e)\big] \ \geq \ c \, \e^{(\beta-1/2)/\beta}
\end{equation}
\citep[cf.][Proposition 2.5.1]{Talagrand:2014}.
Thus, the main point to take away here is that even in the simple case of a linear Gaussian process, the condition $\beta>1/2$ is necessary in order for the discretization error to vanish as $\e\to 0$.

Another benefit of looking at the linear Gaussian case is that the lower bound~\eqref{eqn:lowertext} provides a reference point for assessing our upper bound on the discretization error. For instance, it will follow from Proposition~\ref{prop:key} that
\begin{equation}\label{eqn:specialupper}
 \E[\Delta_n(\e)] \ \leq \ c \e^{(\beta-1/2)/\beta}\log(n).
\end{equation}
 Hence, it is notable that the dependence on $\e$ does not change in comparison to the linear Gaussian case, even though the quadratic nature of the process $\G_n$ causes it \emph{not to be sub-Gaussian} with respect to the metric $\rho$.  Moreover, it also turns out that the dependence on $\e$  even remains the same for $L_q$ norms of $\Delta_n(\e)$ when $q$ is large.

One more point of theoretical interest is that the bound~\eqref{eqn:specialupper} arises in a situation where standard chaining seems to give a slower dependence on $\e$ than a more problem-specific approach. As an example of a standard approach, one might try to show that $\G_n$ is sub-exponential with respect to $\rho$, and then appeal to an entropy integral bound such as in~\citep[][Theorem 2.2.4]{vanderVaart:Wellner:2000}. However, this ultimately leads to an upper bound scaling like $\e^{(\beta-1)/\beta}$, which would require the excessive condition \smash{$\beta>1$} (as opposed to $\beta>1/2$).  Likewise, the development of new techniques for quadratic processes akin to $\G_n$ has attracted interest in the literature, as surveyed in~\citep[][\textsection\,9.3-9.4]{Talagrand:2014}. Nevertheless, it should also be noted that existing results in this direction do not seem to be directly applicable to our analysis of the bootstrap. For instance, the abstract approaches based on Talagrand's $\gamma_1$ and $\gamma_2$ functionals lead to challenges in connection with the bootstrap, because the discrete process $\G_n^*$ induces a random metric on  $\Theta$ that does not lend itself to calculations. On the other hand, the approach taken here allows $\G_n$ and $\G_n^*$ to be treated on nearly equal footing~(cf.~Proposition~\ref{prop:keyboot}).

\section{Application to inference on covariance matrices}\label{sec:cov}
To illustrate the numerical performance of the bootstrap, this section considers two types of inference problems associated with covariance matrices: (1) error estimation for $\hat\Sigma$, and (2) construction of simultaneous confidence intervals for the eigenvalues of $\Sigma$. In particular, all of the numerical results were obtained in a situation where the leading eigenvalue $\lambda_1(\Sigma)$ has high multiplicity.

\paragraph{Simulation settings} The simulations were based on the model described in Assumption~\ref{A:model}, with $n\in\{300,500,700\}$ and $d=p=1,\!000$, giving $n<p$ in every case. Also, the matrix $A$ was constructed to be symmetric so that it can be interpreted as $A=\Sigma^{1/2}$. To specify $A$ in more detail, its singular values (equivalently eigenvalues) were chosen as
\begin{equation*}
\sigma_1(A)=\cdots =\sigma_5(A)=1 \ \ \ \ \ \ \text{ and } \ \ \ \ \ \sigma_j(A)=j^{-\beta} \ \ \text{ for } j\in\{6,\dots,p\},
\end{equation*}
with decay parameter values $\beta\in\{0.75, 1.0, 1.25\}$, and its eigenvectors were taken as the orthogonal factor from a QR decomposition of a $p\times p$ matrix with independent $N(0,1)$ entries.
Next, for each pair $(n,\beta)$, we conducted 5,000 trials in which the $n\times p$ data matrix $X=ZA$ was generated by filling $Z\in\R^{n\times p}$ with independent random variables drawn from $N(0,1)$ or a standardized $t_{20}$ distribution. Lastly, for each trial, we generated 500 bootstrap samples $T_n^*$ by sampling the rows of $X$ with replacement, as described in Section~\ref{sec:main}.

\subsection{Error estimation for $\hat\Sigma$} A natural way to formulate the problem of error estimation for $\hat\Sigma$ is in terms of the $1-\alpha$ quantile of $T_n$, denoted by $q_{1-\alpha}$. By definition, this quantity gives the tightest bound of the form
\begin{equation*}
\|\hat\Sigma-\Sigma\|\op \ \leq \ \frac{q_{1-\alpha}}{\sqrt{n}}
\end{equation*}
that holds with probability at least $1-\alpha$. Likewise, if we let $\hat q_{1-\alpha}$ denote the empirical $(1-\alpha)$-quantile of the bootstrap samples $T_n^*$, then we may regard $\hat q_{1-\alpha}/\sqrt{n}$ as an error estimate for $\hat\Sigma$. 

Alternatively, the estimate $\hat q_{1-\alpha}$ can be viewed as specifying an approximate  $(1-\alpha)$-confidence region for $\Sigma$. That is, if we let $\mathsf{B}_{\text{op}}(r;\hat\Sigma)\subset\R^{p\times p}$ denote the operator-norm ball of radius $r>0$ centered at $\hat\Sigma$, then $q_{1-\alpha}$ is the smallest value of $r$ such that
$$\P\Big(\Sigma\in{\mathsf{B}}\op(r;\hat\Sigma)\Big) \ \geq \ 1-\alpha.$$
Hence, the ideal confidence region may be approximated with ${\mathsf{B}}\op\big(\ts\frac{\hat q_{1-\alpha}}{\sqrt n};\hat\Sigma\big)$.

\begin{table}[!b]
	\caption{Observed coverage probabilities for $\hat q_{1-\alpha}$, with $\alpha=0.1$ and $p=1,000$.
}
\label{tab:summary_coverage}
	%\vspace{-0.5cm}
\begin{subtable}{0.49\textwidth}
		\centering
		\scalebox{0.85}{	
		\begin{tabular}{c  cccccc} %\toprule
			\multirow{2}{*}{\makecell{decay \\ param. $\beta$}}	 & \multicolumn{3}{c}{sample size $n$} \\
			&    300 &  500 & 700 \\
			\midrule%\midrule
			$0.75$   & 92.83\% & 92.26\% & 91.96\% \\
			$1.00$   & 92.66\% & 91.70\% & 91.23\% \\
			$1.25$   & 92.43\% & 91.53\% & 91.16\%  \\
 			\bottomrule 
	\end{tabular}}
	\caption{$N(0,1)$ distribution}
	\end{subtable}%
~~~~~
	\begin{subtable}{0.49\textwidth}
		\centering
		\scalebox{0.85}{	
		\begin{tabular}{c  cccccc}% \toprule
		\multirow{2}{*}{\makecell{decay \\ param. $\beta$}}	 & \multicolumn{3}{c}{sample size $n$} \\
		&   300 &  500 & 700 \\
		\midrule%\midrule
		$0.75$      & 92.90\% & 91.96\% & 91.93\%  \\
		$1.00$ 	    & 92.53\% & 91.76\% & 91.63\%   \\
		$1.25$      & 92.50\% & 91.70\% & 91.46\%   \\
		\bottomrule 
		\end{tabular}}
		\caption{standardized $t_{20}$ distribution}
	\end{subtable}%
\end{table}
To demonstrate the performance of $\hat q_{1-\alpha}$, its observed coverage probabilities have been listed in Table~\ref{tab:summary_coverage}. (Note that these probabilities can be interpreted either with respect to the coverage of the error bound or the confidence region.) The table is organized based on the simulation settings described earlier, and in all cases, the nominal coverage level was set to 90\%. Looking at the results, we see that for sufficiently large sample sizes, the observed coverage comes nearly within 1\% of the desired level. Another important feature of $\hat q_{1-\alpha}$ is that its errors occur very reliably in the conservative direction, with the observed coverage never falling below the nominal level. With regard to the parameter $\beta$, we see the intuitive effect that more spectrum decay yields better coverage, but that this improvement occurs very gradually  as a function of $\beta$, which is understandable in light of Theorem~\ref{THM:MAIN}.

\begin{table}[!b]
	\caption{Observed simultaneous coverage probabilities for $\hat{\mathcal{I}}_1,\dots,\hat{\mathcal{I}}_p$, with  $\alpha=0.1$ and $p=1,000$.
	}
	\label{tab:summary_eigs}
	%\vspace{-0.5cm}
	\begin{subtable}{0.49\textwidth}
		\centering
		\scalebox{0.85}{	
		\begin{tabular}{c  cccccc} %\toprule
			\multirow{2}{*}{\makecell{decay \\ param. $\beta$}}	 & \multicolumn{3}{c}{sample size $n$} \\
			&  300 &  500 & 700 \\
			\midrule%\midrule
			$0.75$   & 94.46\% & 94.26\% & 93.26\%\\
			$1.00$   & 93.13\% & 92.06\% & 91.53\% \\
			$1.25$   & 92.63\% & 91.73\% & 91.23\% \\
		
 			\bottomrule 
	\end{tabular}}
	%\centering
	\caption{$N(0,1)$ distribution }
	\end{subtable}%
~~~~
	\begin{subtable}{0.49\textwidth}
		\centering
		\centering
		\scalebox{0.85}{	
		\begin{tabular}{c  cccccc}% \toprule
		\multirow{2}{*}{\makecell{decay \\ param. $\beta$}}	 & \multicolumn{3}{c}{sample size $n$} \\
		&  300 &  500 & 700 \\
		\midrule%\midrule
		$0.75$  & 94.03\% & 93.87\% & 93.76\% \\
		$1.00$ 	& 92.90\% & 91.66\% & 91.46\%   \\
		$1.25$  & 92.56\% & 91.40\% & 91.38\%   \\	
		\bottomrule 
		\end{tabular}}
		\caption{standardized $t_{20}$ distribution}
	\end{subtable}%
\end{table}

\subsection{Simultaneous confidence intervals}\label{sec:simult} Here, we consider the problem of approximating a collection of random intervals $\mathcal{I}_1,\dots,\mathcal{I}_p$ that satisfy
\begin{equation}\label{eqn:simult}
\P\bigg(\bigcap_{j=1}^p \big\{\lambda_j(\Sigma)\in\mathcal{I}_j\big\}\bigg) \ \geq \ 1-\alpha.
\end{equation}
Our approach is based on Weyl's inequality, which ensures that the condition
\begin{equation*}
|\lambda_j(\hat\Sigma)-\lambda_j(\Sigma)| \ \leq \ \|\hat\Sigma-\Sigma\|\op 
\end{equation*}
holds simultaneously for all $j\in\{1,\dots,p\}$, with probability 1. To proceed, let $q_{1-\alpha}$ again denote the $(1-\alpha)$-quantile of $T_n$, and let $\mathcal{I}_j=[\lambda_j(\hat\Sigma)\pm q_{1-\alpha}/\sqrt{n}]$. Then, Weyl's inequality implies that the condition~\eqref{eqn:simult} must hold. In turn, we may use the bootstrap estimate $\hat q_{1-\alpha}$ to form the approximate intervals defined by $\hat{\mathcal{I}}_j=[\lambda_j(\hat\Sigma)\pm \hat q_{1-\alpha}/\sqrt{n}]$. 

As a way to gain robustness against the effects of eigenvalue multiplicity, the papers~\citep{Hall:2006,Hall:2009} also considered an approach of this type --- but instead using the Frobenius norm, which can lead to potentially much wider intervals than the operator norm. In the latter paper, a further refinement of this approach was developed, and in fact, it would be possible to be combine it with our operator-norm based intervals, but we do not pursue this here for the sake of brevity.

The simulation results for the intervals $\hat{\mathcal{I}}_1,\dots,\hat{\mathcal{I}}_p$ are given in Table~\ref{tab:summary_eigs}, based on the previous settings where $\lambda_1(\Sigma)$ has multiplicity 5. (The entries of the table are the observed simultaneous coverage probabilities for a nominal level of 90\%.) Although the intervals are somewhat conservative due to Weyl's inequality, they are still close enough to the nominal level to be of practical interest, especially for larger values of $\beta$. Also, to put matters into context, it is important to note that a naive application of the bootstrap to the individual sample eigenvalues is known to work poorly in the presence of high multiplicity.  Hence, the user may be willing to tolerate a bit of conservatism in order to avoid the harms of closely spaced population eigenvalues.

\section{Application to randomized numerical linear algebra}\label{sec:sketching} 
Over the past decade, RandNLA has become the focus of intense activity in many fields related to large-scale computation~\citep{mahoney2011randomized,Halko:2011,woodruff2014sketching,Kannan:2017,RandNLA_PCMIchapter_chapter}. Broadly speaking, this new direction of research 
has stemmed from the principle that randomization is a very general mechanism for scaling up algorithms. However, in exchange for scalability, randomized sketching algorithms are typically less accurate than their deterministic predecessors. 
 Therefore, in order to use sketching reliably, it is crucial to verify that the algorithmic error is small, which motivates new applications of the bootstrap beyond its traditional domains.

The purpose of this section is to illustrate how the bootstrap can be applied to estimate operator-norm error for randomized matrix multiplication, which has been a prominent topic in the RandNLA literature~\citep[e.g.,][]{Drineas:2001,Drineas:2006,Magen:2011,Pagh:2013,Ipsen:2015,Cohen:2016,Gupta:2018}. A related study of the bootstrap for this application can also be found in~\citep{Lopes:JMLR}, which differs from the current work insofar as it deals exclusively with the entrywise $\ell_{\infty}$-norm and does not focus the role of spectrum decay.

To proceed, we will first provide a brief review of the algorithmic setting (Section~\ref{sec:basicsketch}), followed by an efficient implementation of the bootstrap using an extrapolation technique~(Section~\ref{sec:bootsketch}). In turn, we will present numerical results for synthetic matrices (Section~\ref{sec:syntheticsketch}), as well as an example concerning spatial modes of temperature variation (Section~\ref{sec:seasketch}).

\subsection{Rudiments of sketching}\label{sec:basicsketch}
Consider a situation involving a very large deterministic matrix $A\in\R^{d\times p}$ with $p\ll d$, where the product $A\ttop A$ is too expensive to compute to high precision. 
For instance, this often occurs when $A$ must be stored on disk because it exceeds the constraints of fast memory. Alternatively, even when memory is not a bottleneck, matrix products can become too expensive if they must be computed frequently as a subroutine of a larger pipeline.

The basic idea of sketching is to work with a shorter version of $A$, referred to as a ``sketch of $A$'', and denoted as $\tilde A\in \R^{n\times p}$ where $n\ll d$. This matrix is defined as
\begin{equation*}
\tilde A = SA,
\end{equation*}
where $S\in\R^{n\times d}$ is a random ``sketching matrix'' that is generated by the user. In particular, the user must choose the ``sketch size'' $n$. Intuitively, the matrix $S$ is intended to shorten $A$ in a way that retains most of the information, so that the inexpensive product $\tilde A\ttop \tilde A$ will provide a good approximation to~$A\ttop A$.  

\paragraph{The sketching matrix} Typically, the action of $S$ upon $A$ is interpreted in either of two ways:  randomly projecting columns from $\R^d$ into $\R^n$, or discretely sampling $n$ among $d$ rows.
In addition, the matrix $S$ is commonly generated by the user so that its rows are i.i.d., and that it satisfies $\E[S\ttop S]=I_n$, which implies that $\tilde A\ttop \tilde A$ is unbiased with respect to  $A\ttop A$. At a high level, these basic properties are sufficient to understand  all of our work below, but numerous types of sketching matrices have been studied in the literature.
For instance, two of the most well-known are the \emph{Gaussian random projection} and \emph{uniform row sampling} types, where the former has i.i.d.~rows drawn from $N(0,\frac{1}{n}I_d)$, and the latter has i.i.d.~rows drawn uniformly from $\{\sqrt{d/n}\,e_1,\dots,\sqrt{d/n}\,e_d\}$. More elaborate examples may be found in the references above.

\paragraph{Cost versus accuracy} Whenever sketching is implemented, the choice of the  sketch size $n$ plays a pivotal role in a tradeoff between computational cost and accuracy. To see this, note that on one hand, the cost to compute $\tilde A\ttop\tilde A$ is generally proportional to $n$, with the number of operations being $\mathcal{O}(np^2)$. On the other hand, the operator-norm error of $\tilde A\ttop \tilde A$ tends to decrease stochastically like $1/\sqrt{n}$, because the difference $\tilde A\ttop \tilde A-A\ttop A$ can be expressed as a sample average of $n$ centered random matrices (i.e.,~in the same way as $\hat\Sigma-\Sigma$).\\

\begin{figure}[H]
	\centering
	\DeclareGraphicsExtensions{.pdf}
	\begin{overpic}[width=0.70\textwidth]{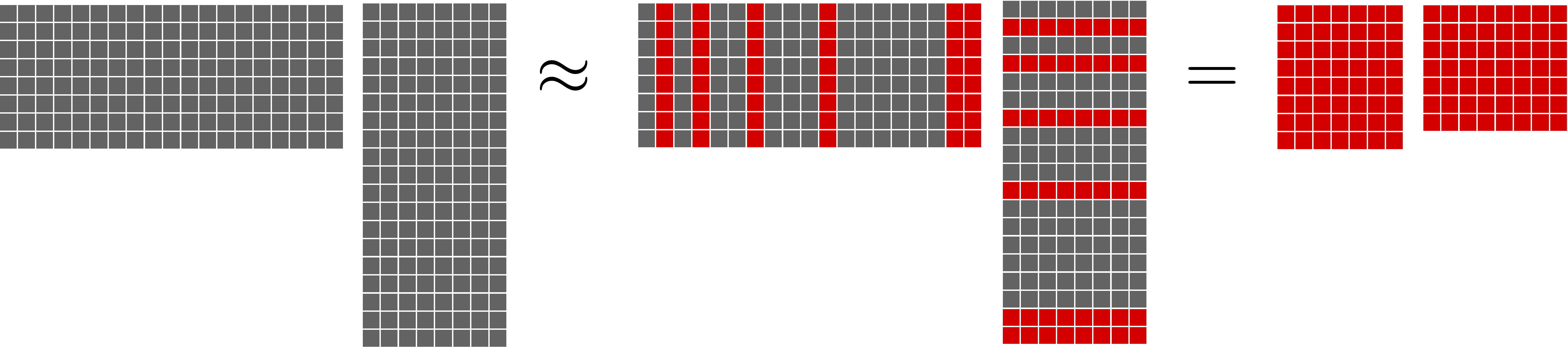} 
		\put(10,23){\color{black}{$\A^\top$}} 
		\put(26,23){\color{black}{$\A$}} 
		\put(49,23){\color{black}{$({SA})^\top$}} 
		\put(66,23){\color{black}{${SA}$}} 
		\put(84,23){\color{black}{${\tilde{A}}^\top$}} 
		\put(94,23){\color{black}{${\tilde{A}}$}} 
	\end{overpic}
	%\vspace{+.10in}	
	\caption{Randomized matrix multiplication with row sampling.}
	\label{fig:sketching}
\end{figure}
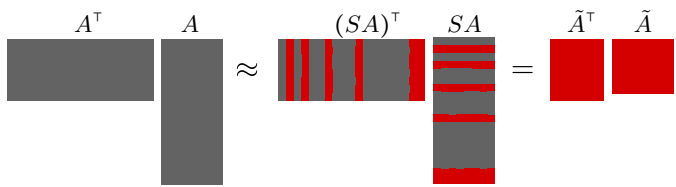

\paragraph{The need for error estimation} Although the choice of the sketch size $n$ has clear importance, this choice also involves practical difficulties that expose a major gap between the theory and practice of sketching. Specifically, these difficulties arise because the value of the sketching error $\|\tilde A\ttop \tilde A -A\ttop A\|\op$ is  \emph{unknown in practice}, as it depends on the unknown product $A\ttop A$. Hence, it is hard for the user to know if any given choice of $n$ will achieve a desired level of accuracy. 

As a way to handle this dilemma, one option is to consult the RandNLA literature on theoretical error bounds for $\|\tilde A\ttop\tilde A-A\ttop A\|\op$, as surveyed in the references above. However, much like in the setting of covariance estimation, these results usually only provide qualitative guidance, and they rarely offer an explicit numerical bound. Most often, this occurs because of unspecified theoretical constants, but there is also a second key limitation: Theoretical error bounds are generally formulated to hold in a worst-case sense, and so they often fail to account for special structure. Due to these issues, we propose instead to directly estimate the error via a computationally efficient bootstrap method. This has the twofold benefit of providing a numerical bound and adapting automatically to the structure of the problem at hand.

\paragraph{Comparison of sketching and covariance estimation} To clarify the relationship
between the sketching error $\|\tilde A\ttop \tilde A-A\ttop A\|\op$ and the covariance estimation error $\|\hat\Sigma-\Sigma\|\op$, let $A$, $\Sigma$, and $\hat\Sigma$ be understood as in the context of the model~\ref{A:model} with $Z\in\R^{n\times d}$ having rows $Z_1,\dots,Z_n$, and let 
$$S=\ts\frac{1}{\sqrt n}\, Z.$$
Under these conditions, the matrix $S$ has the desired properties of a sketching matrix mentioned earlier, and furthermore $\tilde A\ttop \tilde A-A\ttop A= \hat\Sigma-\Sigma$.
However, it is worth highlighting that this formal similarity conceals some operational differences. For instance, the matrices $Z$ and $A$ are unobservable to the user in covariance estimation, whereas the user does have access to $S$ and $A$ in sketching. Secondly, in covariance estimation, the user often does not have the option to increase $n$, but in sketching, it is possible to construct a rough initial sketch of $A$ for inspection, and then take a second sketch to improve performance. Later on, we will show how this second point has an important link with our error estimation method, because it will enable the user to \emph{dynamically predict} the total sketch size needed to reach a given level of accuracy.

\subsection{Error estimation with an extrapolated bootstrap}\label{sec:bootsketch}
The intuition for applying the bootstrap to sketching comes from thinking of the matrix $\tilde A$ as a ``dataset'' whose rows are ``observations''. In particular, this interpretation is supported by the fact that many types of sketching matrices $S$ cause the rows of $\tilde A$ to be i.i.d. Therefore, we may expect that sampling from the rows of $\tilde A$ with replacement will faithfully mimic the process that generated $\tilde A$.

To fix some notation for describing the bootstrap method, let ${\tt{q}}_{1-\alpha}$ denote the $(1-\alpha)$-quantile of the sketching error variable $\|\tilde A\ttop \tilde A-A\ttop A\|\op$, which is the minimal value such that the event
$$\|\tilde A\ttop \tilde A-A\ttop A\|\op \ \leq \ {\tt{q}}_{1-\alpha}$$
 holds with probability at least $1-\alpha$. Our main goal is to construct an estimate $\hat{{\tt{q}}}_{1-\alpha}$ using only the sketch $\tilde A$ as a source of information. Below, we state a basic version of the bootstrap method in Algorithm~\ref{alg:bootstrap}, which will later be accelerated via an extrapolation technique in Section~\ref{sec:extrap}. 
\begin{alg} \label{alg:bootstrap}
	\normalfont
	{(Bootstrap estimate of sketching error).}\\[-0.3cm]
	\hrule
	\vspace{0.1cm}
\noindent {\bf Input}: The number of bootstrap samples $B$, and the sketch $\tilde{\A}\in\R^{n\times p}$.\\[0.2cm]
	{\bf For } $b=1,\dots,B$\; {\bf do}
	\begin{enumerate}
		\item Form $\tilde{\A}^*\in\R^{n\times p}$ by drawing $n$ rows from $\tilde{A}$ with replacement.\\[-0.2cm]
		\item
		Compute the bootstrap sample $\ve^{*}_{b}:=\big\|(\tilde{\A}^*)^\top (\tilde{\A}^*)-\tilde{\A}^\top \tilde{\A}\big\|\op$.
	\end{enumerate}
	{\bf Return:}
	$\hat{{\tt{q}}}_{1-\alpha} \longleftarrow$
	the $(1-\alpha)$-quantile of the values $\ve^{*}_{1},\dots,\ve^{*}_{B}$.
\end{alg}
\vspace{-0.2cm}
\hrule
\vspace{0.2cm}
\paragraph{Remark} Given that the construction of $\tilde A$ is fully controlled by the user, one might ask why bootstrapping is preferable to carrying out many repetitions of the actual sketching process. The answer comes down to the fact that constructing $\tilde A$ requires a computation involving the full matrix $A$, which often incurs high communication costs. In fact, this issue is one of the primary motivations for the whole subject of RandNLA, which is usually deals with situations where it is only feasible to access $A$ at most a handful of times. In contrast to the task of constructing $\tilde A$, Algorithm~\ref{alg:bootstrap} only requires inexpensive access to the much smaller matrix $\tilde A$, and it requires no access to $A$ whatsoever.

\subsubsection{Extrapolation}\label{sec:extrap}
Because the user has the option to increase the sketch size $n$ by performing an extra round of sketching, it becomes possible to accelerate the bootstrap with an extrapolation technique that is often not applicable in covariance estimation. To develop the idea, we should first recall that the fluctuations of $\|\tilde A\ttop \tilde A-A\ttop A\|\op$ tend to scale like $1/\sqrt n$ as a function of $n$, because the difference $\tilde A\ttop \tilde A-A\ttop A$ can be written as a centered sample average of $n$ random matrices. Therefore, if we view the sketching error quantile as a function of $n$, say ${\tt{q}}_{1-\alpha}={\tt{q}}_{1-\alpha}(n)$, then we may expect the following approximate relationship between a small ``initial'' sketch size $n_0$, and a larger ``final'' sketch size $n_1$,
\begin{equation}\label{eqn:extapprox}
{\tt{q}}_{1-\alpha}(n_1) \ \approx \ \sqrt{\ts\frac{n_0}{n_1}}\,\,{\tt{q}}_{1-\alpha}(n_0).
\end{equation}
The significance of this approximation is that ${\tt{q}}_{1-\alpha}(n_0)$ is computationally much easier to estimate than ${\tt{q}}_{1-\alpha}(n_1)$, since the former involves bootstrapping a matrix of size $n_0\times p$, rather than $n_1\times p$. (More general background on the connections between extrapolation and resampling methods can be found in~\citep{Bickel:1988,Bertail:1997,Bertail:2001,Bickel:2002,Lopes:AOSRF}, among others.)

Based on the heuristic approximation \eqref{eqn:extapprox}, we can obtain an inexpensive estimate of ${\tt{q}}_{1-\alpha}(n_1)$ for any $n_1>n_0$ by using
\begin{equation}\label{eqn:ext}
\hat{\tt{q}}_{1-\alpha}^{\text{\ ext}}(n_1) \ := \ \sqrt{\ts\frac{n_0}{n_1}}\,\,\hat{{\tt{q}}}_{1-\alpha}(n_0),
\end{equation}
where $\hat{{\tt{q}}}_{1-\alpha}(n_0)$ is obtained from Algorithm~\ref{alg:bootstrap}. More concretely, if the user has the ultimate intention of achieving \smash{$\|\tilde A\ttop \tilde A-A\ttop A\|\op\leq \e_{\text{tol}}$} for some tolerance $\e_{\text{tol}}$, then extrapolation may be applied in the following way: First, the user should check the condition $\hat{{\tt{q}}}_{1-\alpha}(n_0)\leq \e_{\text{tol}}$ to see if $n_0$ is already large enough. Second, if $n_0$ is too small, then the rule~\eqref{eqn:ext} instructs the user to obtain a final sketch size $n_1$ satisfying $\hat{\tt{q}}_{1-\alpha}^{\text{\ ext}}(n_1) \leq \e_{\text{tol}}$, which is equivalent to
\begin{equation*}
n_1 \ \geq \ \ts\frac{n_0}{\e_{\text{tol}}^2}\, \hat{\tt{q}}_{1-\alpha}(n_0)^2.
\end{equation*}
Furthermore, our numerical results will demonstrate that this simple technique remains highly effective even when $n_1$ is much larger than $n_0$, such as by an \emph{order of magnitude}~(cf.~Sections~\ref{sec:syntheticsketch} and~\ref{sec:seasketch}).

\subsubsection{Assessment of cost} Since the overall purpose of sketching is to reduce computation, it is important to explain why the added cost of the bootstrap is manageable. In particular, the added cost should not be much higher than the cost of sketching itself.
As a simple point of reference, the cost to construct $\tilde A$ and then compute $\tilde A\ttop \tilde A$ with most state-of-the-art sketching algorithms is at least $\bunderline{C}_{\text{sketch}}=\Omega(dp + n_1p^2)$, where $n_1$ refers to the ``final'' sketch size described above.
Next, to assess the cost of the bootstrap, we can take advantage of a small initial sketch size $n_0$ by using extrapolation, as well as the fact that the bootstrap samples can be trivially computed in parallel, with say $m$ processors. When these basic factors are taken into account, the cost of the bootstrap turns out to be at most $\widebar{C}_{\text{\,boot}}=\mathcal{O}(Bn_0p^2/m)$. 

From this discussion of cost, perhaps the most essential point to emphasize is that $\bunderline{C}_{\text{sketch}}$ grows linearly with $d$, whereas  $\widebar{C}_{\text{\,boot}}$ is \emph{independent of d}. Indeed, this of great importance for scalability, because randomized matrix multiplication is of primary interest in situations where $d$ is extremely large. Beyond this high-level observation, we can also take a more detailed look to see that the condition $\widebar{C}_{\text{\,boot}}=\mathcal{O}(\bunderline{C}_{\text{sketch}})$ occurs when $B \ = \ \mathcal{O}((\ts\frac{n_1}{n_0}+\frac{d}{pn_0})m)$. Furthermore, such a condition on $B$ can be considered realistic in light of our experiments, since the modest choice of $B=50$ is shown to yield  good results.

\subsection{Numerical results for synthetic matrices}\label{sec:syntheticsketch}
We now demonstrate the performance of the bootstrap estimate $\hat{\tt{q}}_{1-\alpha}$ in a range of conditions, both with and without extrapolation. Most notably, the numerical results for extrapolation are quite encouraging.

\paragraph{Simulation settings} The choices for the matrix $A\in\R^{d\times p}$ were developed in analogy with those in Section~\ref{sec:cov}, except that in this context, the matrix is very tall with $d=10,000$ and $p=1,000$. If we let $\A = U D V^\top$ denote the singular value decomposition, then the singular vectors were specified by taking \smash{$U\in\R^{d\times p}$} and $V\in\R^{p\times p}$ to be orthonormal factors from QR decompositions of matrices filled with independent $N(0,1)$ entries. In addition, the singular values were chosen as $\sigma_1(A)=\cdots=\sigma_5(A)=1$ and $\sigma_j(A)=j^{-\beta}$ for $j\in\{6,\dots,p\}$, with decay parameters $\beta\in\{0.75, 1.0, 1.25\}$. In particular, these values were chosen in order to show that the bootstrap can work even when there are no gaps among the leading singular values.

\begin{figure}[H]
	
	\centering
	\begin{subfigure}{1\textwidth}	
		\centering
		\DeclareGraphicsExtensions{.pdf}
		\begin{overpic}[width=0.31\textwidth]{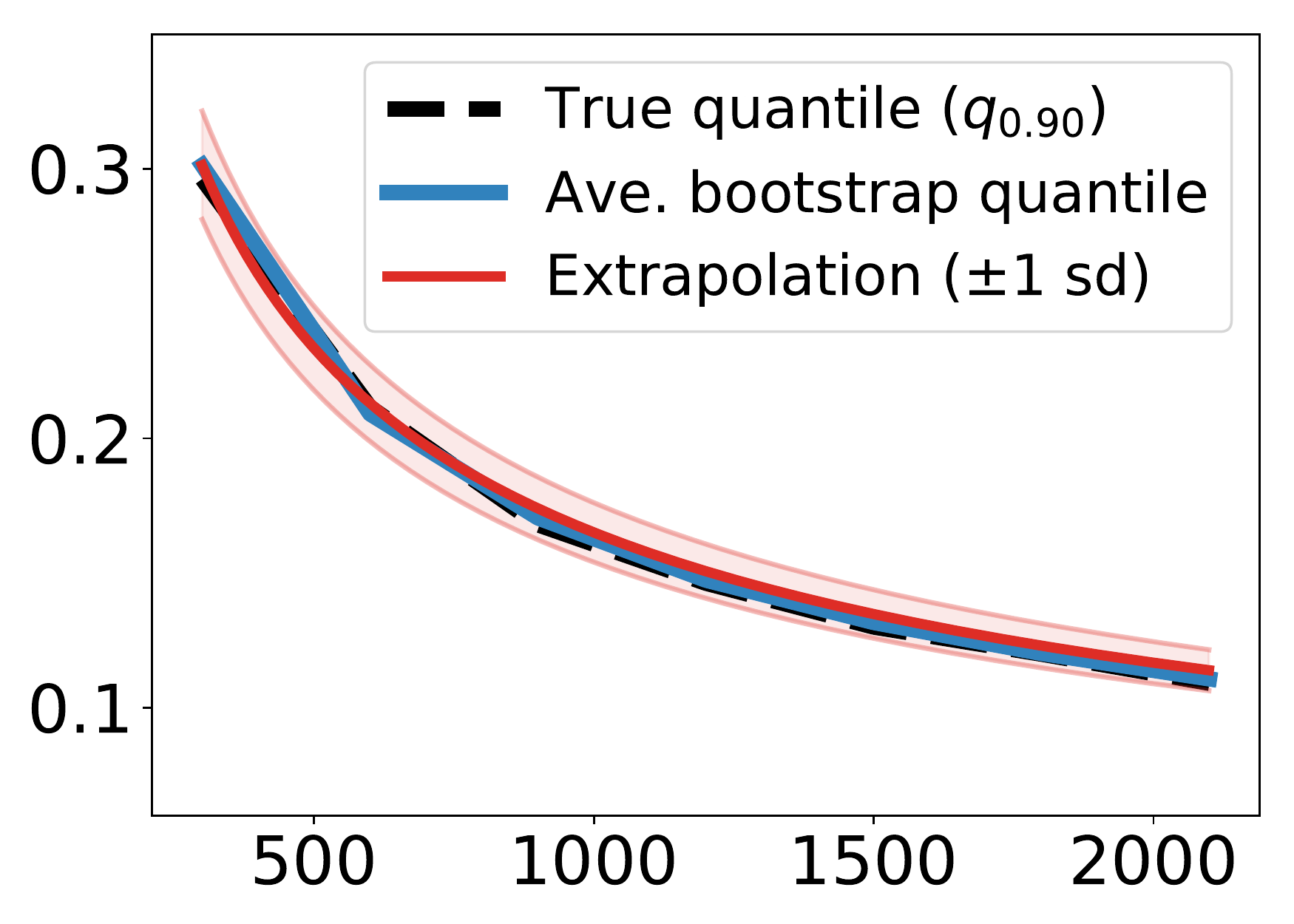} 
			\put(40,75){\color{black}{\scriptsize $\beta=0.75$}} 
			%\put(45,-1){\color{black}{\small sketch size}}   
			\put(-7,11){\rotatebox{90}{\small op.~norm error}}
		\end{overpic}\hspace*{-0.2cm}	
		~
		\begin{overpic}[width=0.31\textwidth]{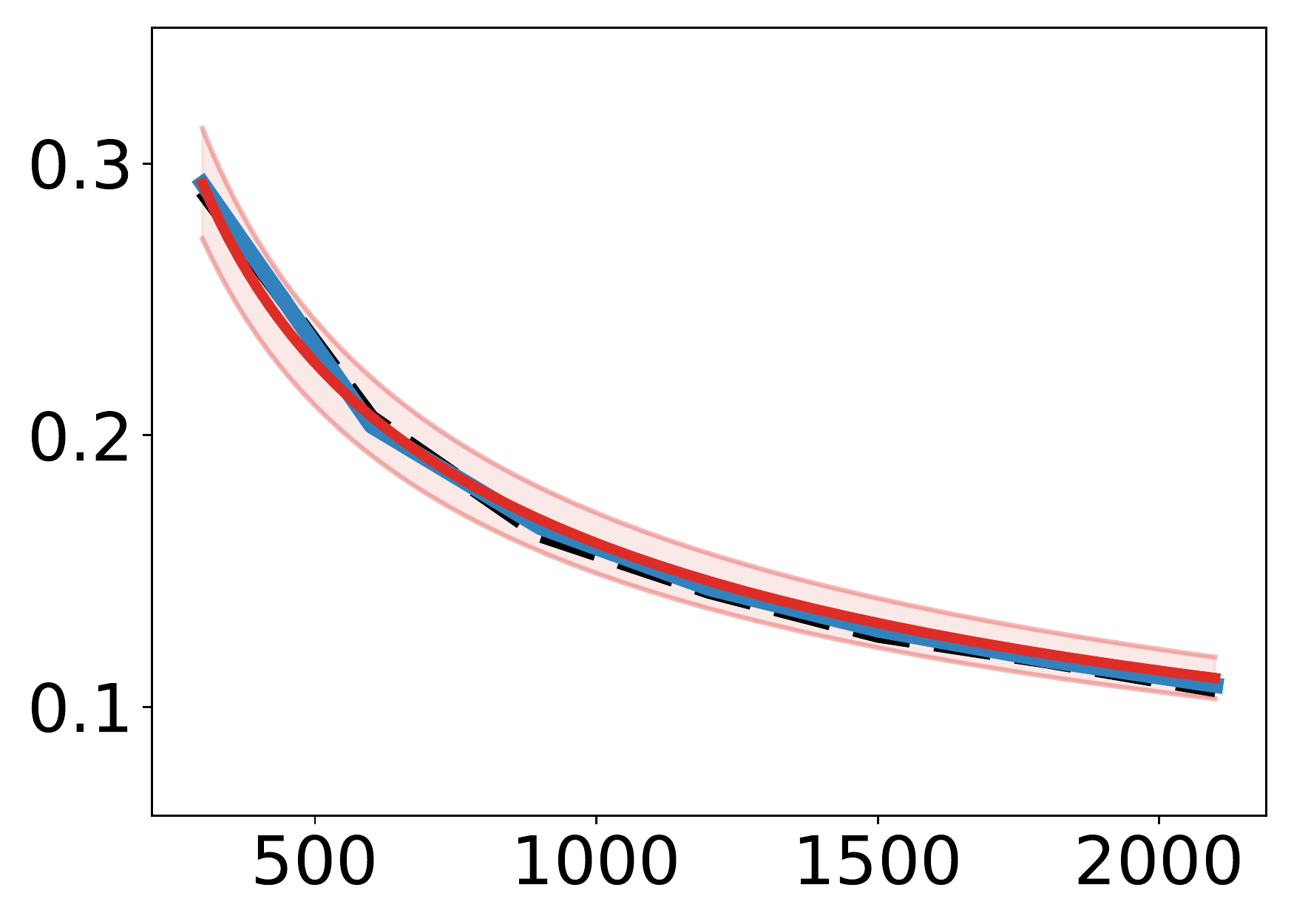} 
			\put(40,75){\color{black}{\scriptsize $\beta=1.0$}} 
			%\put(45,-1){\color{black}{\small sketch size}}   
			%\put(-2,15){\rotatebox{90}{\small $\ell_2$ norm error}}
		\end{overpic}\hspace*{-0.2cm}		
		~
		\begin{overpic}[width=0.31\textwidth]{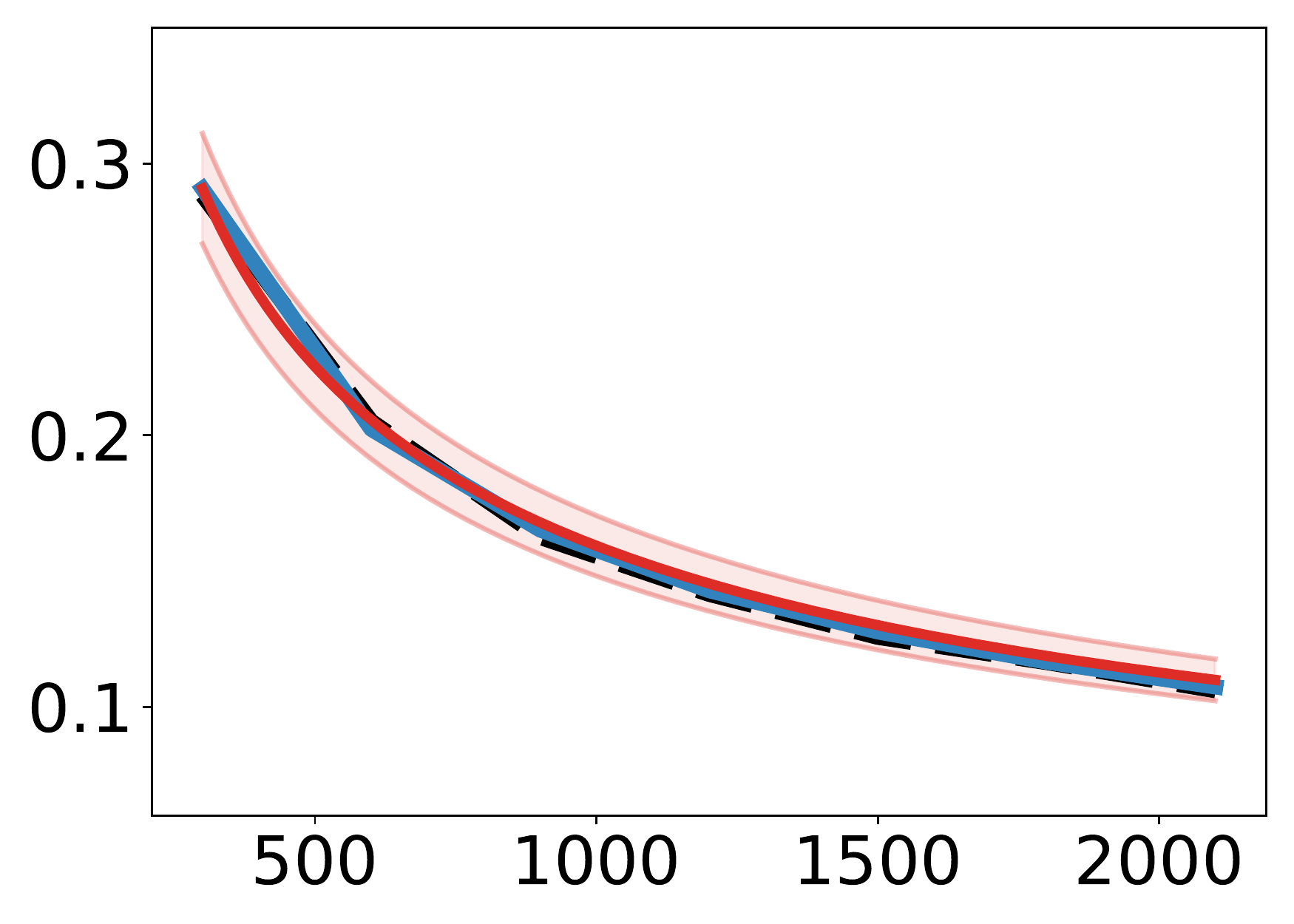} 
			\put(40,75){\color{black}{\scriptsize $\beta=1.25$}} 
		\end{overpic}
		\caption{Sketching with Gaussian random projections.} \label{fig:results_mat_sketching_gaussian}
	\end{subfigure}	
	
	\begin{subfigure}{1\textwidth}	
		\centering
		\DeclareGraphicsExtensions{.pdf}
		\begin{overpic}[width=0.31\textwidth]{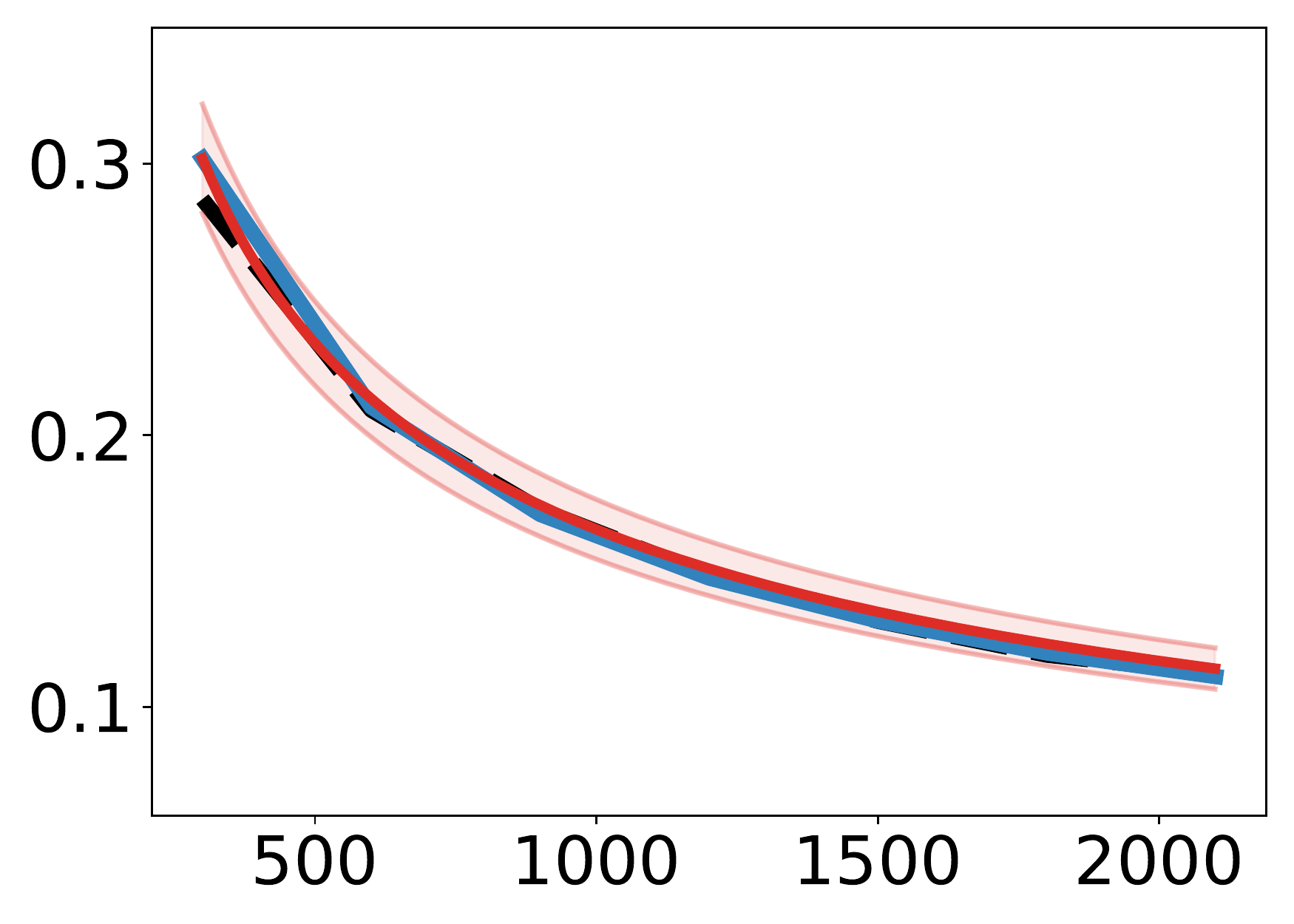} 
			\put(35,-8){\color{black}{\small sketch size $n$}}   
			\put(-7,11){\rotatebox{90}{\small op.~norm error}}
		\end{overpic}\hspace*{-0.2cm}
		~
		\begin{overpic}[width=0.31\textwidth]{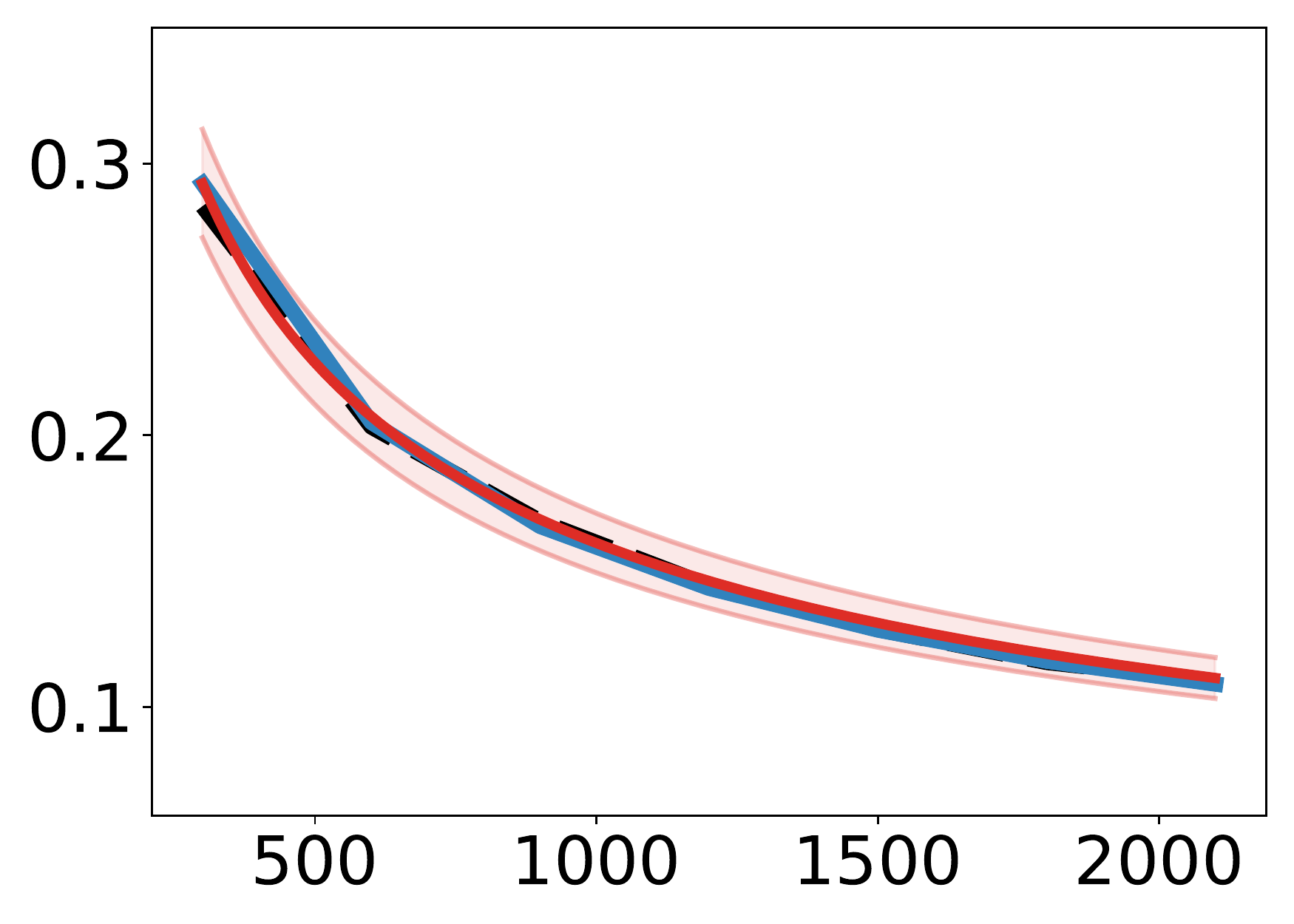} 
			\put(35,-8){\color{black}{\small sketch size $n$}}   
		\end{overpic}\hspace*{-0.2cm}
		~
		\begin{overpic}[width=0.31\textwidth]{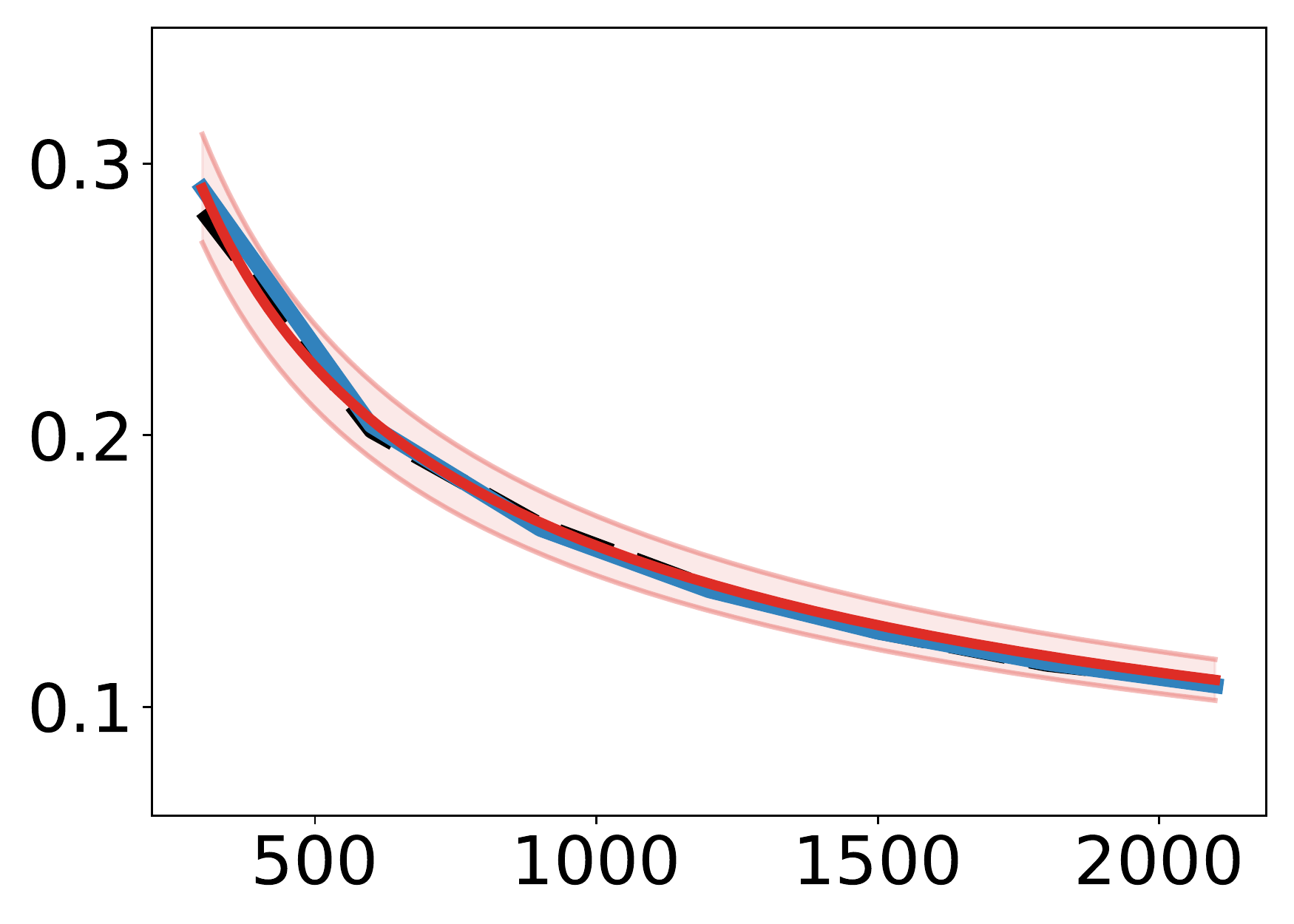} 
			\put(35,-8){\color{black}{\small sketch size $n$}}   
		\end{overpic}
		\vspace{+0.3cm}		
		\caption{Sketching with uniform row sampling.}	\label{fig:results_mat_sketching_uniform}
	\end{subfigure}
		\vspace{-0.3cm}
	%\vspace{+.10in}	
	\caption{Bootstrap estimates for the 90\% quantile of the error $\|\tilde A\ttop \tilde A-A\ttop A\|\op$.}
	\label{fig:results_mat_sketching}
\end{figure}

\paragraph{Design of simulations} The design of the simulations can be understood in terms of Figure~\ref{fig:results_mat_sketching}. For each value of $\beta$, and sketch size $n\in\{300,\dots,2,\!100\}$, we performed 1,000 trials of sketching to compute independent copies $\tilde A\in\R^{n\times p}$ using two different types of sketching matrices: Gaussian random projection, and uniform row sampling, as defined in Section~\ref{sec:basicsketch}. In turn, the actual values of $\|\tilde A\ttop \tilde A-A\ttop A\|\op$ in these trials yielded a high quality approximation to the true 90\% quantile ${\tt{q}}_{0.9}={\tt{q}}_{0.9}(n)$, plotted as a function of $n$ with the black dashed line.

With regard to Algorithm~\ref{alg:bootstrap}, it was applied during each trial to compute  $\hat{\tt{q}}_{1-\alpha}$ using $B=50$ bootstrap samples. The average of these estimates is plotted as a function of $n$ with the solid blue line. In addition, the performance of the extrapolation rule~\eqref{eqn:ext} was studied by applying it to each estimate $\hat{\tt{q}}_{0.9}(n_0)$ computed at $n_0=300$. The average of the extrapolated curves is plotted in solid red, with the pink envelope signifying $\pm1$ standard deviation.

%\vspace{0.7cm}

\vspace{-0.6cm}
\paragraph{Comments on results} Figure~\ref{fig:results_mat_sketching} shows that on average, the bootstrap estimates are nearly equal to the true quantile over the entire range of sketch sizes $n\in\{300,\dots,2,\!100\}$, both with and without extrapolation. Indeed, the performance of the extrapolated estimate is especially striking, because it shows that bootstrapping a rough initial sketch $\tilde A$ of size $300\times1,000$ can be used to accurately predict the error of a much larger sketch of size $2,\!100\times 1,\!000$. To put this into context, we should also remember that the original matrix $A$ is of size $10,\!000\times 1,\!000$, and hence the initial sketch is able to provide quite a bit of information about the sketching task for a small computational price. Moreover, the fact that the extrapolation works up to the larger sketch size of $2,\!100$ means that a 7-fold speedup can be obtained in comparison to naively applying Algorithm~\ref{alg:bootstrap} to the larger sketch. (In fact, the plots seems to suggest that the extrapolation would remain accurate for sketch sizes beyond 2,\!100, and that even larger speedups are attainable.) Lastly, it is worth noting that even though a small choice $B=50$ bootstrap samples was used, the standard deviation of the extrapolated estimate is rather well-behaved, as indicated by the pink envelope.

%\vspace{-0.2cm}
\subsection{Sea surface temperature measurements}\label{sec:seasketch}
Large-scale dynamical systems are ubiquitous in the physical sciences, and advances in technology for measuring these systems have led to rapidly increasing volumes of data. Consequently, it is often too costly to apply standard tools of exploratory data analysis in a direct manner, and there has been growing interest to use sketching as a data-reduction strategy that preserves the essential information~\citep[e.g.][among others]{Brunton:2015,Erichson:2017,ribeiro2019randomized,bai2019randomized,saibaba2019randomized,Bjarkason:2019,Tropp:2019}.

\begin{figure}[H]
	\centering
	\begin{subfigure}{0.45\textwidth}	
		\centering
		\DeclareGraphicsExtensions{.png}
		\begin{overpic}[width=.99\textwidth,trim=0.0cm 0.0cm 0.0cm 0.0cm, clip]{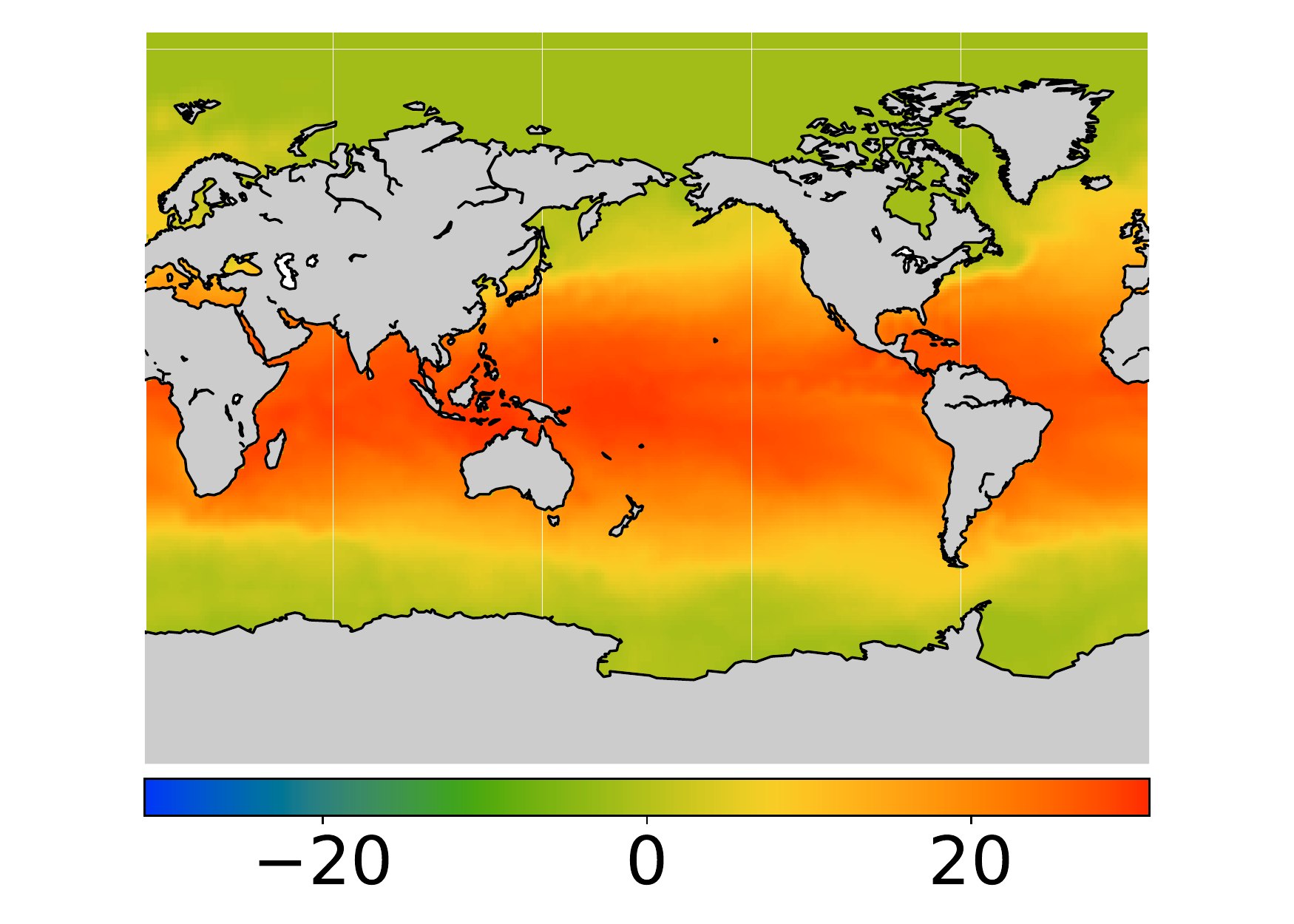} 
			\linethickness{1.0pt}
            \put(33,32){\line(1,0){51}}
            \put(84,59.5){\line(0,-1){28}}
            \put(33,59){\line(1,0){51}}
            \put(33,59.5){\line(0,-1){28}}
			\put(20,-4){\color{black}{ temperature fluctuations}}   
		\end{overpic}\vspace*{+0.2cm}
		\caption{ENSO region}\label{fig:sst_example}
	\end{subfigure}		
	~~~
		\begin{subfigure}{0.45\textwidth}	
		\centering
		\DeclareGraphicsExtensions{.png}
		\begin{overpic}[width=.99\textwidth, trim=0.0cm 0.0cm 0.0cm 0.0cm, clip]{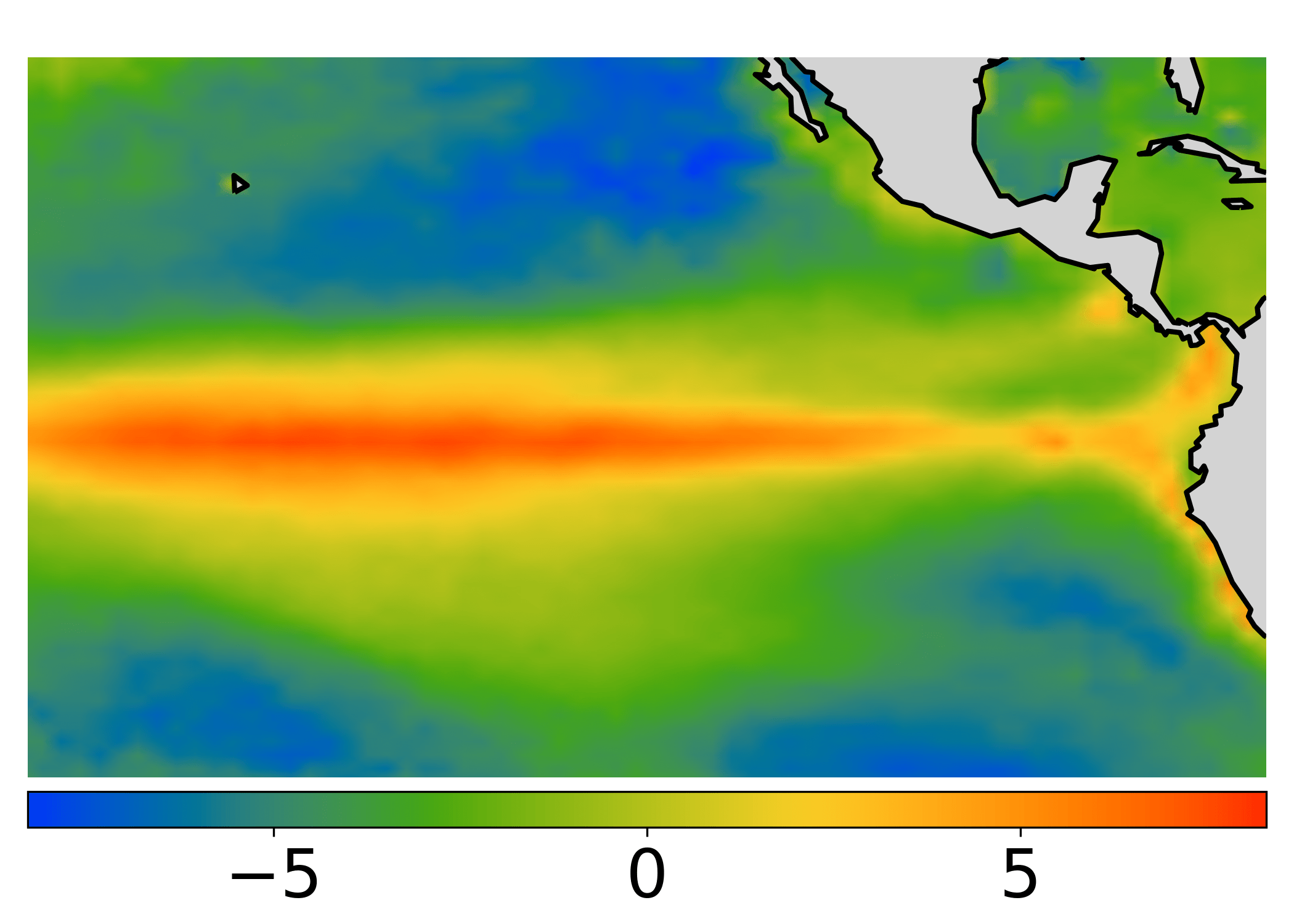} 
			\put(20,-4){\color{black}{temperature fluctuations}}   			
		\end{overpic}\vspace*{+0.2cm}
		\caption{exact ENSO mode}\label{fig:sst_true}
	\end{subfigure}	
\vspace{0.0cm}
\caption{(a) The relevant ENSO region, marked with a rectangule. (b) The true ENSO mode, obtained by exact computation with the full product $A\ttop A$. 
	}\label{fig:sstmid}
\end{figure}
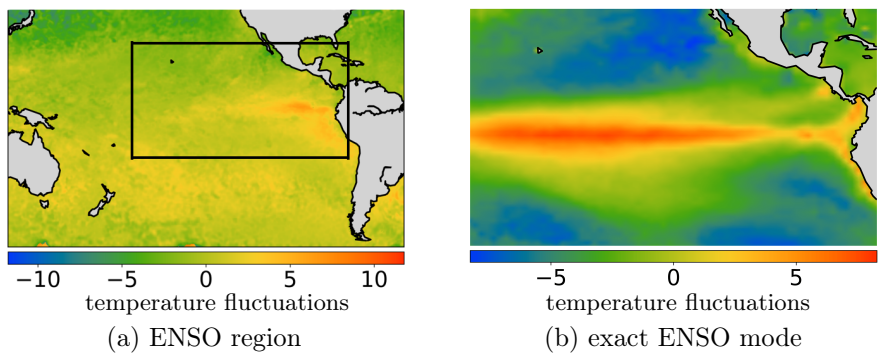

This type of situation is especially common in fields such as climate science and fluid dynamics,
where we may be presented with a very large matrix $A\in\R^{d\times p}$ whose rows 
form a long sequence of ``snapshots''  that represent a dynamical system at time points $1,\dots,d$. 
As a concrete example, we consider satellite recordings of sea surface temperature that have been collected over the time period 1981-2018, and are available from the National Oceanic and Atmospheric Administration~(\citeauthor{NOAA})~\citep[cf.][]{reynolds2002improved}.
More specifically, we deal with a particular subset of the data corresponding to $d=13,\!271$ temporal snapshots at $p=3,\!944$ spatial grid points in the eastern Pacific Ocean, shown in Figure~(\ref{fig:sst_example}). From the standpoint of climate science, this region important for studying the phenomenon known as the El Ni\~{n}o Southern Oscillation (ENSO).

%\vspace{0.5cm}
%
%%%%%%%%%%%%%%%%%%%%%%%%% FIRST LAYOUT

%\vspace{-0.5cm}

\begin{figure}[H]
	\begin{subfigure}{0.45\textwidth}	
		\centering
		\DeclareGraphicsExtensions{.pdf}
		\begin{overpic}[width=0.9\textwidth]{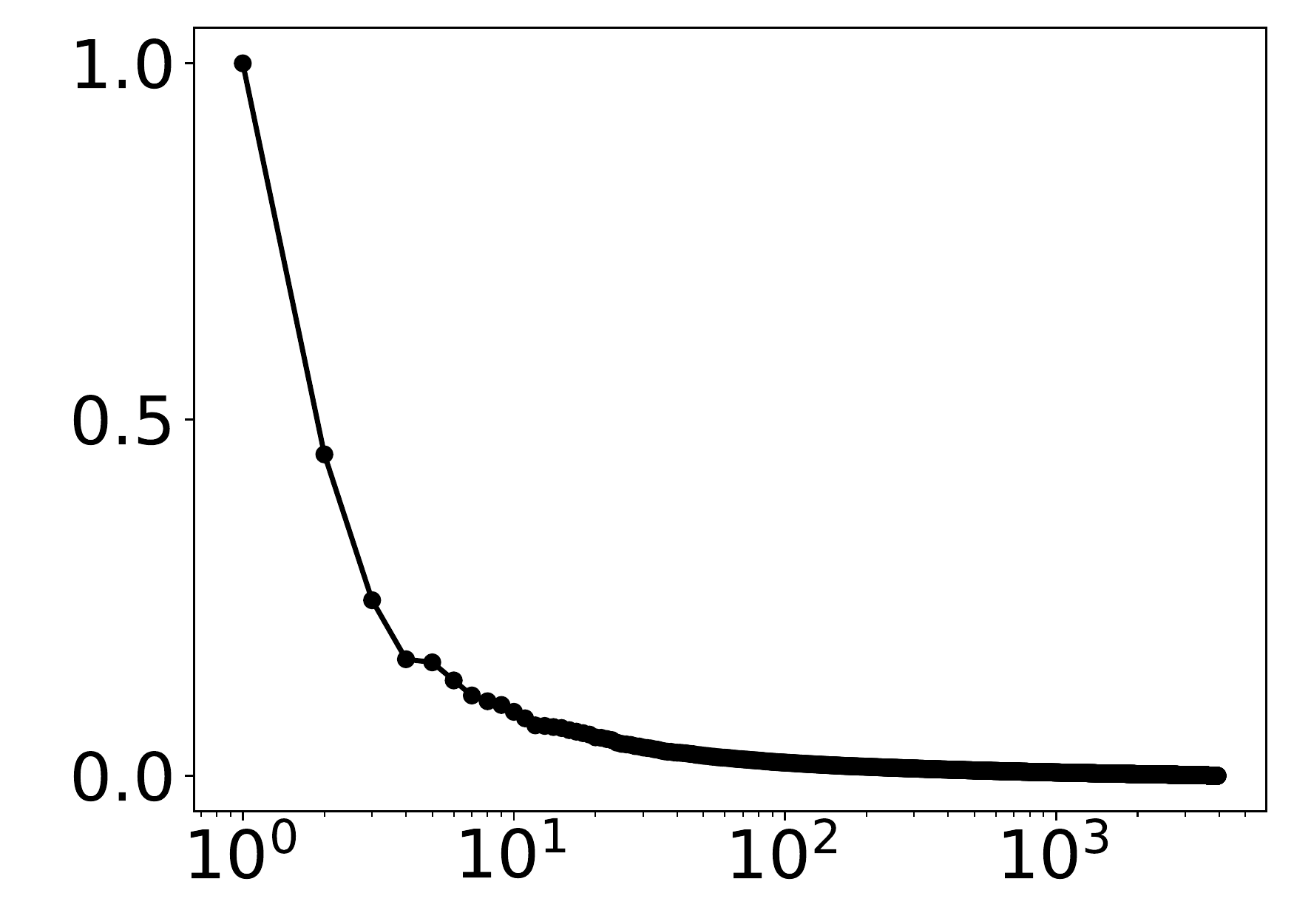} 
			\put(-9,13){\rotatebox{90}{eigenvalue $\lambda_j(A\ttop A)$}}
			\put(44,-3){\color{black}{ index $j$}}
			
			\linethickness{2pt}
			\put(45,49){\color{red} \vector(-1,-2){13}}
			\put(36,51){\color{black}{\footnotesize ENSO mode}}

		\end{overpic}\vspace*{+0.2cm}	
		\caption{spectrum of $A\ttop A$}\label{fig:sst_spectrum}		
	\end{subfigure}	
	~~~~~~
	\begin{subfigure}{0.45\textwidth}	
	\centering
	\DeclareGraphicsExtensions{.pdf}
	\begin{overpic}[width=0.9\textwidth]{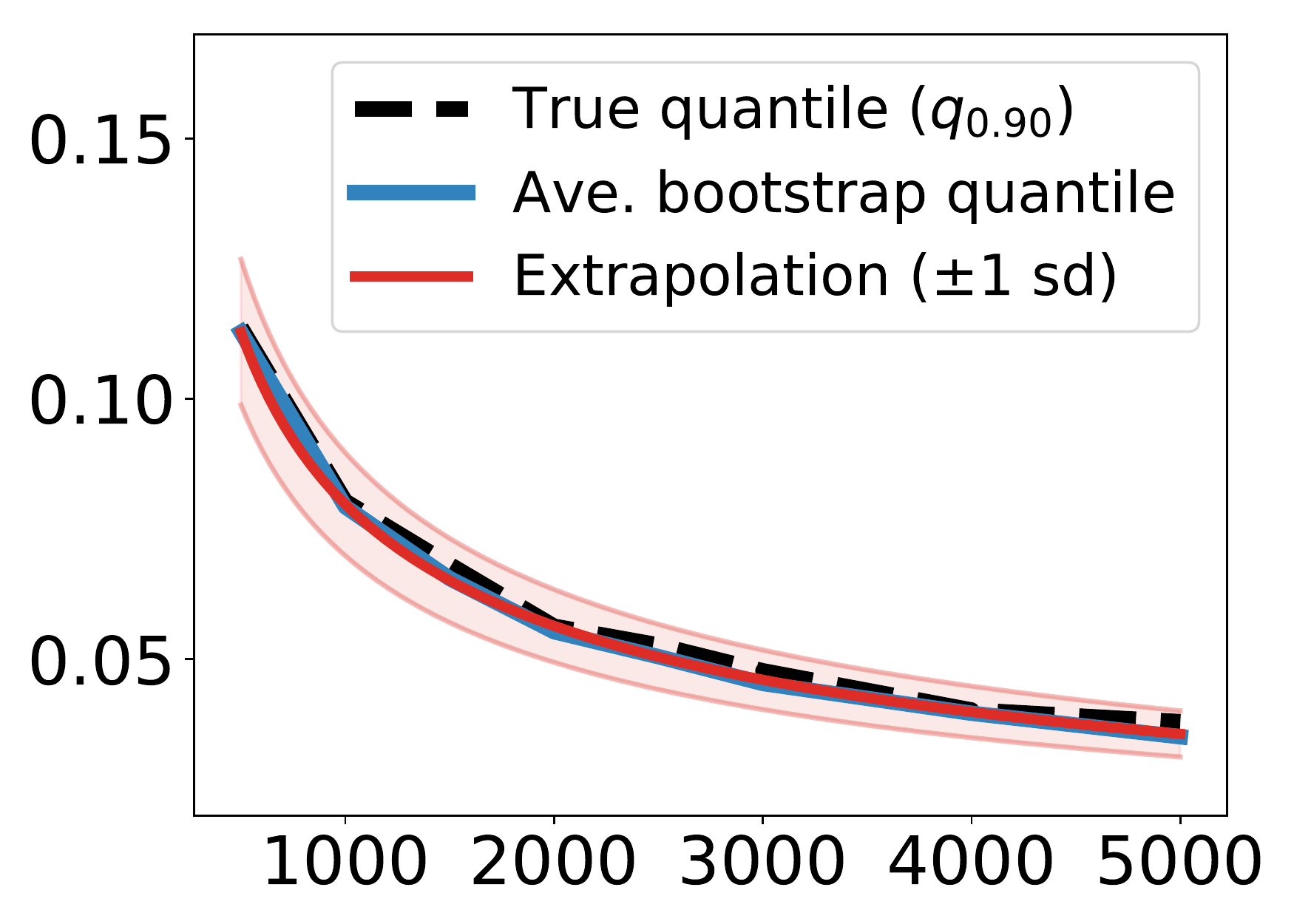} 
		\put(-5,21){\rotatebox{90}{op.~norm error}}
		\put(42,-4){\color{black}{sketch size $n$}}   			
	\end{overpic}\vspace*{+0.2cm}	
	\caption{error estimation}\label{fig:sst_sim}
\end{subfigure}		
	\caption{The left panel displays the decaying eigenvalues of $A\ttop A$, where the x-axis is logarithmic. The right panel demonstrates that the extrapolated and non-extrapolated bootstrap methods accurately estimate the 90\% quantile of the sketching error $\|\tilde A\ttop \tilde A-A\ttop A\|\op$ over a wide range of sketch sizes. In particular, the extrapolation rule gives accurate results at a final sketch size $n_1=5,\!000$ that is 10 times larger than the initial sketch size $n_0=500$.}\label{fig:sst}
\end{figure}
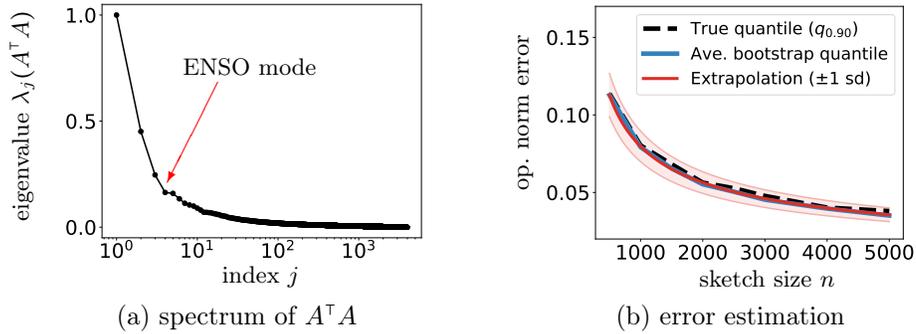

This example is relevant to our discussion of sketching for several reasons. First, the matrix product $A\ttop A$ is of interest because it describes spatial modes of temperature variation through its eigenstructure. In particular, the fourth eigenvector (mode) of $A\ttop A$ identifies the intermittent El Ni\~no and La Ni\~na warming events that are influential global weather patterns, as displayed in Figure~(\ref{fig:sst_true})~\citep{erichson2018sparse}. 
Second, the singular values of $A$ have a natural decay profile, which is illustrated in Figure~(\ref{fig:sst_spectrum}). Lastly, the example demonstrates the need for error estimation in order to guide the choice of sketch size. This can be seen in Figures~(\ref{fig:sst_100}) and~(\ref{fig:sst_800}) below, where it is shown that an insufficient sketch size can heavily distort the ENSO mode in comparison to the exact form given in Figure~(\ref{fig:sst_true}).

\begin{figure}[H]
	\begin{subfigure}{0.45\textwidth}	
		\centering
		\DeclareGraphicsExtensions{.png}
		\begin{overpic}[width=0.99\textwidth]{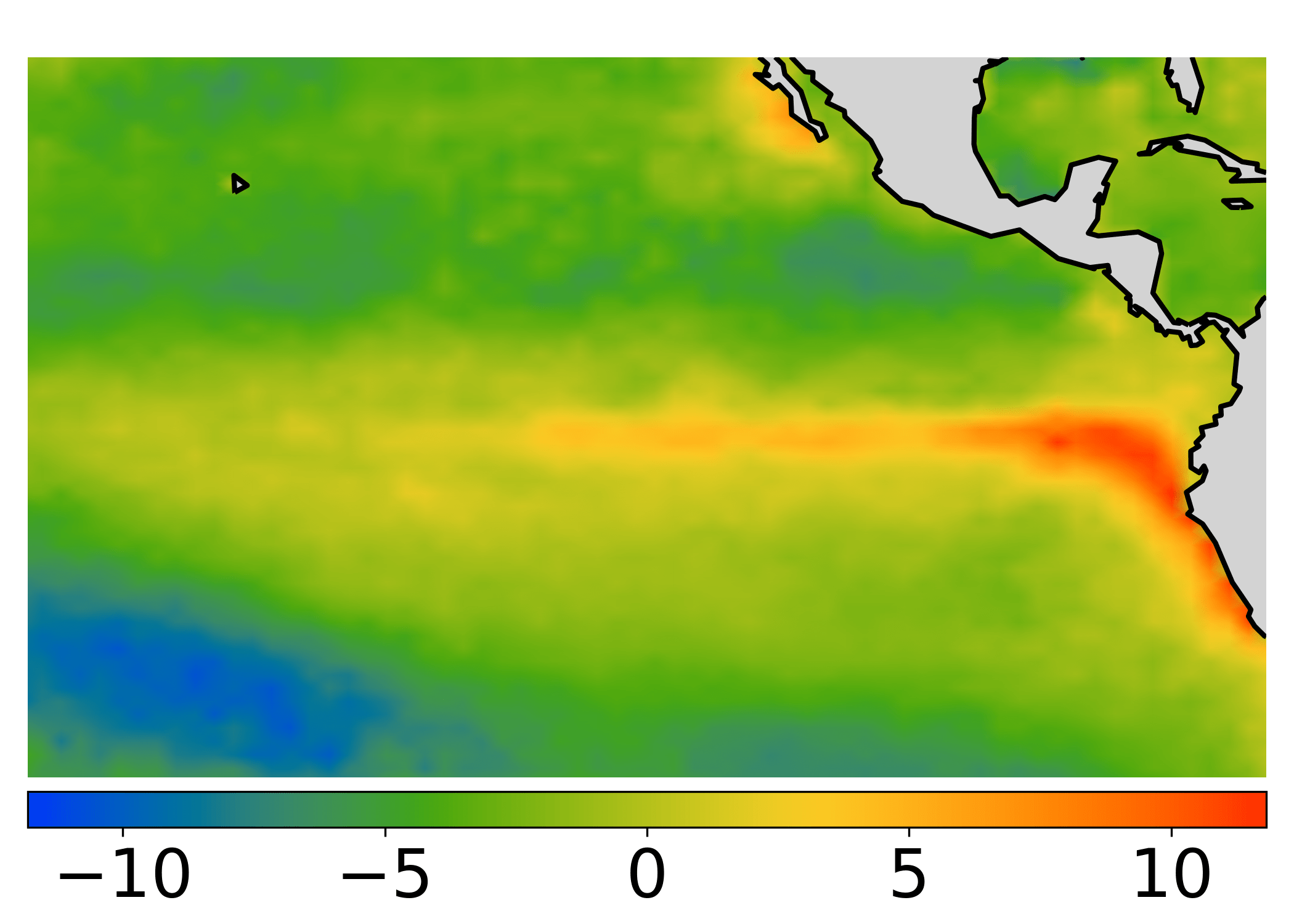} 
			\put(20,-5){\color{black}{\small temperature fluctuations}}   			
		\end{overpic}\vspace*{+0.3cm}
		\caption{ approximate ENSO mode\\ \centering (sketch size $n=500$)}\label{fig:sst_100}	
	\end{subfigure}	
	~
	\begin{subfigure}{0.45\textwidth}	
		\centering
		\DeclareGraphicsExtensions{.png}
		\begin{overpic}[width=0.99\textwidth]{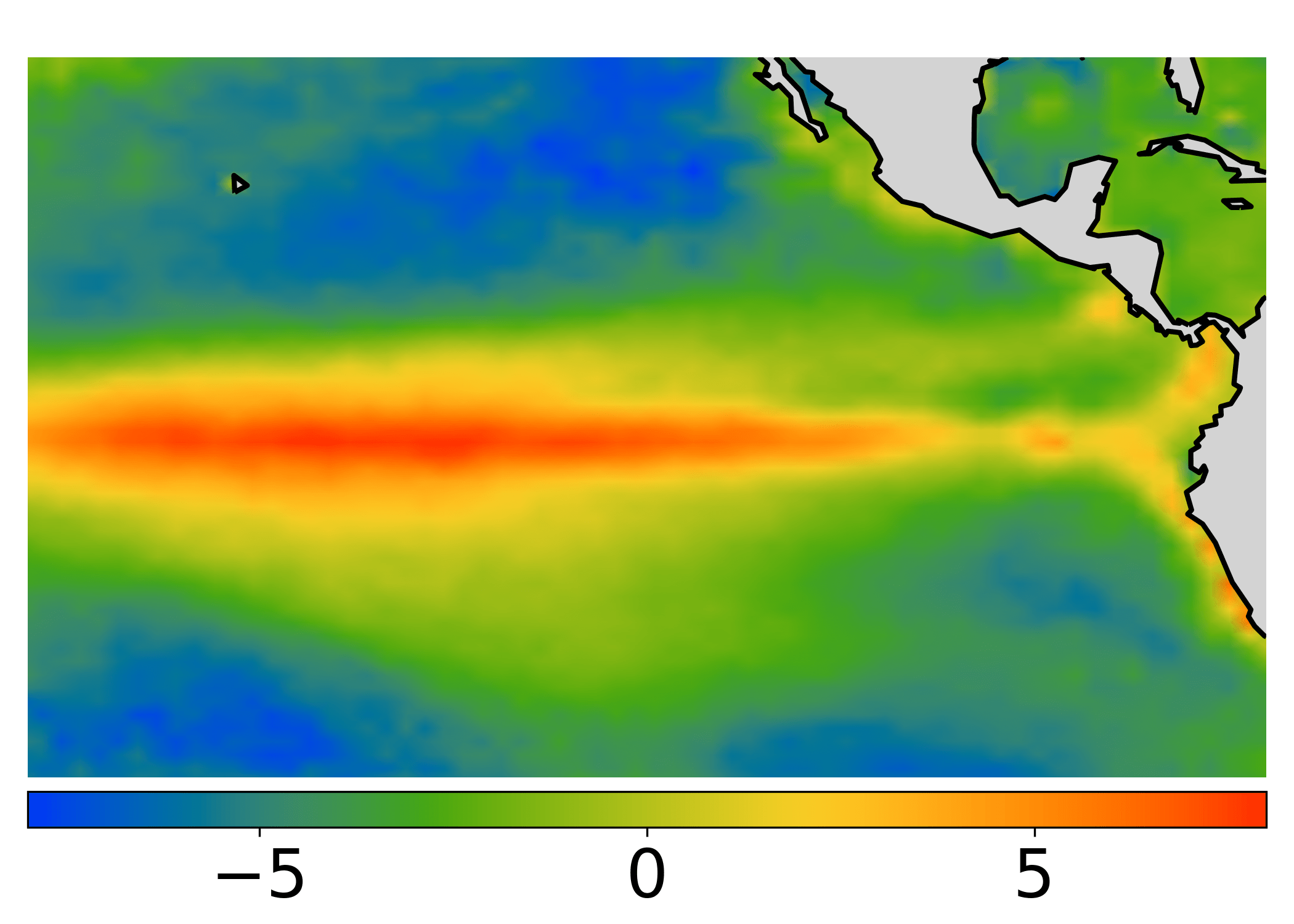} 
			\put(20,-5){\color{black}{\small temperature fluctuations}}   			
		\end{overpic}\vspace*{+0.3cm}
		\centering
		\caption{ approximate ENSO mode\\ \centering (sketch size $n=3,\!000$)}\label{fig:sst_800}
	\end{subfigure}	
	\caption{The left and right panels show approximations to the ENSO mode based on the approximate product $\tilde A\ttop\tilde A$, obtained from Gaussian random projections with sketch sizes $n=500$ and $n=3,\!000$. A comparison with the exact ENSO mode in Figure~(\ref{fig:sst_true}) above shows that an insufficient sketch size can lead to a substantial distortion.
	}\label{fig:sstlast}
\end{figure}

%\vspace{-1.2cm}

To conclude this example, we present numerical results for the bootstrap error estimates. Analogously to~Section~\ref{sec:syntheticsketch}, we consider the task of estimating the 90\% quantile ${\tt{q}}_{0.9}(n)$ of the sketching error, viewed as a function of $n$. The full matrix $A$ is of size $13,\!271\times 3,\!944$, as described earlier, except that it was normalized to satisfy $\sigma_1(A)=1$, so that the results here can be easily compared on the same scale with the previous results in~Section~\ref{sec:syntheticsketch}. Also, the results shown here in Figure~(\ref{fig:sst_sim}) are plotted in the same format, with the number of trials being 1,\!000, the number of bootstrap samples being $B=50$, and the sketching matrices being Gaussian random projections.

From looking at Figure~(\ref{fig:sst_sim}), we see that the performance of the bootstrap in the case of the naturally generated matrix $A$ is very similar to that in the previous cases of synthetic matrices. Namely, the averages of both the extrapolated and non-extrapolated estimates virtually overlap with the true curve, and furthermore, the fluctuations of the extrapolated estimates are well controlled. 
Lastly, the extrapolation rule accurately estimates the quantile value ${\tt{q}}_{0.9}(n_1)$ at a final sketch size $n_1=5,\!000$ that is 10 times larger than the initial sketch size $n_0=500$, which shows the potential of this rule to accelerate computations without sacrificing the quality of estimation.

%\vspace{-0.7cm}

\bibliography{op_norm_bib_shared,sketching}

\newpage

\begin{center}
\large
SUPPLEMENTARY MATERIAL
\end{center}

\appendix

\section{Outline for the proof of Theorem~\ref{THM:MAIN}}
In this section, we define several objects that will recur in our arguments, and then explain how the main components of the proof fit together. As essential pieces of notation, recall the ellipsoidal boundary set in $\R^d$,
$$\mathcal{E} = \{Au \, | \, u\in\mathbb{S}^{p-1}\},$$
as well as its signed version 
$$\Theta=\mathcal{E}\times \{\pm 1\}.$$ 
This set is equipped with the metric $\rho(\theta,\tilde\theta)=\|v-\tilde v\|_2+|s-\tilde s|$ between generic elements $\theta=(v,s)$ and $\tilde\theta=(\tilde v,\tilde s)$.
\subsection{Bootstrap and Gaussian processes} Let $(Z_1^*,\dots,Z_n^*)$ be sampled with replacement from $(Z_1,\dots,Z_n)$, and define the bootstrap counterpart of $\G_n$ as
\begin{equation*}
 \G_n^*(\theta)= \ts\frac{s}{\sqrt n} \displaystyle\sum_{i=1}^n \langle v,Z_i^*\rangle^2-\E[\langle v,Z_i^*\rangle^2|X],
\end{equation*}
where we note that $\E[\langle v,Z_i^*\rangle^2|X]=\ts\frac{1}{n}\sum_{i=1}^n \langle v,Z_i\rangle^2$. This definition of $\G_n^*$ allows $T_n^*$ to be expressed as
$$T_n^*=\sup_{\theta\in\Theta}\,\G_n^*(\theta).$$
In addition, we define $G_n$ as the centered Gaussian process on $\Theta$ whose covariance structure matches that of $\G_n$, 
\begin{equation*}
\cov(G_n(\theta),G_n(\tilde\theta)) \ = \ \cov(\G_n(\theta),\G_n(\tilde \theta)) \text{ \ \ \ for all \ \ \ } \theta,\tilde\theta\in\Theta.
\end{equation*}

\subsection{Subsets of indices} In order to define some special subsets of $\mathcal{E}$ and $\Theta$, let
\begin{equation*}
\begin{split}
\ell_n &=\Big\lceil \big(1\vee \log(n)^3\big)\wedge p\Big\rceil \ \ \ \ \ \text{ and } \ \ \ \ \ k_n =\Big\lceil \Big(\ell_n\vee \log(n)^{\frac{6\beta+4}{\beta-1/2}}\Big)\wedge p\Big\rceil,\\[0.2cm]
\end{split}
\end{equation*}
which always satisfy $1\leq \ell_n\leq k_n\leq p$. Also, let the columns of $V_{k_n}\in\R^{p\times k_n}$ contain the leading $k_n$ right singular vectors of $A$. 
Based on these items, we define $\mathcal{E}_n^{\uparrow}$ as a subset of $\mathcal{E}$ arising from vectors in $\mathbb{S}^{p-1}$ that are ``partially aligned'' with the columns of $V_{k_n}$,
$$\mathcal{E}_n^{\uparrow}=\Big\{Au \ \Big| \  u\in\mathbb{S}^{p-1} \text{ and }  \|V_{k_n}\ttop u\|_2> \ts\frac{1}{2}k_n^{-\beta+1/2}\Big\}.$$
Likewise, by analogy with the definition of $\Theta$, let 
$$\Theta_n^{\uparrow}=\mathcal{E}_n^{\uparrow}\times \{\pm 1\}.$$

The next piece of notation is an $\e$-net for $\Theta_n^{\uparrow}$ with respect to the metric $\rho$. This net is denoted as $\Theta_n^{\uparrow}(\e)\subset\Theta_n^{\uparrow}$ and has the defining property that for any $\theta\in\Theta_n^{\uparrow}$, there is at least one point $\theta'\in\Theta_n^{\uparrow}(\e)$ with $\rho(\theta,\theta')\leq \e$. Throughout the proofs, we will mostly use the particular choice $\e=\e_n$ with
\begin{equation}\label{eqn:endef}
\e_n=n^{-\beta/(6\beta+4)}.
\end{equation}
Lastly, due to classical bounds on the metric entropy of ellipsoids (as recorded in Lemma~\ref{lem:mitjagin}), it is possible to choose an $\e_n$-net for $\Theta_n^{\uparrow}(\e_n)$ with respect to $\rho$  so that its cardinality satisfies
$\log\text{card}(\Theta_n^{\uparrow}(\e_n)) \ \lesssim \ \e_n^{-1/\beta}$.

\subsection{Decomposition into six main terms}
We will bound the Kolmogorov distance between $\mathcal{L}(T_n)$ and $\mathcal{L}(T_n^*|X)$ with six terms,
$$d_{\textup{K}}\big( \mathcal{L}(T_n)\, , \, \mathcal{L}(T_n^*|X)\big) \ \leq \ \I_n \ + \ \II_n \ + \ \IIIo_n \ + \ \III_n \ + \ \IV_n \ + \ \V_n,$$
which are defined below. 
The essential novelty of the proof deals with the four terms $(\I_n, \II_n, \IV_n,\V_n)$, and almost all of the effort will be focused on these. 

\begin{enumerate}
\small
\item Localizing the maximizer of $\G_n$: 
\begin{equation*}
\I_n =d_{\textup{K}}\Big(\mathcal{L}\big(\ts\sup_{\theta\in\Theta} \G_n(\theta)\big) \, , \, \mathcal{L}\big(\ts\sup_{\theta\in\Theta_n^{\uparrow}} \G_n(\theta)\big)\Big).
\end{equation*}
(We use the phrase ``localizing the maximizer of $\G_n$'', because the problem of showing that $\I_n$ is small amounts to showing that the maximizing index for $\G_n$ is likely to fall $\Theta_n^{\uparrow}$.)\\[0.2cm]
\item Discrete approximation of $\G_n$: 
$$\II_n=d_{\textup{K}}\Big(\mathcal{L}\big(\ts\sup_{\theta\in\Theta_n^{\uparrow}} \G_n(\theta)\big) \, , \, \mathcal{L}\big(\ts\sup_{\theta\in\Theta_n^{\uparrow}(\e_n)} \G_n(\theta)\big)\Big)\\[0.2cm]$$
\item Gaussian approximation:
$$\IIIo_n=d_{\textup{K}}\Big(\mathcal{L}\big(\ts\sup_{\theta\in\Theta_n^{\uparrow}(\e_n)} \G_n(\theta)\big)\, , \, \mathcal{L}\big(\ts\sup_{\theta\in\Theta_n^{\uparrow}(\e_n)} G_n(\theta)\big|X\big)\Big)\\[0.2cm]$$
\item Bootstrap approximation:
$$\III_n=d_{\textup{K}}\Big(\mathcal{L}\big(\ts\sup_{\theta\in\Theta_n^{\uparrow}(\e_n)} G_n(\theta)\big)\, , \, \mathcal{L}\big(\ts\sup_{\theta\in\Theta_n^{\uparrow}(\e_n)} \G_n^*(\theta)\big|X\big)\Big)\\[0.2cm]$$
\item Discrete approximation of $\G_n^*$:
$$\IV_n=d_{\textup{K}}\Big( \mathcal{L}\big(\ts\sup_{\theta\in\Theta_n^{\uparrow}(\e_n)} \G_n^*(\theta)\big|X\big) \, , \, \mathcal{L}\big(\ts\sup_{\theta\in\Theta_n^{\uparrow}} \G_n^*(\theta)\big|X\big)\Big)\\[0.2cm]$$
\item Localizing the maximizer of $\G_n^*$:
$$\V_n=d_{\textup{K}}\Big( \mathcal{L}\big(\ts\sup_{\theta\in\Theta_n^{\uparrow}} \G_n^*(\theta)\big|X\big) \, , \,  \mathcal{L}\big(\ts\sup_{\theta\in\Theta} \G_n^*(\theta)\big|X\big)\Big)$$
\end{enumerate}

\noindent Altogether, the six terms are handled consecutively in Appendices~\ref{sec:localizeI} through~\ref{sec:Itilde}, with each appendix corresponding to a different term (except for $\IIIo_n$ and $\III_n$, which are handled together).

\subsection{A general-purpose bound for sample covariance matrices}
Below, we provide a supporting result that will help to streamline some of the proofs later on. Notably, the result can be applied to any sequence of i.i.d.~vectors whose $\ell_2$-norms have well-controlled moments.

\begin{proposition}\label{prop:newcontraction}
Let $\xi_1,\dots,\xi_n\in\R^p$ be i.i.d.~random vectors, and for any $q\geq  3$, define the quantity
\begin{equation}\label{eqn:newrdef}
{\tt{r}}(q) = q\cdot\frac{\Big(\E[\|\xi_1\|_2^{2q}\big]\Big)^{\frac{1}{q}}}{\big\|\E[\xi_1\xi_1\ttop]\big\|\op}.
\end{equation}
 Then, there is an absolute constant $c>0$ such that 
\begin{equation*}
\bigg(\E\bigg\|\ts\frac{1}{n}\displaystyle\sum_{i=1}^n \xi_i\xi_i\ttop -\E[\xi_i\xi_i\ttop]\bigg\|\op^q\bigg)^{1/q} \ \leq \ c\cdot \big\|\E[\xi_1\xi_1\ttop]\big\|\op\cdot \Big(\sqrt{\ts\frac{{\tt{r}}(q)}{n^{1-3/q}}}\,\vee\,\ts\frac{{\tt{r}}(q)}{n^{1-3/q}}\Big).
\end{equation*}
\end{proposition}

\paragraph{Remarks} The proof is given in Appendix~\ref{sec:support}. To convert this result into a convenient high-probability bound, consider the choice $q=\log(n)\vee 3$ and the Chebyshev inequality 
\vspace{-0.2cm}
\begin{equation}\label{eqn:generalchebyshev}
\P\big(\, |Y| \, \geq\, e\|Y\|_q \big) \ \leq \ e^{-q}
\end{equation}
 for a generic random variable $Y$. Then, the event
\begin{equation*}
 \bigg\|\ts\frac{1}{n}\displaystyle\sum_{i=1}^n \xi_i\xi_i\ttop -\E[\xi_i\xi_i\ttop]\bigg\|\op  \ \leq \ c\cdot\big\|\E[\xi_1\xi_1\ttop]\big\|\op\cdot \Big(\sqrt{\ts\frac{{\tt{r}}(q)}{n}}\,\vee\,\ts\frac{{\tt{r}}(q)}{n}\Big)
\end{equation*}
 holds with probability at least $1-\frac{1}{n}$. Also, it will sometimes be useful to consider the special case where the random variable $\|\xi_1\|_2$ can be described in terms of its $\psi_2$-norm. This gives
\begin{equation}\label{eqn:subGspecial}
{\tt{r}}(q)\leq c\cdot q^2\cdot \frac{\big\|\|\xi_1\|_2\big\|_{\psi_2}^2}{\big\|\E[\xi_1\xi_1\ttop]\big\|\op},
\end{equation}
which can be obtained from the facts about Orlicz norms summarized in Lemmas~\ref{lem:orlicz} and~\ref{lem:hansonwright2}.

\section{\!The term $\I_{\lowercase{n}}$: localizing the maximizer of $\G_{\lowercase{n}}$}
\label{sec:localizeI}
The following proposition is the main result of this section, and it will be established with several lemmas later on.
\begin{proposition}\label{prop:Imain}
Suppose that Assumption~\ref{A:model} holds. Then, there is a constant $c>0$ not depending on $n$ such that
\begin{equation*}
\I_n \ \lesssim  \ n^{-\frac{\beta-1/2}{6\beta+4}}\log(n)^c.
\end{equation*}
\end{proposition}

\proof  Here, we only explain how the main pieces fit together, with the details being given in the remainder of this section. Observe that for any $t\in\R$, we have
\begin{equation*}
\Big |\P\Big(\ts\sup_{\theta\in\Theta} \G_n(\theta)\leq t\Big) - \P\Big(\ts\sup_{\theta\in\Theta_n^{\uparrow}}\G_n(\theta)\leq t\Big)\Big | =\P\Big(\mathcal{A}(t)\cap \mathcal{B}(t)\Big),
\end{equation*}
where we define the events
\begin{equation*}
\mathcal{A}(t)=\Big\{\ts\sup_{\theta\in \Theta_n^{\uparrow}}\G_n(\theta)\leq t\Big\} \text{ \ \ \ and \ \ \ } 
\mathcal{B}(t)=\Big\{\ts\sup_{\theta\in \Theta\setminus\Theta_n^{\uparrow}}\G_n(\theta)> t\Big\}.
\end{equation*}
For any pair of real numbers $t_{1,n}$ and $t_{2,n}$ satisfying $t_{1,n}\leq t_{2,n}$, it is straightforward to check that the inclusion 
$(\mathcal{A}(t)\cap \mathcal{B}(t)) \ \subset \ (\mathcal{A}(t_{2,n})\cup \mathcal{B}(t_{1,n}))$
holds simultaneously for all $t\in\R$.
Applying a union bound, and then taking the supremum over $t\in\R$, we obtain
$$\I_n\,\leq \, \P(\mathcal{A}(t_{2,n}))\,+\, \P(\mathcal{B}(t_{1,n})).$$
The difficult part of the proof is carried out below in Lemmas~\ref{lem:lowertail} and~\ref{lem:PB}. In those results, we will determine values of $t_{1,n}$ and $t_{2,n}$ for which the probabilities $\P(\mathcal{A}(t_{2,n}))$ and $\P(\mathcal{B}(t_{1,n}))$ are at most of order $n^{-\frac{\beta-1/2}{6\beta+4}}\log(n)^c$. Furthermore, the chosen values of $k_n$ and $\ell_n$ will ensure that the inequality $t_{1,n}\leq t_{2,n}$ holds for all large $n$.\qed

\begin{remark}\label{rem:pkclarification}
Note that in the special case where $k_n=p$, the matrix $V_{k_n}$ is a square orthogonal matrix, which implies $\Theta_n^{\uparrow}=\Theta$, and as a result, the terms $\I$ and $\V_n$ become exactly 0. Therefore, in the proofs that handle the terms $\I_n$ and $\V_n$, we may assume without loss of generality that $k_n<p$. This small reduction will be needed for specifying how quickly $\ell_n$ and $k_n$ grow as a function of $n$, namely $\ell_n\asymp \log(n)^3$ and $k_n\asymp \log(n)^{\frac{6\beta+4}{\beta-1/2}}$.
\end{remark}

\subsection{Bounding the probability $\P(\mathcal{A}(t_{2,n}))$} In this subsection, we will need to introduce another special subset of $\Theta$. Namely, let $v_1,\dots,v_{\ell_n}\in\R^p$ denote the $\ell_n$ leading right singular vectors of $A$, and define the set 
\begin{equation}\label{eqn:thetaellndef}
\Theta_{\ell_n}=\{(Av_1,1),\dots,(Av_{\ell_n},1)\},
\end{equation}
which satisfies $\Theta_{\ell_n}\subset\Theta_n^{\uparrow}$. Furthermore, we may add the points in $\Theta_{\ell_n}$ to the net $\Theta_n^{\uparrow}(\e_n)$ while preserving the condition $\log\text{card}(\Theta_n^{\uparrow}(\e_n)) \ \lesssim \ \e_n^{-1/\beta}$, since $\ell_n$ grows logarithmically in $n$, whereas $\e_n^{-1/\beta}$ grows algebraically in $n$. For this reason, it will be possible to assume the condition $\Theta_{\ell_n}\subset\Theta_n^{\uparrow}(\e_n)$ without loss of generality in our work below.

\begin{lemma}\label{lem:lowertail} Suppose that Assumption~\ref{A:model} holds. Then, there exists a constant $c_2>0$ not depending on $n$ such that the choice
$$t_{2,n}=c_2\ell_n^{-2\beta}\sqrt{\log(\ell_n)}$$
implies
\begin{equation*}
\P(\mathcal{A}(t_{2,n})) \ \lesssim \ n^{-\frac{\beta-1/2}{6\beta+4}}\log(n)^c.
\end{equation*}
\end{lemma}
\proof  By the definitions of $\mathcal{A}(t_{2,n})$ and $\IIIo_n$, we have
\begin{align*}
\small
\P(\mathcal{A}(t_{2,n}))  
& \ \leq \ \P\Big(\ts\sup_{\theta\in \Theta_n^{\uparrow}(\e_n)}\G_n(\theta)\leq t_{2,n}\Big) \ \ \  \ \ \  \ \  \ \ \ \ \ \ \ (\text{since $\Theta_n^{\uparrow}(\e_n)\subset\Theta_n^{\uparrow}$})\\[0.2cm]
& \ \leq \ \P\Big(\ts\sup_{\theta\in \Theta_n^{\uparrow}(\e_n)}G_n(\theta)\leq t_{2,n}\Big) \ + \ \IIIo_n\\[0.2cm]
& \ \leq \ \P\Big(\ts\sup_{\theta\in \Theta_{\ell_n}}G_n(\theta)\leq t_{2,n}\Big) \ + \ \IIIo_n \ \ \  \ \ \ \ \ (\text{since $\Theta_{\ell_n}\subset\Theta_n^{\uparrow}(\e_n)$}).
\end{align*}
Lemma~\ref{lem:CCK} will show that the term $\IIIo_n$ is at most of order $n^{-\frac{\beta-1/2}{6\beta+4}}\log(n)^c$, and so it remains to control supremum of $G_n$ over $\Theta_{\ell_n}$. This is a substantial task, involving several ingredients that are developed in subsequent lemmas, and so we only explain how the ingredients are combined here. To proceed, define the standardized version of $G_n$ as
$$\bar{G}_n(\theta)=G_n(\theta)/\varsigma_n(\theta) \ \ \ \text{ where} \ \ \ \varsigma_n(\theta)=\sqrt{\var(G_n(\theta))}.$$
Also, define the minimum standard deviation $\varsigma_n^{\circ}=\inf_{\theta\in\Theta_{\ell_n}}\varsigma_n(\theta)$, which is shown to satisfy the following lower bound in Lemma~\ref{lem:minvar},
$$\varsigma_n^{\circ} \ \gtrsim \ \ell_n^{-2\beta}.$$
Based on the definition of $\bar G_n$, it follows that
\begin{equation*}
\small
\begin{split}
\P\Big(\ts\sup_{\theta\in\Theta_{\ell_n}}G_n(\theta)\leq t_{2,n}\Big) & \ 
\leq \ \P\Big(\ts\sup_{\theta\in\Theta_{\ell_n}}\bar{G}_n(\theta) \leq \ts\frac{t_{2,n}}{\varsigma_n^{\circ}}\Big) \\[0.2cm]
& \ \leq   \ \P\Big(\ts\sup_{\theta\in\Theta_{\ell_n}}\bar{G}_n(\theta)\leq \sqrt{\log(\ell_n)}\Big),
\end{split}
\end{equation*}
where the last step can be arranged by the choice of $c_2$ in the definition of $t_{2,n}$.
Next, in Lemmas~\ref{lem:slepian2} and~\ref{lem:covcondition} below,  we show that the last probability satisfies the bound
\begin{equation*}
\begin{split}
  \P\Big(\ts\sup_{\theta\in\Theta_{\ell_n}}\bar{G}_n(\theta)\leq \sqrt{\log(\ell_n)} \Big) & \ \lesssim \ \exp\Big(-\ts\frac{1}{2}\ell_n^{1/3}\Big).
  \end{split}
\end{equation*}
Finally, as explained in the remark above Lemma~\ref{lem:lowertail}, we may assume $k_n<p$, which implies $\ell_n\geq \log(n)^3$, and hence $\exp(-\ts\frac{1}{2}\ell_n^{1/3})\leq n^{-1/2}\leq n^{-\frac{\beta-1/2}{6\beta+4}}\log(n)^c$. 
This completes the proof.\qed

\paragraph{Remark} The following lemma was developed in the paper~\citep[][Lemma B.2]{Lopes:AOS}. Although there are many upper-tail bounds for the maxima of Gaussian processes, there are relatively few lower-tail bounds, which is the notable aspect of this result.

\begin{lemma}\label{lem:slepian2} 
For each integer $m\geq 1$, let $\mathsf{R}=\mathsf{R}(m)$ be a correlation matrix in $\R^{m\times m}$, and let $\mathsf{R}^+=\mathsf{R}^+(m)$ denote the matrix with $(i,j)$ entry given by $\max\{\mathsf{R}_{ij},0\}$. Suppose the matrix $\mathsf{R}^+$ is positive semidefinite for all $m$, and that there are constants $\e_1\in(0,1)$ and $c>0$, not depending on $m$, such that the inequalities
\begin{align*}
\sum_{i\neq j} \mathsf{R}_{ij}^+ & \  \leq  \ c\, m\\[0.2cm]
\ \max_{i\neq j}\mathsf{R}_{ij}^+  & \  \leq \, 1-\e_1
\end{align*}
hold for all $m$.
Lastly, let $(\zeta_1,\dots,\zeta_m)$ be a Gaussian vector drawn from $N(0,\mathsf{R})$. Then, there is a constant $C>0$, not depending on $m$, such that the following inequality holds for all $m\geq 1$,
\begin{equation}\label{eqn:slepian2result}
\P\Big(\max_{1\leq j\leq m} \zeta_j\leq \sqrt{\log(m)}\Big) \ \leq \ C\exp\big(-\ts\frac{1}{2}m^{1/3}\big).
\end{equation}
\end{lemma}

\paragraph{Remark} In essence, the next lemma shows that if we restrict $\G_n$ to the finite set $\Theta_{\ell_n}$, then the correlation matrix of the resulting vector satisfies the conditions of Lemma~\ref{lem:slepian2} (as needed for the completion of the proof of Lemma~\ref{lem:lowertail}).

\begin{lemma}\label{lem:covcondition}
Let the elements of $\Theta_{\ell_n}$ be written as $\{\theta_1,\dots,\theta_{\ell_n}\}$. Also, let $R(\ell_n)\in\R^{\ell_n\times \ell_n}$ denote the correlation matrix of $(\G_n(\theta_1),\dots,\G_n(\theta_{\ell_n}))$, and define the matrix \mbox{$R^+(\ell_n)\in\R^{\ell_n\times \ell_n}$} as
\begin{equation*}
R_{ij}^+(\ell_n)=\max\{R_{ij}(\ell_n),0\}.
\end{equation*}
Under these conditions, it follows that $R^+(\ell_n)$ is positive semidefinite, and there is a constant $c>0$ not depending on $n$ such that
\begin{equation}\label{eqn:firstcorbound}
\sum_{ i\neq j} R_{ij}^+(\ell_n)\ \leq \ c\,\ell_n.
\end{equation}
Furthermore, there is a constant $\e_1\in(0,1)$ not depending on $n$ such that
\begin{equation}\label{eqn:secondcorbound}
\max_{i\neq j} R_{ij}^+(\ell_n)  \ \leq \ 1-\e_1.
\end{equation}
\end{lemma}
\proof Below, we will write $R_{ij}=R_{ij}(\ell_n)$ to ease notation, and likewise for $R^+_{ij}$. By Lemma~\ref{lem:minvar} the following identity 
holds for all $i,j\in\{1,\dots,\ell_n\}$, where we let $u_1,\dots,u_{\ell_n}$ denote the leading $\ell_n$ left singular vectors of $A$,
\begin{equation*}
\begin{split}
\textup{cov}(\G_n(\theta_i),\G_n(\theta_{j})) \ = \
2\sigma(A)_i^2\sigma_{j}^2(A)\Big(1\{i=j\}+\ts\frac{(\kappa-3)}{2} \sum_{l=1}^d \langle e_l,u_i\rangle^2\langle e_l,u_{j}\rangle^2\Big).
\end{split}
\end{equation*}
For distinct $i$ and $j$, this leads to 
\begin{equation}\label{eqn:Rijformula}
R_{ij} =\frac{(\kappa-3) \sum_{l=1}^d \langle e_l,u_i\rangle^2\langle e_l,u_j\rangle^2}{\sqrt{2+(\kappa-3)\sum_{l=1}^d \langle e_l,u_i\rangle^4}\sqrt{2+(\kappa-3)\sum_{l=1}^d \langle e_l,u_j\rangle^4}}.
\end{equation}

In the case when $\kappa\leq 3$, we have $R_{ij}^+=0$ for $i\neq j$, and so the matrix $R^+$ is clearly positive semidefinite. Furthermore, both of the bounds~\eqref{eqn:firstcorbound} and~\eqref{eqn:secondcorbound} hold in this case. To consider the opposite case when $\kappa>3$, observe that $R^+=R$, and so again, the matrix $R^+$ is positive semidefinite. In addition, the formula~\eqref{eqn:Rijformula} implies
\begin{equation*}
 R_{ij}^+  \ \leq \ \ts\frac{\kappa-3}{2} \sum_{l=1}^d \langle e_l,u_i\rangle^2\langle e_l,u_j\rangle^2,
\end{equation*}
and so
\begin{equation*}
\begin{split}
\sum_{1\leq i\neq j\leq \ell_n} \!\!\!R_{ij}^+  & \ \leq \ \ts\frac{\kappa-3}{2} \displaystyle\sum_{i=1}^{\ell_n}\sum_{l=1}^d \langle e_l,u_i\rangle^2 \sum_{j=1}^{\ell_n} \langle e_l,u_j\rangle^2\\[0.3cm]
& \leq \ \ts\frac{\kappa-3}{2} \displaystyle\sum_{i=1}^{\ell_n}\sum_{l=1}^d  \langle e_l,u_i\rangle^2\\[0.3cm]
& \ = \ \ts\frac{\kappa-3}{2} \ell_n,
\end{split}
\end{equation*}
where we have used the fact that $\sum_{j=1}^{\ell_n} \langle e_l,u_j\rangle^2\leq 1$. This proves the bound~\eqref{eqn:firstcorbound}.
Turning to the second bound~\eqref{eqn:secondcorbound}, we may again assume $\kappa>3$. If we let $a=\sum_{l=1}^d \langle e_l,u_i\rangle^4$ and  $b=\sum_{l=1}^d \langle e_l,u_j\rangle^4$, then an application of the Cauchy-Schwarz inequality to~\eqref{eqn:Rijformula} gives
\begin{equation*}
\begin{split}
R_{ij}^+&\leq \frac{\sqrt{ab}}{\sqrt{\big(\ts\frac{2}{\kappa-3}+a\big)\big(\frac{2}{\kappa-3}+b\big)}},\\[0.2cm]
&\leq \frac{\sqrt{ab}}{\sqrt{\big(\ts\frac{2}{\kappa-3}\big)^2+ab}}\\[0.2cm]
&\leq \frac{1}{\sqrt{(\ts\frac{2}{\kappa-3}\big)^2+1}},
\end{split}
\end{equation*}
where the last step follows from the the fact that $ab\leq 1$.
This proves~\eqref{eqn:secondcorbound}.
 \qed

\subsection{Bounding the probability $\P(\mathcal{B}(t_{1,n}))$}

\paragraph{Remark} Before handling the probability $\P(\mathcal{B}(t_{1,n}))$ in Lemma~\ref{lem:PB} below, it is necessary to state a lemma involving a bit of matrix analysis. The proof is straightforward and is hence omitted. For notation, let $A=UDV\ttop$ denote the s.v.d.~\!of $A$, where $U\in\R^{d\times p}$ and $V\in \R^{p\times p}$ have orthonormal columns, and $D\in\R^{p\times p}$ is diagonal. 

\begin{lemma}\label{lem:diagbound}
Fix any $\delta\in(0,1)$, and any symmetric matrix $M\in\R^{d\times d}$. Also, let $V_{k_n}\in\R^{p\times k_n}$ denote the first $k_n$ columns of $V$, and define the map \smash{$\mathsf{T}_{k_n}^{\delta}\!:\R^{p\times p}\to \R^{p\times p}$} that scales the first $k_n$ diagonal entries of a matrix by $\delta$, and leaves all other entries unchanged.
 Then, there is an absolute constant $c>0$ such that
\begin{equation*}
\sup_{\|w\|_2\leq 1, \, \|V_{k_n}\ttop w\|_2\leq \delta}\ |w\ttop A\ttop M Aw| \ \leq \ c\Big\|\mathsf{T}_{k_n}^{\delta}\!\!(D)\,U\ttop M U \, \mathsf{T}_{k_n}^{\delta}\!\!(D) \Big\|\op.
\end{equation*}
\end{lemma}

\noindent We now complete this section with the following bound on $\P(\mathcal{B}(t_{1,n}))$.

\begin{lemma}\label{lem:PB}  Suppose that Assumption~\ref{A:model} holds. Then, there exists a constant $c_1>0$ not depending on $n$ such that the choice
$$ t_{1,n}=c_1 k_n^{-2\beta+1} \log(n)$$
implies
\begin{equation*}
\P(\mathcal{B}(t_{1,n})) \ \lesssim \ \ts\frac{1}{n}.
\end{equation*}
\end{lemma}
\proof 
Let $Z\in\R^{n\times d}$ be the matrix whose rows are $Z_1,\dots,Z_n$. Using the definition of $\Theta_n^{\uparrow}$ and the variational representation of the operator norm, we obtain an expression for $\sup_{\theta\in \Theta\setminus\Theta_n^{\uparrow}}\G_n(\theta)$  by letting $\delta=\ts\frac{1}{2}k_n^{-\beta+1/2}$, 
$$\displaystyle\sup_{\theta\in \Theta\setminus\Theta_n^{\uparrow}}\G_n(\theta) \ = \ \displaystyle\sup_{\|w\|_2\leq 1, \, \|V_{k_n}\ttop w\|_2\leq \delta}\ \sqrt n\, |w\ttop A\ttop \big(\ts\frac{1}{n}Z\ttop Z-I_d\big) Aw|.$$
This expression allows us to apply Lemma~\ref{lem:diagbound}, which gives
\begin{equation}\label{eqn:applymatrixanalysis}
\small
\begin{split}
\sup_{\theta\in \Theta\setminus\Theta_n^{\uparrow}}\G_n(\theta)  & \ \ \leq \ \ c\sqrt{n}\Big\| \mathsf{T}_{k_n}^{\delta}\!\!(D)U\ttop\Big(\ts\frac{1}{n}Z\ttop Z-I_d\Big)U\mathsf{T}_{k_n}^{\delta}\!\!(D)\Big\|\op.
%
%& \ \ = \ \ \sup_{\|w\|_2\leq 1} \ts\frac{c}{\sqrt n}\bigg|\displaystyle\sum_{i=1}^n \langle Z_i,U\mathsf{T}_{k_n}^{\delta}\!\!(D) w\rangle^2-\E[\langle Z_i,U\mathsf{T}_{k_n}^{\delta}\!\!(D) w\rangle^2]\bigg|
\end{split}
\end{equation}
Next, we apply the form of Proposition~\ref{prop:newcontraction} given by the bound~\eqref{eqn:subGspecial}, along with the choices $\xi_i= \mathsf{T}_{k_n}^{\delta}\!\!(D)U\ttop Z_i$ and $q=\log(n)\vee 3$. This gives
\begin{equation}\label{eqn:applygeneral}
\small
\begin{split}
\Big\| \ts\sup_{\theta\in \Theta\setminus\Theta_n^{\uparrow}}\G_n(\theta)\Big\|_{q} 
& \ \lesssim  \  \big\|\E[\xi_1\xi_1\ttop]\big\|\op\cdot \Big(\sqrt{{\tt{r}}(q)}\,\vee\,\ts\frac{{\tt{r}}(q)}{\sqrt n}\Big)\\[0.2cm]
 & \ \lesssim \ \Big(q \cdot \big\|\E[\xi_1\xi_1\ttop]\big\|\op^{1/2}\cdot \big\| \|\xi_1\|_2\big\|_{\psi_2}\Big)\bigvee \Big(\ts\frac{q^2}{\sqrt n}\cdot \big\|\|\xi_1\|_2\big\|_{\psi_2}^2\Big).
\end{split}
\end{equation}
To simplify this bound, note that 
$$\big\|\E[\xi_1\xi_1\ttop]\big\|\op^{1/2} \ = \  \|\mathsf{T}_{k_n}^{\delta}(D)U\ttop\|\op \ \lesssim \ k_n^{-\beta+1/2}.$$
 Also, a background fact in Lemma~\ref{lem:hansonwright2} gives
\begin{equation}\label{eqn:Frobcalc2}
\begin{split}
 \big\|\|\xi\|_2\big\|_{\psi_2} \ & = \ \big\|\|\mathsf{T}_{k_n}^{\delta}\!\!(D)U\ttop Z_1\|_2\big\|_{\psi_2}\\[0.2cm]
& \ \lesssim \  \|\mathsf{T}_{k_n}^{\delta}\!\!(D)U\ttop\|_F \\[0.2cm]
& \ \lesssim \ \bigg(\sum_{j=1}^{k_n} \delta^2\sigma_j^2(A) \ + \ \sum_{j=k_n+1}^{p} \sigma_j^2(A)\bigg)^{1/2}\\[0.2cm]
& \ \lesssim \  k_n^{-\beta+1/2}.
 \end{split}
\end{equation}
So, combining with the earlier bound~\eqref{eqn:applygeneral}, we have
\begin{equation*}
\Big\| \ts\sup_{\theta\in \Theta\setminus\Theta_n^{\uparrow}}\G_n(\theta)\Big\|_{q}  \ \lesssim \ k_n^{-2\beta+1}\log(n),
\end{equation*}
which leads to the stated result by Chebyshev's inequality.\qed

\section{The term $\II_{\lowercase{n}}$: discrete approximation of $\G_{\lowercase{n}}$}\label{sec:II}

\begin{lemma}\label{lem:II}
Suppose Assumption~\ref{A:model} holds. Then, there is a constant $c>0$ not depending on $n$ such that
\begin{equation*}
\II_n \ \lesssim \ n^{-\frac{\beta-1/2}{6\beta+4}}\log(n)^c.
\end{equation*}
\end{lemma}
\proof The approach is based on the fact that the Kolmogorov metric can always be bounded in two parts: a coupling term and an anti-concentration term. More specifically, for any two random variables $\xi$ and $\zeta$ defined on the same probability space, the following inequality holds for any $r>0$,
\begin{equation}\label{eqn:coupleanti}
\dK(\LL(\xi),\LL(\zeta)) \ \leq \ \sup_{t\in\R}\P\big(|\zeta-t|\leq r\big)  \ + \  \P(|\xi-\zeta|\geq r).
\end{equation}
For the present context, we will let $\zeta$ play the role of $\sup_{\theta\in\Theta_n^{\uparrow}(\e_n)}\G_n(\theta)$ and let $\xi$ play the role of $\sup_{\theta\in\Theta_n^{\uparrow}}\G_n(\theta)$.  

With regard to the coupling inequality, it will be shown in Proposition~\ref{prop:key} that if $r$ is proportional to $\e_n^{1-\frac{1}{2\beta}}\log(n)$, then
\begin{equation*}
 \P(|\xi-\zeta|\geq r) \leq \ts\frac{c}{n}.
\end{equation*}
Next, with regard to the anti-concentration inequality, we will approximate $\zeta$ with another random variable, say $\chi$, and then use an anti-concentraiton inequality for $\chi$ instead. To do this, it is simple to verify from the definition of the Kolmogorov metric that the inequality
\begin{equation}\label{eqn:reversecoupleanti}
 \sup_{t\in\R}\P\big(|\zeta-t|\leq r\big)  \ \leq \ \sup_{t\in\R}\P\big(|\chi-t|\leq 2r\big)  \ + \ 2\dK(\LL(\zeta),\LL(\chi)).
\end{equation}
Hence, if we choose $\chi=\sup_{\theta\in\Theta_n^{\uparrow}(\e_n)}G_n(\theta)$, then the distance $\dK(\LL(\zeta),\LL(\chi))$ is the same as $\IIIo_n$, which is shown to be of order $cn^{-\frac{\beta-1/2}{6\beta+4}}\log(n)^c$ in Lemma~\ref{lem:CCK}. Furthermore, by using the stated choice of $r$ and the fact that $\e_n^{1-\frac{1}{2\beta}}=n^{-\frac{\beta-1/2}{6\beta+4}}$, it follows from Lemma~\ref{lemma:Gaussiananti} that
\begin{equation*}
\sup_{t\in\R}\P\big(|\chi-t|\leq 2r\big)  \ \lesssim n^{-\frac{\beta-1/2}{6\beta+4}}\log(n)^c,
\end{equation*}
which completes the proof.\qed

\vspace{0.1cm}

\subsection{A coupling inequality for $\G_n$}
The main goal for this subsection is to establish the following coupling inequality~\eqref{eqn:keycouple} between the random variables $\ts\sup_{\theta\in\Theta_n^{\uparrow}}\G_n(\theta)$ and $\sup_{\theta\in\Theta_n^{\uparrow}(\e_n)}\G_n(\theta)$. To clarify the notation in the result, the supremum in~\eqref{eqn:supdisplay} is interpreted as being over the set of pairs $\{(\theta,\tilde\theta)\in\Theta^2 \, | \, Ма\rho(\theta,\tilde\theta)\leq \e_n\}$.

\begin{proposition}\label{prop:key}
Let $q=\log(n)\vee 3$, and suppose that Assumption~\ref{A:model} holds. Then, 
\begin{equation}\label{eqn:supdisplay}
\begin{split}
\bigg\| \sup_{\rho(\theta,\tilde\theta)\leq \e_n} \big|\G_n(\theta)-\G_n(\tilde \theta)\big| 
& \, \bigg\|_{q} \ \lesssim  \ \ \e_n^{1-\frac{1}{2\beta}}\, \log(n).
 \end{split}
\end{equation}
Furthermore, there is a constant $c>0$ not depending on $n$ such that
\begin{equation}\label{eqn:keycouple}
\P\bigg(\ \Big| \ts\sup_{\theta\in\Theta_n^{\uparrow}}\G_n(\theta) \, - \, \sup_{\theta\in\Theta_n^{\uparrow}(\e_n)}\G_n(\theta)\Big| \ \geq \ c\, \e_n^{1-\frac{1}{2\beta}}\, \log(n)\bigg) \ \leq \ \ts\frac{c}{n}.
\end{equation}

\end{proposition}

\proof We only prove the first statement, because the second statement is essentially a consequence of the Chebyshev inequality~\eqref{eqn:generalchebyshev}. As an initial observation, note that if the indices $\theta=(v,s)$ and $\tilde\theta=(\tilde v,\tilde s)$ satisfy the condition $\rho(\theta,\tilde\theta)\leq \e_n<1$, then the corresponding signs $s$ and $\tilde s$ must be equal.
This leads to the algebraic identity
\begin{equation*}
\begin{split}
\big|\G_n(\theta)-\G_n(\tilde \theta) \big|
& =  \frac{1}{\sqrt n}\bigg|\sum_{i=1}^n \langle v+\tilde v,Z_i\rangle \langle v-\tilde v,Z_i\rangle - \E\big[ \langle v+\tilde v,Z_i\rangle \langle v-\tilde v,Z_i\rangle\big]\bigg|.
\end{split}
\end{equation*}
To rewrite the quadratic forms in terms of a symmetric matrix, let
$$\mathcal{Q}=\ts\frac{1}{2}\Big((v+\tilde v)(v-\tilde v)\ttop +(v-\tilde v)(v+\tilde v)\ttop\Big),$$
so that
\begin{equation*}
\big|\G_n(\theta)-\G_n(\tilde \theta) \big| \ = \ \frac{1}{\sqrt n}\bigg|\sum_{i=1}^n Z_i\ttop \mathcal{Q} Z_i -\E[Z_i\ttop \mathcal{Q} Z_i]\bigg|.
\end{equation*}
Next, let $t>0$ denote a free parameter to be chosen later, and define the vectors
\begin{equation*}
\omega(t)=\ts\frac{1}{2}\Big[t(v+\tilde v)+\ts\frac{1}{t}(v-\tilde v)\Big] \text{ \ \ \ \  and \ \ \ \  } \tilde \omega(t)=\ts\frac{1}{2}\Big[t(v+\tilde v)-\ts\frac{1}{t}(v-\tilde v)\Big].
\end{equation*}
In turn, it can be checked that these vectors give the following representation of $\mathcal{Q}$,
\begin{equation}\label{eqn:Qrep}
\mathcal{Q}=  \omega(t)\omega(t)\ttop -\, \tilde \omega(t)\tilde\omega(t)\ttop,
\end{equation}
which has a certain invariance property, insofar as it holds for every $t>0$, while $\mathcal{Q}$ itself does not depend on $t$.
The utility of this representation is that it will allow us to work with sums of squares, and also, to optimize with respect to the choice of $t$.

To proceed, we will define a particular ellipsoid that contains the vectors $\omega(t)$ and $\tilde\omega(t)$, and then take a supremum over this ellipsoid to derive a stochastic upper bound on  $\sup_{\rho(\theta,\tilde\theta)\leq \e_n} \big|\G_n(\theta)-\G_n(\tilde \theta)\big|$. For this purpose let $A(\e_n)\in\R^{d\times p}$ be the matrix with the same s.v.d.~as $A$, except that the singular value $\sigma_j(A)$ is replaced with $\sqrt 2\min\{\sigma_j(A),\e_n/2\}$ for every \smash{$j\in\{1,\dots,p\}$}. Also, define $A(t,\e_n)\in\R^{d\times 2p}$ as the column concatenation
\begin{equation*}
A(t,\e_n) = \Big[tA \, , \, \ts\frac{1}{t}A(\e_n)\Big].
\end{equation*}
With this matrix in hand, it can be shown that both vectors $\omega(t)$ and $\tilde\omega(t)$ lie in the ellipsoid $A(t,\e_n)(\mathsf{B}_{2p}(2))$, where $\mathsf{B}_{2p}(2)$ denotes the $\ell_2$-ball of radius 2 in $\R^{2p}$. (For the details, see Lemma~\ref{lem:ellipsoiddetails}.) In particular, this ellipsoid does not depend on the indices $\theta$ and $\tilde\theta$ underlying $\omega(t)$ and $\tilde\omega(t)$.
As a result, we have
\small
\begin{equation*}
 \sup_{\rho(\theta,\tilde\theta)\leq \e_n} \big|\G_n(\theta)-\G_n(\tilde \theta)\big|  \  \leq \ \sup_{w\in \mathsf{B}_{2p}(2)} \ts\frac{2}{\sqrt n}\bigg|\displaystyle\sum_{i=1}^n \langle Z_i,A(t,\e_n)w\rangle^2 - \E\big[\langle Z_i,A(t,\e_n)w\rangle^2\big]\bigg|.
 \end{equation*}
\normalsize
We now apply Proposition~\ref{prop:newcontraction} with $\xi_i=A(t,\e_n)\ttop Z_i$, which gives 
 \begin{equation}\label{eqn:keyconseq}
\bigg\| \sup_{\rho(\theta,\tilde\theta)\leq \e_n} \big|\G_n(\theta)-\G_n(\tilde \theta)\big|\bigg\|_{q} \ \lesssim \ 
c\cdot \sqrt n \cdot \big\|A(t,\e_n)\big\|\op^2\cdot \Big(\sqrt{\ts\frac{{\tt{r}}(q)}{n^{1-3/q}}}\,\vee\,\ts\frac{{\tt{r}}(q)}{n^{1-3/q}}\Big).
 \end{equation}
 Due to the choice $q=\log(n)\vee 3$, we have $n^{1-3/q}\asymp n$, and also, the bound~\eqref{eqn:subGspecial} implies
 \begin{equation*}
{\tt{r}}(q) \ \lesssim \ q^2\frac{\big\|\|A(t,\e_n)\ttop Z_1\|_2\big\|_{\psi_2}^2}{\big\|A(t,\e_n)\big\|\op^2}.
\end{equation*}
Furthermore, Lemma~\ref{lem:hansonwright2} gives
\begin{equation}\label{eqn:psi2ell2}
\begin{split}
 \,\Big\| \|A(t,\e_n)\ttop Z_1\|_2\Big\|_{\psi_2} \ \lesssim \ \|A(t,\e_n)\|_F,
\end{split}
 \end{equation}
and then combining with~\eqref{eqn:keyconseq} leads to
 \footnotesize
\begin{equation*}
 \bigg\| \sup_{\rho(\theta,\tilde\theta)\leq \e_n} \big|\G_n(\theta)-\G_n(\tilde \theta)\big|\bigg\|_{q} \ \lesssim \ 
\Big(q\cdot \|A(t,\e_n)\|\op\cdot \|A(t,\e_n)\|_F\Big) \,\bigvee \Big(\ts\frac{q^2}{\sqrt n}\|A(t,\e_n)\|_F^2\Big).
\end{equation*}
\normalsize
Hence, to complete the proof, it remains to bound the norms of $A(t,\e_n)$ and then specify a value of $t$. From the definition of $A(t,\e_n)$ and a short calculation, we have
\begin{equation}\label{eqn:frobcalc}
\begin{split}
\|A(t,\e_n)\|_F & \ \leq \  t\|A\|_F +\ts\frac{1}{t}\|A(\e_n)\|_F\\[0.2cm]
& \ \lesssim \ t \ + \ \ts\frac{1}{t}\e_n^{1-\frac{1}{2\beta}},
\end{split}
\end{equation}
as well as
\begin{equation}\label{eqn:opcalc}
\begin{split}
\|A(t,\e_n)\|\op & \ \leq \ t\|A\|\op+\ts\frac{1}{t}\|A(\e_n)\|\op\\[0.2cm]
& \ \lesssim \ t \ + \ts\frac{1}{t}\e_n.
\end{split}
\end{equation}
Taking $t=\e_n^{\frac{1}{2}-\frac{1}{4\beta}}$ leads to the stated result.
\qed

\subsection{Anti-concentration inequality for $G_n$}

\begin{lemma}\label{lemma:Gaussiananti}
Suppose that Assumption~\ref{A:model} holds, and let $\{\delta_n\}\subset(0,1)$ be any numerical sequence with $\log(1/\delta_n)\lesssim \log(n)$. Then, there is a constant $c>0$ not depending on $n$ such that
\begin{equation*}
\sup_{t\in\R}\,\P\,\bigg(\Big|\ts\sup_{\theta\in\Theta_n^{\uparrow}(\e_n)}G_n(\theta)-t \,\Big|\,\leq \,\delta_n\bigg) \ \lesssim \ \delta_n \log(n)^c.
\end{equation*}
\end{lemma}
\proof For each $\theta\in\Theta$, let
\begin{equation*}
\begin{split}
\varsigma_n(\theta) &=\sqrt{\var(G_n(\theta))}
\end{split}
\end{equation*}
as well as 
$$ \bar\varsigma_n =\sup_{\theta\in\Theta_n^{\uparrow}(\e_n)}\varsigma_n(\theta) \ \ \ \ \text{ and } \ \ \ 
\vunderline{\varsigma}_n =\inf_{\theta\in\Theta_n^{\uparrow}(\e_n)} \varsigma_n(\theta).
$$
In addition, define the expected supremum
$$\mu_n=\E\Big[\ts\sup_{\theta\in\Theta_n^{\uparrow}(\e_n)}G_n(\theta)/\varsigma_n(\theta)\Big].$$
As a consequence of the anti-concentration inequality in Theorem 3 of~\citep{CCK:PTRF}, we have
\begin{equation*}
\small
 \sup_{t\in\R}\, \P\,\bigg(\Big|\ts\sup_{\theta\in\Theta_n^{\uparrow}(\e_n)}G_n(\theta)-t \,\Big|\leq \delta_n\bigg) \ \, \lesssim \,  \ \ts\frac{\bar{\varsigma}_n}{\vunderline{\varsigma}_n^2}\cdot \delta_n\cdot\Big(\mu_n+\sqrt{1\vee\log(\vunderline{\varsigma}_n/\delta_n)}\Big),
\end{equation*}
where $c>0$ is a constant not depending on $n$.
(Note that in the paper~\citep{CCK:PTRF}, the dependence of the bound on $\bar{\varsigma}_n$ and $\vunderline{\varsigma}_n$ is not given explicitly, but a scan through the proof shows that it is sufficient to use a prefactor of $\ts\frac{\bar{\varsigma}_n}{\vunderline{\varsigma}_n^2}$.) To control the dependence on $\bar{\varsigma}_n$ and $\vunderline{\varsigma}_n$, we may use Lemma~\ref{lem:minvar} to obtain
\begin{equation*}
\ts\frac{\bar{\varsigma}_n}{\vunderline{\varsigma}_n^2} \ \lesssim \ \log(n)^c \ \ \ \text{ and }\ \ \  \log(\vunderline{\varsigma}_n/\delta_n) \ \lesssim \ \log(n),
\end{equation*}
for some constant $c>0$ that does not depend on $n$. 

To complete the proof, we must bound $\mu_n$. An initial step is to work with the unstandardized process $G_n$ by using the bound
\begin{equation*}
\begin{split}
\mu_n  \ &  \ \leq  \ \ts\frac{1}{\vunderline{\varsigma}_n}\,\E\Big[\ts\sup_{\theta\in\Theta_n^{\uparrow}(\e_n)}|G_n(\theta)|\Big],\\[0.2cm]
& \ \leq  \ \ts\frac{2}{\vunderline{\varsigma}_n}\,\E\Big[\ts\sup_{\theta\in\Theta_n^{\uparrow}(\e_n)}G_n(\theta)\Big] \ + \ \ts\frac{1}{\vunderline{\varsigma}_n}\,\E\Big[|G_n(\theta_0)|\Big] \ \ \ \ \ (\text{some } \theta_0\in\Theta_n^{\uparrow}(\e_n)),\\[0.2cm]
& \ \leq \ \ts\frac{2}{\vunderline{\varsigma}_n}\,\E\Big[\ts\sup_{\theta\in\Theta_n^{\uparrow}(\e_n)}G_n(\theta)\Big] \ + \ \ts\frac{\bar{\varsigma}_n}{\vunderline{\varsigma}_n},
\end{split}
\end{equation*}
where the second step is a general fact about processes that are symmetric about the origin~\citep[cf.][p.14]{Talagrand:2014}. Next, we will compare $G_n$ with a simpler Gaussian process, whose expected supremum can be analyzed more easily. Due to the Sudakov-Fernique inequality~\citep[Proposition A.2.6]{vanderVaart:Wellner:2000}, if we can construct a centered Gaussian process on $\Theta_n^{\uparrow}(\delta_n)$, say $\Gamma_n(\theta)$, that satisfies the condition
\begin{equation}\label{eqn:slepiancond}
\E\Big[\big(G_n(\theta)-G_n(\tilde\theta)\big)^2\Big] \ \leq \ \E\Big[\big(\Gamma_n(\theta)-\Gamma_n(\tilde\theta)\big)^2\Big]
\end{equation}
for all $\theta,\tilde\theta\in\Theta_n^{\uparrow}(\e_n)$, then the following bound will hold
\begin{equation*}
\E\Big[\ts\sup_{\theta\in\Theta_n^{\uparrow}(\delta_n)} G_n(\theta)\Big]
 \ \leq \  \E\Big[\ts\sup_{\theta\in\Theta_n^{\uparrow}(\delta_n)} \Gamma_n(\theta)\Big].
\end{equation*}
For this purpose, way may apply Lemma~\ref{lem:minvar} to obtain the following formula any $\theta=(v,s)$ and $\tilde\theta=(\tilde v,\tilde s)$,
\begin{equation*}
\begin{split}
\E\Big[\big(G_n(\theta)-G_n(\tilde\theta)\big)^2\Big] & \ = \ 2 \big(\|v\|_2^4+\|\tilde v\|_2^4\big) + (\kappa-3)\big(\|v\|_4^4+\|\tilde v\|_4^4\big)\\[0.2cm]
& \ \ \ \ \ \ -4s\tilde s\langle v,\tilde v\rangle^2-2(\kappa-3)s\tilde s \sum_{l=1}^d  \langle e_l,v\rangle^2\langle e_l,\tilde v\rangle^2.
\end{split}
\end{equation*}
To simplify this expression, let $w,\tilde w\in\R^d$ be vectors with respective $l$th coordinates equal to $\langle e_l,v\rangle^2$ and $\langle e_l,\tilde v\rangle^2$. It can then be checked that
\footnotesize
\begin{equation*}
\begin{split}
\E\Big[\big(G_n(\theta)-G_n(\tilde\theta)\big)^2\Big] & \  = \ 2\Big(s\|v\|_2^2-\tilde s\|\tilde v\|_2^2\Big)^2+(\kappa-3)\|sw-\tilde s\tilde w\|_2^2 + 4s\tilde s\Big(\|v\|_2^2\|\tilde v\|_2^2-\langle v,\tilde v\rangle^2\Big).
\end{split}
\end{equation*}
\normalsize
Letting the three terms on the right be denoted as $J_1,J_2,$ and $J_3$, we can obtain the following bounds by using the fact that all vectors $v\in\mathcal{E}$ satisfy $\|v\|_2^2\leq\|\Sigma\|\op\lesssim 1$. For $J_1$, we have
\begin{equation*}
\begin{split}
 J_1 & \ \leq \  4(s-\tilde s)^2\|v\|_2^4+4(\|v\|_2+\|\tilde v\|_2)^2(\|v\|_2-\|\tilde v\|_2)^2\\[0.2cm]
 & \ \lesssim \ (s-\tilde s)^2 \ + \ \|v-\tilde v\|_2^2.
 \end{split}
 \end{equation*}
Next, for $J_2$, we have
\begin{equation*}
\begin{split}
J_2 & \ \leq \ 2(\kappa-3)_+ \|w\|_2^2 (s-\tilde s)^2 \ + \  2(\kappa-3)_+ \|w-\tilde w\|_2^2\\[0.2cm]
& \ \lesssim \ (s-\tilde s)^2 \ + \ \|v-\tilde v\|_2^2.
\end{split}
\end{equation*}
Lastly, for the third term, it can be checked that $J_3 \ \lesssim \ \|v-\tilde v\|_2^2$, and then combining leads to
\begin{equation*}
\begin{split}
\E\Big[\big(G_n(\theta)-G_n(\tilde\theta)\big)^2\Big] & \ \leq \ c_0(s-\tilde s)^2+c_0\|v-\tilde v\|_2^2,
\end{split}
\end{equation*}
for some constant $c_0>0$ that does not depend on $n$.
Next, we define a centered Gaussian process $\Gamma_n(\theta)$ for any $\theta=(v,s)$ according to
$$\Gamma_n(\theta) \ = \ \sqrt{c_0}\, s\,\zeta_0 \ + \  \sqrt{c_0}\langle v,\zeta\rangle,$$
where $\zeta\in\R^d$ is a standard Gaussian vector, and $\zeta_0\in\R$ is an independent standard Gaussian variable. This yields
\begin{equation*}
 \E\Big[\big(\Gamma_n(\theta)-\Gamma_n(\tilde\theta)\big)^2\Big] \ = \ c_0(s-\tilde s)^2 \ + \ c_0\|v-\tilde v\|_2^2,
\end{equation*}
which shows that the condition~\eqref{eqn:slepiancond} indeed holds. Finally, the expected supremum of $\Gamma_n$ over $\Theta_n^{\uparrow}(\e_n)$ satisfies
\begin{equation*}
\begin{split}
\E\Big[\ts\sup_{\theta\in\Theta_n^{\uparrow}(\e_n)} \Gamma_n(\theta)\Big] & \ \lesssim \ \E\big[|\zeta_0|\big] \ + \ \E\Big[\ts\sup_{\|u\|_2=1}\langle Au,\zeta\rangle\Big]\\[0.2cm]
&\lesssim \ 1 \ + \ \E\Big[\|A\ttop \zeta\|_2\Big]\\[0.2cm]
&\leq \ 1 + \|A\|_F\\[0.2cm]
& \lesssim  \ 1,
\end{split}
\end{equation*}
which completes the proof.\qed

\section{ \small The terms $\IIIo_{\lowercase{n}}$ and $\III_{\lowercase{n}}$: Gaussian and bootstrap approximation}\label{sec:CCK}

The following lemma is obtained as an application of the Gaussian and bootstrap approximation results in~\citep{CCK:AOP}. More recently, the paper~\citep{Deng:2017} demonstrated that under certain conditions, it can be beneficial to avoid a Gaussian approximation step, and to instead directly compare the supremum of an empirical process to its bootstrap counterpart. However, for reasons that seem to be quite technical, it is not clear if this benefit can be carried over to our setting, and accordingly, we proceed with an approach based on Gaussian approximation.
\begin{lemma}\label{lem:CCK}
Suppose that Assumption~\ref{A:model} holds. Then, there is a constant $c>0$ not depending on $n$ such that
\begin{equation}\label{eqn:keyIIIobound}
\IIIo_n \ \leq \ cn^{-\frac{\beta-1/2}{6\beta+4}} \,\log(n)^c,
\end{equation}
and the event
\begin{equation}\label{eqn:keyIIItildebound}
\III_n \ \leq  \   c\, n^{-\frac{\beta-1/2}{6\beta+4}} \,\log(n)^c
\end{equation}
holds with probability at least $1-\ts\frac cn$.
\end{lemma}

\proof We first establish~\eqref{eqn:keyIIIobound}, and then turn to~\eqref{eqn:keyIIItildebound} at the end of the proof. Let $m=\text{card}(\Theta_n^{\uparrow}(\e_n))$, and define i.i.d.~vectors $\xi_1,\dots,\xi_n\in\R^m$ as follows. Let $\{\theta_1,\dots,\theta_m\}$ be an enumeration of $\Theta_n^{\uparrow}(\e_n)$, with $j$th element represented as $\theta_j=(v_j,s_j)$. Next,
for each $i\in\{1,\dots,n\}$ and $j\in\{1,\dots,m\}$, define the random variable
$$\xi_{ij}=k_n^{4\beta-1}s_j\Big(\langle Z_i,v_j\rangle^2 \, - \, \E[\langle Z_j,v_j\rangle^2]\Big).$$
(The scale factor $k_n^{4\beta-1}$ will only play a technical role in order to prevent the variance of $\xi_{ij}$ from becoming too small.)
This definition gives the relation
\begin{equation}\label{eqn:tempcounterpart}
\sup_{\theta\in\Theta_n^{\uparrow}(\e_n)} k_n^{4\beta-1} \G_n(\theta) \ = \ \max_{1\leq j\leq m}\sqrt n \,\bar{\xi}_j
\end{equation}
where $\bar\xi_j=\ts\frac{1}{n}\sum_{i=1}^n \xi_{ij}$.
Although the left side of this relation is a scaled version of $\sup_{\theta\in\Theta_n^{\uparrow}(\e_n)} \G_n(\theta)$, it is important to note that the Kolmogorov metric is scale-invariant, and so the distance between the suprema of $k_n^{4\beta-1} \G_n$ and $k_n^{4\beta-1}G_n$ is equivalent to the distance between the suprema of $\G_n$ and $G_n$.

The proof of~\eqref{eqn:keyIIIobound} is completed by applying Proposition 2.1 in~\citep{CCK:AOP} to $\max_{1\leq j\leq m}\sqrt n\,\bar{\xi}_j$. In order to apply this result, it is enough to note that the vectors $\xi_1,\dots,\xi_n$ are centered, i.i.d., and satisfy the following conditions, which can be verified using Lemmas~\ref{lem:minvar},~\ref{lem:orlicz}, and~\ref{lem:hansonwright2}:
\begin{align}
\min_{1\leq j\leq m} \var(\xi_{1j}) & \ \gtrsim \ 1\label{eqn:CCKconditions1}\\[0.2cm]
\max_{1\leq j\leq m} \|\xi_{1j}\|_{\psi_1}& \ \lesssim \ k_n^{4\beta-1}.\label{eqn:CCKconditions2}\\[0.2cm]
\max_{1\leq j\leq m} \E[|\xi_{1j}|^{2+l}] & \ \lesssim \ k_n^{l(12\beta-3)}\ \ \ \ \text{ for $l\in\{1,2\}$}.\label{eqn:CCKconditions3}
\end{align}
Based on these conditions, as well as $k_n\lesssim \log(n)^{c}$ and $\log(m)\lesssim \e_n^{-1/\beta}$ (by Lemma~\ref{lem:mitjagin}), it follows from Proposition 2.1 in~\citep{CCK:AOP} that there is a constant $c>0$ not depending on $n$ such that
\begin{equation*}
\footnotesize
d_{\textup{K}}\bigg(\mathcal{L}\Big(\ts\sup_{\theta\in\Theta_n^{\uparrow}(\e_n)} k_n^{4\beta-1} \G_n(\theta)\Big)\, , \,\mathcal{L}\Big(\sup_{\theta\in\Theta_n^{\uparrow}(\e_n)} k_n^{4\beta-1} G_n(\theta)\Big)\bigg) \ \lesssim \ n^{-\frac{1}{6}}\e_n^{-\frac{7}{6\beta}} \,\log(n)^c.
\end{equation*}
Substituting in the choice $\e_n=n^{-\frac{\beta}{6\beta+4}}$ leads to~\eqref{eqn:keyIIIobound}.

Finally, to prove~\eqref{eqn:keyIIItildebound}, let $(\xi_1^*,\dots,\xi_n^*)$ be drawn with replacement from $(\xi_1,\dots,\xi_n)$, and let $\bar\xi_j^*=\ts\frac{1}{n}\sum_{i=1}^n \xi_{ij}^*$ for any $j\in\{1,\dots,m\}$. This gives the bootstrap counterpart of the relation~\eqref{eqn:tempcounterpart},
\begin{equation*}
 \sup_{\theta\in\Theta_n^{\uparrow}(\e_n)} k_n^{4\beta-1} \G_n^*(\theta) \ = \ \max_{1\leq j\leq m}\sqrt n( \bar{\xi}_j^* -\bar\xi_j) .
\end{equation*}
Due to this relation, Proposition 4.3 in the aforementioned paper shows that under the conditions~\eqref{eqn:CCKconditions1}-\eqref{eqn:CCKconditions3}, there is a constant $c>0$ not depending on $n$ such that the event
\footnotesize
\begin{equation*}
d_{\textup{K}}\bigg(\mathcal{L}\Big(\ts\sup_{\theta\in\Theta_n^{\uparrow}(\e_n)} k_n^{4\beta-1} \G_n^*(\theta)\Big|X\Big)\, , \,\mathcal{L}\Big(\sup_{\theta\in\Theta_n^{\uparrow}(\e_n)} k_n^{4\beta-1} G_n(\theta)\Big)\bigg) \ \lesssim \ n^{-\frac{1}{6}}\e_n^{-\frac{7}{6\beta}} \,\log(n)^c,
\end{equation*}
\normalsize
holds with probability at least $1-\frac cn$. As before, substituting in the choice $\e_n=n^{-\frac{\beta}{6\beta+4}}$ leads to~\eqref{eqn:keyIIItildebound}.
\qed

\section{The term $\IV_{\lowercase{n}}$: discrete approximation of $\G_{\lowercase{n}}^*$}

\begin{lemma}\label{lem:IItilde} Suppose that Assumption~\ref{A:model} holds. Then, there is a constant $c>0$ not depending on $n$ such that the event
\begin{equation*}
\IV_n \ \leq \ cn^{-\frac{\beta-1/2}{6\beta+4}}\log(n)^c
\end{equation*}
occurs with probability at least $1-\frac cn$.
\end{lemma}
\proof
Recall the Kolmogorov distance can always be bounded in terms of an anti-concentration term and a coupling term, as in~\eqref{eqn:coupleanti}. Using such an approach, we have
\begin{equation*}
\IV_n \ \leq \ \IV_n' \ + \ \IV_n''
\end{equation*}
where we define the following terms for a fixed number $\delta>0$,
\begin{equation*}
\IV_n' \ = \ \sup_{t\in\R}\,\P\bigg(\Big|\ts\sup_{\theta\in\Theta_n^{\uparrow}(\e_n)} \G_n^*(\theta)-t\,\Big| \ \leq \,\delta\, \bigg| \, X\bigg) 
\end{equation*}
and
\begin{equation*}
 \IV_n'' \ = \  \ \P\bigg(\Big|\ts\sup_{\theta\in\Theta_n^{\uparrow}(\e_n)} \G_n^*(\theta) \, - \, \sup_{\theta\in\Theta_n^{\uparrow}} \G_n^*(\theta)\Big| \, \geq \, \delta \, \bigg|\, X\bigg).
\end{equation*}
When $\delta$ is proportional to $\e_n^{1-\frac{1}{2\beta}}\log(n)^c$, we will show in Proposition~\ref{prop:keyboot} below that $\IV_n''$ is at most  $c/n$ with probability at least $1-c/n$.

To address the anti-concentration term, recall the inequality~\eqref{eqn:reversecoupleanti}, which implies
\begin{equation*}
\IV_n' \ \leq \ \sup_{t\in\R}\P\bigg(\Big|\ts\sup_{\theta\in\Theta_n^{\uparrow}(\e_n)} G_n(\theta)-t\Big| \ \leq \ 2\delta \,\bigg)\ + \ 2\III_n.
\end{equation*}
For the stated choice of $\delta$, it is shown in Lemma~\ref{lemma:Gaussiananti} that the first term on the right side is at most of order $\e_n^{1-\frac{1}{2\beta}}\log(n)^c= n^{-\frac{\beta-1/2}{6\beta+4}}\log(n)^c$. Finally, Lemma~\ref{lem:CCK} shows that the event $\III_n\leq c\, n^{-\frac{\beta-1/2}{6\beta+4}} \,\log(n)^c
$ holds with probability at least $1-c/n$, which completes the proof.\qed

\vspace{0.2cm}

\begin{remark}\label{rem:conditionalnorm} To introduce another piece of notation, for any $q\geq 1$ and random variable $\xi$, define the conditional norms
\begin{equation}\label{eqn:conditionalorlicz}
\|\xi\|_{q\,|X} = \ \E\big[\,|\xi|^q\,\big|X\big]^{1/q} \ \ \ \text{ and } \ \ \  \|\xi\|_{\psi_2|X}= \inf\Big\{r>0\,\Big|\, \E[\psi_2(|\xi|/r)|X]\leq 1\Big\}.
\end{equation}
\end{remark}

\begin{proposition}\label{prop:keyboot}
Let $q=\log(n)\vee 3$, and suppose that Assumption~\ref{A:model} holds. Then, there is a constant $c>0$ not depending on $n$ such that the event
\begin{equation}\label{eqn:supdisplayboot}
\begin{split}
\bigg\| \sup_{\rho(\theta,\tilde\theta)\leq \e_n} \big|\G_n^*(\theta)-\G_n^*(\tilde \theta)\big| & \, \bigg\|_{q\,|X} \ \leq  \ \ c\, \e_n^{1-\frac{1}{2\beta}}\log(n)^{5/2}
 \end{split}
\end{equation}
holds with probability at least $1-\ts\frac{c}{n}$, and the event
\begin{equation}\label{eqn:bootcoupling2}
\P\bigg( \Big| \ts\sup_{\theta\in\Theta_n^{\uparrow}}\G_n^*(\theta) \, - \, \sup_{\theta\in\Theta_n^{\uparrow}(\e_n)}\G_n^*(\theta)\Big| \ \geq \ c\, \e_n^{1-\frac{1}{2\beta}}\log(n)^{5/2}\bigg| \, X\bigg) \ \leq \ \ts\frac{c}{n}
\end{equation}
also holds with probability at least $1-\ts\frac{c}{n}$.
\end{proposition}

\proof We only prove~\eqref{eqn:supdisplayboot}, since~\eqref{eqn:bootcoupling2} is essentially a direct consequence. To begin, note that the first half of the of the proof of Proposition~\ref{prop:key} can be repeated to show that
\small
\begin{equation*}
 \sup_{\rho(\theta,\tilde\theta)\leq \e_n} \big|\G_n^*(\theta)-\G_n^*(\tilde \theta)\big|  \  \leq \ \sup_{w\in\mathsf{B}_{2p}(2)} \ts\frac{2}{\sqrt n}\bigg|\displaystyle\sum_{i=1}^n \langle Z_i^*,A(t,\e_n)w\rangle^2 - \E\big[\langle Z_i^*,A(t,\e_n)w\rangle^2\big|X\big]\bigg|,
 \end{equation*}
 \normalsize
where $(Z_1^*,\dots,Z_n^*)$ are i.i.d.~samples with replacement from $(Z_1,\dots,Z_n)$, and we retain the definition of $A(t,\e_n)$ from that proof.
Next, we use the shorthand $\xi_i^*=A(t,\e_n)\ttop Z_i^*$ and apply Proposition~\ref{prop:newcontraction} to the right side above, yielding
 \small
 \begin{equation}\label{eqn:tempsupinc}
\bigg\| \sup_{\rho(\theta,\tilde\theta)\leq \e_n} \big|\G_n^*(\theta)-\G_n^*(\tilde \theta)\big|\bigg\|_{q\,|X} \ \leq \ c\cdot n^{1/2}\cdot \big\|\E[\xi_1^*(\xi_1^*)\ttop|X]\big\|\op\cdot \bigg(\sqrt{\ts\frac{\hat{\tt{r}}(q)}{n^{1-3/q}}}\,\vee\,\ts\frac{\hat{\tt{r}}(q)}{n^{1-3/q}}\bigg),
 \end{equation}
 \normalsize
 where we let
 \begin{equation}\label{eqn:rdef}
\hat{\tt{r}}(q) = q\cdot\frac{\E\big[\|\xi_1^*\|_2^{2q}\big|X\big]^{\frac{1}{q}}}{\ \big\|\E[\xi_1^*(\xi_1^*)\ttop|X]\big\|\op}.
\end{equation}

To simplify the previous bound, note that since $Z_1^*$ is drawn uniformly from $(Z_1,\dots,Z_n)$, it follows that the inequality
 \begin{equation*}
\E\big[ \|\xi_1^*\|_2^{2q}\big|X\big]^{\frac{1}{q}} \ \leq \ \max_{1\leq i\leq n} \|A(t,\e_n)\ttop Z_i\|_2^2
 \end{equation*}
 holds almost surely, and similarly
 \begin{equation*}
 \begin{split}
 \big\|\E[\xi_1^*(\xi_1^*)\ttop|X]\big\|\op & \ = \ \Big\|\ts\frac{1}{n}\sum_{i=1}^n\xi_i\xi_i\ttop\Big\|\op\\[0.2cm]
 & \ \leq \ \max_{1\leq i\leq n} \|A(t,\e_n)\ttop Z_i\|_2^2.
 \end{split}
 \end{equation*}

 Hence, to complete the proof, it suffices to derive a high-probability bound on $\max_{1\leq i\leq n} \|A(t,\e_n)Z_i\|_2^2$. Using the facts in Lemmas~\ref{lem:orlicz} and~\ref{lem:hansonwright2}, as well as the earlier bounds~\eqref{eqn:psi2ell2} and~\eqref{eqn:frobcalc}, we have
 \begin{equation*}
 \begin{split}
 \Big\|\max_{1\leq i\leq n} \|A(t,\e_n)\ttop Z_i\|_2^2\Big\|_{\psi_1} 
& \ \lesssim \ \log(n)\,\big\|\|A(t,\e_n)\ttop Z_1\|_2\big\|_{\psi_2}^2\\[0.2cm]
& \ \lesssim \ \log(n)\,\|A(t,\e_n)\|_F^2\\[0.2cm]
 & \ \lesssim \ \log(n)\, \Big(t+\ts\frac{1}{t}\e_n^{1-\frac{1}{2\beta}}\Big)^2.
\end{split}
\end{equation*}
Therefore, taking $t=\e_n^{-\frac{1}{2}+\frac{1}{4\beta}}$ implies
\begin{equation*}
 \P\bigg(\max_{1\leq i\leq n} \|A(t,\e_n)Z_i\|_2^2 \ \geq \ c\log(n)^2\e_n^{1-\frac{1}{2\beta}} \, \bigg) \ \leq \ \ts\frac{c}{n},
\end{equation*}
which leads to~\eqref{eqn:supdisplayboot} after combining with~\eqref{eqn:tempsupinc}.\qed

\section{The term $\tilde{\mathbf{I}}_{\lowercase{n}}$: localizing the maximizer of $\G_{\lowercase{n}}^*$}\label{sec:Itilde}

\begin{lemma}
Suppose that Assumption~\ref{A:model} holds. Then, there is a constant $c>0$, not depending on $n$, such that the event
\begin{equation*}
\V_n \ \leq \ c n^{-\frac{\beta-1/2}{6\beta+4}}\log(n)^c
\end{equation*}
holds with probability at least $1-c/n$.
\end{lemma}

\proof Observe that
\begin{equation*}
\begin{split}
\V_n &  \ = \ \sup_{t\in\R}\Big |\,  \P\Big(\ts\sup_{\theta\in\Theta_n^{\uparrow}}\G_n^*(\theta)\leq t\, \Big|X\Big) - \P\Big(\ts\sup_{\theta\in\Theta} \G_n^*(\theta)\leq t\, \Big| X\Big) \, \Big |\\[0.2cm]
& \ = \ \sup_{t\in\R} \P\Big(\mathcal{A}'(t)\cap \mathcal{B}'(t)\,\Big| X\Big),
\end{split}
\end{equation*}
where we define the events
\begin{equation*}
\mathcal{A}'(t)=\Big\{\ts\sup_{\theta\in \Theta_n^{\uparrow}}\G_n^*(\theta)\leq t\Big\}
\end{equation*}
and
\begin{equation*}
\mathcal{B}'(t)=\Big\{\ts\sup_{\theta\in \Theta\setminus\Theta_n^{\uparrow}}\G_n^*(\theta)> t\Big\}.
\end{equation*}
By repeating the argument at the beginning of Section~\ref{sec:localizeI}, the following inequality holds for any real numbers $t_{n,1}'$ and $t_{n,2}'$ satisfying $t_{1,n}'\leq t_{2,n}'$,
$$\V_n \ \leq \ \P(\mathcal{A}(t_{2,n}')|X) \ + \ \P(\mathcal{B}(t_{1,n}')|X).$$
To complete the proof, it remains to show there are choices of $t_{1,n}'$ and $t_{2,n}'$ such that $t_{1,n}'\leq t_{2,n}'$ for all large $n$, and the quantities
 $\P(\mathcal{A}'(t_{2,n})|X)$ and $\P(\mathcal{B}'(t_{1,n})|X)$ are at most of order $n^{-\frac{\beta-1/2}{6\beta+4}}\log(n)^c$ with probability at least $1-c/n$. Such choices of $t_{1,n}'$ and $t_{2,n}'$ are established below in Lemma~\ref{lem:lowertailboot}. Note also that the condition $t_{1,n}'\leq t_{2,n}'$ only needs to be established in the case when $k_n<p$, due to the considerations in Remark~\ref{rem:pkclarification}.\qed

\begin{lemma}\label{lem:lowertailboot}
Suppose that Assumption~\ref{A:model} holds. 
Then, there are positive constants $c_1$, $c_2$, and $c$, not depending on $n$,  for which the following statement is true: \\[-0.3cm]

If $t_{1,n}'$ and $t_{2,n}'$ are chosen as
\begin{align}
t_{1,n}'&=c_1 k_n^{-2\beta+1} \log(n)^{3}\label{eqn:t1primedef}\\[0.2cm]
t_{2,n}'& = c_2\ell_n^{-2\beta}\sqrt{\log(\ell_n)},\label{eqn:t2primedef}
\end{align}
then the events
\begin{equation}\label{eqn:Aprimebound}
 \P(\mathcal{A}'(t_{2,n}')|X) \ \leq cn^{-\frac{\beta-1/2}{6\beta+4}}\log(n)^c
\end{equation}
and
\begin{equation*}
 \P(\mathcal{B}'(t_{1,n}')|X) \ \leq \ \ts\frac{c}{n}
\end{equation*}
both occur with probability at least $1-\frac cn$.
\end{lemma}
\proof Based on the definitions of $\IV_n$, $\III_n$, and $\IIIo_n$, we have
\begin{equation*}
\begin{split}
\P(\mathcal{A}'(t_{n,2}')|X) &  \ = \ \P\big(\ts\sup_{\theta\in \Theta_n^{\uparrow}}\G_n^*(\theta)\leq t_{n,2}'\, \big|X\big)\\[0.2cm]
 & \ \leq \ \P\big(\ts\sup_{\theta\in \Theta_n^{\uparrow}}\G_n(\theta)\leq t_{n,2}'\big) +  \ \IV_n \ + \  \III_n \ + \ \IIIo_n \ + \ \II_n.
 \end{split}
\end{equation*}
With regard to the first term of the last line, Lemma~\ref{lem:lowertail} shows that the following holds for a suitable choice of $c_2$ in~\eqref{eqn:t2primedef},
\begin{equation*}
\P\big(\ts\sup_{\theta\in \Theta_n^{\uparrow}}\G_n(\theta)\leq t_{n,2}'\big) \lesssim \ cn^{-\frac{\beta-1/2}{6\beta+4}}\log(n)^c.
\end{equation*}
Combining this with the bounds on $\IV_n$, $\III_n$, $\IIIo_n$, and $\II_n$ in Lemmas~\ref{lem:IItilde},~\ref{lem:CCK}, and~\ref{lem:II}, we reach the stated result in~\eqref{eqn:Aprimebound}.

We now turn to controlling $\P(\mathcal{B}'(t_{1,n}')|X)$. Letting $q=\log(n)\vee 3$, the basic goal is to identify a number $b_n$ that satisfies
\begin{equation}\label{eqn:bndef}
\Big\|\ts\sup_{\theta\in \Theta\setminus\Theta_n^{\uparrow}}\G_n^*(\theta)\Big\|_{q\,|X} \leq \ b_n
\end{equation}
with probability at least $1-\frac cn$. If this can be established, then Chebyshev's inequality will imply that the bound
$$\P\Big(\ts\sup_{\theta\in \Theta\setminus\Theta_n^{\uparrow}}\G_n^*(\theta) \ \geq \ e\,b_n\,\Big| X\Big) \ \leq \ e^{-q}$$
holds with probability at least $1-\frac{c}{n}$. Hence, the number $b_n$ corresponds to $t_{1,n}'$, and also, our choice of $q$ gives $e^{-q}\leq 1/n$.

 To proceed with the details, recall that the s.v.d.~of $A$ is written as $A=UDV\ttop$, and for any $\delta>0$, the map $\mathsf{T}_{k_n}^{\delta}\!:\R^{p\times p}\to\R^{p\times p}$ is defined to act on a matrix by scaling the first $k_n$ diagonal entries by $\delta$ and leaving all other entries unchanged.  Also, let $\delta=\frac{1}{2}k_n^{-\beta+1/2}$ and let $\xi_i^*=\mathsf{T}_{k_n}^{\delta}(D)U\ttop Z_i^*$, where $(Z_1^*,\dots,Z_n^*)$ are i.i.d.~samples with replacement from $(Z_1,\dots,Z_n)$. 
 
 With this notation in place, the argument leading up to~\eqref{eqn:applygeneral} in the proof of Lemma~\ref{lem:PB} can be repeated for the process $\G_n^*$ to obtain
 \small
 \begin{equation}\label{eqn:conditionalqnorm}
  \Big\|\ts\sup_{\theta\in \Theta\setminus\Theta_n^{\uparrow}}\G_n^*(\theta) \Big\|_{q\,|X} \ \leq \ \Big(q \cdot \big\|\E[\xi_1^*(\xi_1^*)\ttop|X]\big\|\op^{1/2}\cdot \big\| \|\xi_1^*\|_2\big\|_{\psi_2|X}\Big)\bigvee \Big(\ts\frac{q^2}{\sqrt n}\cdot \big\|\|\xi_1^*\|_2\big\|_{\psi_2|X}^2\Big).
 \end{equation}
 \normalsize
 To simplify this bound, first notice that if we let $\xi_i=\mathsf{T}_{k_n}^{\delta}(D)U\ttop Z_i$, then  $$\E[\xi_1^*(\xi_1^*)\ttop|X]=\ts\frac{1}{n}\sum_{i=1}^n \xi_i\xi_i\ttop,$$
which leads to 
\begin{equation*}
 \big\|\E[\xi_1^*(\xi_1^*)\ttop|X]\big\|\op^{1/2} \ \leq \ \max_{1\leq i\leq n} \|\mathsf{T}_{k_n}^{\delta}\!(D)U\ttop Z_i\|_2.
\end{equation*}
Similarly, we have
 \begin{equation*}
 \big\|\|\xi_1^*\|_2\big\|_{\psi_2|X} \ \leq \ c\max_{1\leq i\leq n} \|\mathsf{T}_{k_n}^{\delta}\!(D)U\ttop Z_i\|_2,
 \end{equation*}
 for some constant $c>0$ that does not depend on $n$. In turn, we may use the facts about Orlicz norms given in Lemmas~\ref{lem:orlicz} and~\ref{lem:hansonwright2} to obtain
 \begin{equation*}
 \begin{split}
 \Big\|\max_{1\leq i\leq n} \|\mathsf{T}_{k_n}^{\delta}\!(D)U\ttop Z_i\|_2\Big\|_{\psi_2} 
& \ \lesssim \ \sqrt{\log(n)}\big\|\|\mathsf{T}_{k_n}^{\delta}\!(D)U\ttop Z_1\|_2\big\|_{\psi_2}\\[0.2cm]
& \ \lesssim \ \sqrt{\log(n)}\|\mathsf{T}_{k_n}^{\delta}\!(D)\|_F\\[0.2cm]
 %
% & \ = \ \sqrt{\log(n)} \bigg(\ts\sum_{j=1}^{k_n} \delta^2\sigma_j^2(A) \ + \ \ts\sum_{j=k_n+1}^{p} \sigma_j^2(A)\bigg)^{1/2}\\[0.2cm]
%
& \ \lesssim \ \sqrt{\log(n)} k_n^{-\beta+1/2},
\end{split}
\end{equation*}
where the last step re-uses the calculation from~\eqref{eqn:Frobcalc2}. This implies that the bounds
\begin{equation*}
 \big\|\E[\xi_1^*(\xi_1^*)\ttop|X]\big\|\op^{1/2} \ \leq \ c\log(n)k_n^{-\beta+1/2}
\end{equation*}
and
\begin{equation*}
  \big\|\|\xi_1^*\|_2\big\|_{\psi_2|X} \ \leq \ \ c\log(n)k_n^{-\beta+1/2}
\end{equation*}
simultaneously hold with probability at least $1-c/n$. Therefore, combining with the bound~\eqref{eqn:conditionalqnorm} shows that the event
\begin{equation*}
 \Big\|\ts\sup_{\theta\in \Theta\setminus\Theta_n^{\uparrow}}\G_n^*(\theta) \Big\|_{q\,|X} \ \leq \ c\log(n)^3k_n^{-2\beta+1}
\end{equation*}
holds with probability at least $1-c/n$.
Hence, the number $b_n$ in~\eqref{eqn:bndef} may be taken proportional to $\log(n)^3k_n^{-2\beta+1}$, which completes the proof.\qed

\section{Supporting results and proofs}
This section contains the proof of Proposition~\ref{prop:newcontraction}, as well as a lemma summarizing facts about the covariance structure of the process $\G_n$ (Lemma~\ref{lem:minvar}), and technical details involved in the proof of Proposition~\ref{prop:key}.

\subsection{Proof of Proposition~\ref{prop:newcontraction}}\label{sec:support}
Before proceeding directly to the proof of the proposition, we need a preparatory lemma, which is a slightly relaxed version of a result from~\citep[][p.63]{Rudelson:1999}. Also, recall that the Schatten-$q$ norm of a generic real matrix $M$ is defined as $\|M\|_{S_q}=\tr((M\ttop M)^{q/2})^{1/q}$.

\begin{lemma}\label{lem:NC}
Let $x_1,\dots,x_n$ be fixed vectors in $\R^p$, and let $\ve_1,\dots,\ve_n$ be independent Rademacher random variables. Then, there is an absolute constant $c>0$ such that for any $q\geq 2$,
\begin{equation}\label{eqn:new}
\bigg(\E\,\bigg\|\sum_{i=1}^n \ve_i x_ix_i\ttop \bigg\|_{S_q}^q\,\bigg)^{1/q} \ \leq \ c\cdot n^{1/q}\cdot \sqrt{q}\cdot\Big(\max_{1\leq i\leq n}\|x_i\|_2\Big) \cdot \bigg\|\sum_{i=1}^n x_ix_i\ttop\bigg\|\op^{1/2}.
\end{equation}
\end{lemma}
\proof The first step of the proof is to make use of a non-commutative Khinchine inequality established in~\citep{Lust:1986} (see also~\cite{Pisier:2016} Theorem 14.6), 
\begin{equation*}
 \bigg(\E\,\bigg\|\sum_{i=1}^n \ve_i x_ix_i\ttop \bigg\|_{S_q}^q\,\bigg)^{1/q} \ \leq  \ c\sqrt{q}\,\bigg\|\bigg(\sum_{i=1}^n \|x_i\|_2^2\cdot x_ix_i\ttop\bigg)^{1/2}\bigg\|_{S_q}.
\end{equation*}
The next step is to note that any matrix $M$ with rank at most $r$ satisfies $\|M\|_{S_q}\leq r^{1/q}\cdot \|M\|\op $, and that the matrix $\sum_{i=1}^n \|x_i\|_2^2\cdot x_ix_i\ttop$ has rank at most $n$. This leads to
\begin{equation*}
\begin{split}
\bigg\|\bigg(\sum_{i=1}^n \|x_i\|_2^2\cdot x_ix_i\ttop\bigg)^{1/2}\bigg\|_{S_q}& \ \leq \ n^{1/q}\,\bigg\|\bigg(\sum_{i=1}^n \|x_i\|_2^2\cdot x_ix_i\ttop\bigg)^{1/2}\bigg\|\op,
 \end{split}
 \end{equation*}
and it is straightforward to check that this implies the bound~\eqref{eqn:new}.\qed
~\\

\textsc{Proof of Proposition~\ref{prop:newcontraction}.}
The proof extends the approach developed in~\citep{Rudelson:Vershynin:2007} to the case of unbounded random vectors. Using a standard symmetrization argument, we have
\begin{equation}\label{eqn:Lorig}
\bigg(\E\bigg\|\ts\frac{1}{n}\displaystyle\sum_{i=1}^n \xi_i\xi_i\ttop -\E[\xi_i\xi_i\ttop]\bigg\|\op^q\bigg)^{1/q} \ \leq  \  c \bigg(\E\bigg\|\ts\frac{1}{n}\displaystyle\sum_{i=1}^n \ve_i\xi_i\xi_i\ttop \bigg\|\op^q\bigg)^{1/q},
\end{equation}
where $\ve_1,\dots,\ve_n$ are independent Rademacher variables that are also independent of $\xi_1,\dots,\xi_n$. Next, since $\|\cdot\|\op\leq \|\cdot\|_{S_q}$, it follows from Lemma~\ref{lem:NC} that
\small
\begin{equation}\label{eqn:crux}
\bigg(\E\bigg\|\ts\frac{1}{n}\displaystyle\sum_{i=1}^n \ve_i\xi_i\xi_i\ttop \bigg\|\op^q\bigg)^{1/q} \ \leq  \ c \cdot n^{1/q}\cdot\sqrt{\ts\frac{q}{n}}\cdot\E\Big[\max_{1\leq i\leq n}\|\xi_i\|_2^{2q}\Big]^{\frac{1}{2q}} \cdot \bigg(\E\bigg\|\ts\frac{1}{n}\displaystyle\sum_{i=1}^n \xi_i\xi_i\ttop\bigg\|\op^q\bigg)^{\frac{1}{2q}}.
\end{equation}
\normalsize
Using a standard bound for the $L_q$ norm of a maximum, we have
\begin{equation*}
\begin{split}
\E\Big[\max_{1\leq i\leq n}\|\xi_i\|_2^{2q}\Big]^{\frac{1}{2q}} &   \ \leq \ n^{\frac{1}{2q}}\Big(\E[\|\xi_1\|_2^{2q}\big]\Big)^{\frac{1}{2q}}.
\end{split}
\end{equation*}
For the last factor in the bound~\eqref{eqn:crux}, the triangle inequality gives
\begin{equation*}
\begin{split}
\bigg(\E\bigg\|\ts\frac{1}{n}\displaystyle\sum_{i=1}^n \xi_i\xi_i\ttop\bigg\|\op^q\bigg)^{1/q} \ & \leq  \ \bigg(\E\bigg\|\ts\frac{1}{n}\displaystyle\sum_{i=1}^n \xi_i\xi_i\ttop-\E[\xi_i\xi_i\ttop]\bigg\|\op^q\bigg)^{1/q}  \ + \ \big\|\E[\xi_1\xi_1\ttop]\|\op.
\end{split}
\end{equation*}
Hence, if $L$ denotes the left side of~\eqref{eqn:Lorig}, then
\begin{equation*}
 L \ \leq  \ c\cdot n^{\frac{3}{2q}}\cdot\sqrt{\ts\frac{q}{n}}\cdot\Big(\E\big[\|\xi_1\|_2^{2q}\big]\Big)^{\frac{1}{2q}} \cdot\sqrt{L+\big\|\E[\xi_1\xi_1\ttop]\big\|\op}.
\end{equation*}
Finally, by putting
\begin{equation*}
K_1 =  c\cdot n^{\frac{3}{2q}}\cdot\sqrt{\ts\frac{q}{n}}\cdot \Big(\E[\|\xi_1\|_2^{2q}\big]\Big)^{\frac{1}{2q}}\cdot\big\|\E[\xi_1\xi_1\ttop]\big\|\op^{1/2} \ \ \ \ \ \text{ and } \ \ \ \ K_2=\big\|\E[\xi_1\xi_1\ttop]\big\|\op^{-1}
\end{equation*}
we may solve the quadratic inequality
\begin{equation*}
L \ \leq \ K_1\sqrt{K_2L+1}
\end{equation*}
to reach
\begin{equation*}
L \ \leq \ c\Big(K_1 \, \vee \, K_1^2K_2 \Big),
\end{equation*}
which is the stated result.\qed

\subsection{The covariance structure of $\G_n$}
The next result summarizes several the facts about the covariance structure of the process $\G_n$, and provides upper and lower bounds on the parameter $\var(\G_n(\theta))$ over certain subsets of $\Theta$.

\begin{lemma}\label{lem:minvar} Suppose that Assumption~\ref{A:model} holds, and let two generic elements of $\Theta$ be denoted as $\theta=(v,s)$ and $\tilde\theta=(\tilde v,\tilde s)$. Then, 
\begin{equation}\label{eqn:covformula0}
\textup{cov}\big(\G_n(\theta),\G_n(\tilde\theta)\big) \ = \  2s\tilde s\langle v,\tilde v\rangle^2 + (\kappa-3)s\tilde s \sum_{l=1}^d  \langle e_l,v\rangle ^2\langle e_l, \tilde v\rangle^2,
\end{equation}
as well as
\begin{equation}\label{eqn:minvar}
\inf_{\theta\in\Theta_n^{\uparrow}}\sqrt{\var(\G_n(\theta))} \ \gtrsim  \ k_n^{-4\beta+1} \text{ \ \ \ \ \  and \ \ \ \ \ \ } \sup_{\theta\in\Theta_n^{\uparrow}}\sqrt{\var(\G_n(\theta))} \ \lesssim  \ 1.
\end{equation}
Furthermore, if $u_i$ and $v_i$ are the left and right singular vectors of $A$ corresponding to the singular value $\sigma_i(A)$, and we let $\theta_i=(A v_i,1)$, then
\begin{equation}\label{eqn:covformula}
\begin{split}
\textup{cov}(\G_n(\theta_i),\G_n(\theta_{j})) \ = \
2\sigma_i^2(A)\sigma_j^2(A)\bigg(1\{i=j\}+\ts\frac{(\kappa-3)}{2} \displaystyle\sum_{l=1}^d \langle e_l, u_i\rangle^2\langle e_l, u_j\rangle^2\bigg).
\end{split}
\end{equation}
Lastly, if we let $\Theta_{\ell_n}$ be as defined in~\eqref{eqn:thetaellndef}, then
\begin{equation}\label{eqn:varell_n}
\inf_{\theta\in\Theta_{\ell_n}}\sqrt{\var(\G_n(\theta))} \ \gtrsim \  \ell_n^{-2\beta}.
\end{equation}

\end{lemma}

\proof We start with the basic identity
$$\cov(\G_n(\theta),\G_n(\tilde\theta))  = \cov\big(s Z_1\ttop vv\ttop Z_1 \, , \, \tilde s Z_1\ttop \tilde v\tilde v\ttop Z_1\big).$$
Since the entries of $Z_1$ are standardized and independent~with kurtosis $\kappa$, it follows from~\citep[eqn.~9.8.6]{Bai:Silverstein:2010} that 
\begin{equation*}
\cov\big(s Z_1\ttop vv\ttop Z_1 \, , \, \tilde s Z_1\ttop \tilde v\tilde v\ttop Z_1\big) = 2s\tilde s\tr(vv\ttop \tilde v\tilde v\ttop)+(\kappa-3)s\tilde s \sum_{l=1}^d  \langle e_l, v\rangle^2\langle e_l,\tilde v\rangle^2,
\end{equation*}
which implies both~\eqref{eqn:covformula0} and~\eqref{eqn:covformula}.

To establish the lower bound in~\eqref{eqn:minvar}, observe that the previous paragraph gives
\begin{equation}\label{eqn:varformulaetc}
\begin{split}
\var(\G_n(\theta)) & = 2 \|v\|_2^4+(\kappa-3)\|v\|_4^4.
%\[0.2cm]
%
\end{split}
\end{equation}
Due to the assumption $\kappa>1$, there is some fixed $\e_0\in (0,1)$ not depending $n$ such that $\kappa-3\geq -(2-\e_0)$. Consequently, the basic inequality $\|v\|_2\geq \|v\|_4$ implies $(\kappa-3)\|v\|_4^4\geq -(2-\e_0)\|v\|_2^4$, and hence
\begin{equation}\label{eqn:varsimplelower}
\var(\G_n(\theta))\geq \e_0\|v\|_2^4.
\end{equation}
Next, observe that for any $v\in\Theta_n^{\uparrow}$, there is some $w\in\mathbb{S}^{p-1}$ satisfying $v=Aw$ and $\|V_{k_n}\ttop w\|_2 >\frac{1}{2}k_n^{-\beta+1/2}$. Hence, the spectral decomposition $A\ttop A=V D^2V\ttop$, with $v_l$ denoting the $l$th column of $V$ leads to
\begin{equation*}
\begin{split}
\|v\|_2^2  & \ = \ w\ttop A\ttop A w\\[0.2cm]
&\ \geq \ \sum_{l=1}^{k_n} \sigma_l^2(A) \langle w,v_l\rangle ^2\\[0.2cm]
& \ \geq \ \sigma_{k_n}^2(A) \|V_{k_n}\ttop w\|_2^2\\[0.2cm]
& \ \gtrsim \ k_n^{-4\beta+1}.
\end{split}
\end{equation*}
This implies the lower bound in~\eqref{eqn:minvar}. Meanwhile, the upper bound in~\eqref{eqn:minvar} follows from~\eqref{eqn:varformulaetc} and the fact that $\|v\|_2^2\leq \|A\ttop A\|\op\lesssim 1$. Finally, the  lower bound~\eqref{eqn:varell_n} follows from~\eqref{eqn:varsimplelower}.\qed

\subsection{Details for the proof of Proposition~\ref{prop:key}}

\begin{lemma}\label{lem:ellipsoiddetails}
Let the vectors $\omega(t)$ and $\tilde\omega(t)$ as well as the matrix $A(t,\e_n)$ be as defined in the proof of Proposition~\ref{prop:key}. Then, the vectors $\omega(t)$ and $\tilde\omega(t)$ both lie in the ellipsoid $A(t,\e_n)(\mathsf{B}_{2p}(2))$.
\end{lemma}

\proof The proof amounts to showing that the vector $\ts\frac{1}{2}(v+\tilde v)$ lies in $A(\mathsf{B}_p(1))$ and the vector $\ts\frac{1}{2}(v-\tilde v)$ lies in $A(\e_n)(\mathsf{B}_p(1))$. (Note that concatenating two vectors in $\mathsf{B}_p(1)$ yields a vector in $\mathsf{B}_{2p}(2)$.) To proceed, recall that $v$ and $\tilde v$ can be represented as $ v=Aw$ and $\tilde v=A\tilde w$ for some unit vectors $w$ and $\tilde w$. Therefore, the vector
$\ts\frac{1}{2}(v+\tilde v)=A\big(\ts\frac{1}{2}w+\ts\frac{1}{2}\tilde w\big)$
clearly lies in $A(\mathsf{B}_p(1))$.

 Now we turn to $\ts\frac{1}{2}(v-\tilde v)$. From the context of the proof of Proposition~\ref{prop:key}, note that $\|\ts\frac{1}{2}(v-\tilde v)\|_2\leq \ts\frac{1}{2}\e_n$. Also, let the s.v.d.~of $A$ is written as $A=UDV\ttop$, where $U\in\R^{d\times p}$, $D\in\R^{p\times p}$, and $V\in\R^{p\times p}$. Hence, if we let $u_1,\dots,u_p$ denote the columns of $U$, then 
 $$\sum_{l=1}^p \ts\frac{\langle u_l,(v-\tilde v)/2\rangle^2}{\e_n^2/4} \ \leq  \  1.$$
Meanwhile, considering the expression $v-\tilde v=A(w-\tilde w) = UDV\ttop(w-\tilde w)$ gives
\begin{equation*}
\begin{split}
\sum_{l=1}^p \ts\frac{\langle u_l,(v-\tilde v)/2\rangle^2}{\sigma_l^2(A)} & \ = \ \sum_{l=1}^p \ts\frac{\sigma_l^2(A)\langle e_l,V\ttop(w-\tilde w)/2\rangle^2}{\sigma_l^2(A)} \ = \|\ts\frac{1}{2}(w-\tilde w)\|_2^2 \ \leq \ 1.
% \\[0.2cm]
%
\end{split}
\end{equation*}
Hence, if we let $\breve\sigma_l(A)=\sqrt{2}\min\{\sigma_l(A),\e_n/2\}$, then combining leads to
 $$\sum_{l=1}^p \ts\frac{\langle u_l,(v-\tilde v)/2\rangle^2}{\breve\sigma_l^2(A)} \ \leq  \  1.$$
Likewise, if we let $\breve{D}=\text{diag}(\breve\sigma_1(A),\dots,\breve\sigma_p(A))$, then the previous display shows that the vector $\breve D^{-1}U\ttop (v-\tilde v)/2$ lies in the ball $\mathsf{B}_p(1)$, and since $V$ is orthogonal, the vector $x:=V\breve D^{-1}U\ttop(v-\tilde v)/2$ also lies in $\mathsf{B}_p(1)$. In turn, we have
 \begin{equation*}
 \begin{split}
 A(\e_n)x &=\big( U\breve D V\ttop \big) V\breve D^{-1}U\ttop(v-\tilde v)/2\\
 & =UU\ttop (v-\tilde v)/2\\
 &=(v-\tilde v)/2
 \end{split}
 \end{equation*}
where the last step follows from the fact that $v-\tilde v$ lies in the image of $U$. Altogether, this means that $(v-\tilde v)/2$ lies in the ellipsoid $A(\e_n)(\mathsf{B}_p(1))$.\qed

\section{Background results}\label{app:background}

\begin{lemma}[Facts about Orlicz norms]\label{lem:orlicz}
Let $\xi,\xi_1,\dots,\xi_m$ be any sequence of random variables, and let $q\geq 1$, $x>0$, and $r\in\{1,2\}$. Then, there are absolute constants $c,c_0>0$ such that the following hold
		\begin{equation}\label{psi1psi2sq}
		\|\xi^2\|_{\psi_1}\,=\,\|\xi\|_{\psi_2}^2,
		\end{equation}
		\begin{equation}\label{psi1psi2}
		\|\xi\|_{\psi_1}\leq c\,\|\xi\|_{\psi_2},
		\end{equation}
		\begin{equation}\label{lppsi1}
		\|\xi\|_q\leq c\, q^{\frac{1}{r}}\,\|\xi\|_{\psi_r},
		\end{equation}
	\begin{equation}
		\P\big(|\xi|\geq x\big)\leq c \exp\Big\{-c_0\ts\frac{x^r}{\|\xi\|_{\psi_r}^r}\Big\},
		\end{equation}
		and
		\begin{equation}
			\Big\|\max_{1\leq j\leq m} \xi_j\Big\|_{\psi_r} \leq c\,\log(m+1)^{\frac{1}{r}}\,\max_{1\leq j\leq m} \|\xi_j\|_{\psi_r}.
			\end{equation}

\end{lemma}

\proof The first four statements follow from Lemmas 2.7.6 and 2.7.7, as well as Propositions 2.5.2 and 2.7.1 in~\citep{Vershynin:2018}. The fifth statement can be found in Lemma 2.2.2 of~\citep{vanderVaart:Wellner:2000}. 
\qed

\begin{lemma}\label{lem:hansonwright2}
Fix any matrix $M\in\R^{d\times d}$ and vector $v\in\R^d$, and let the random vector $Z_1\in\R^d$ be as in Assumption~\ref{A:model}. Then, there is a constant $c>0$ not depending on $n$ such that 
\begin{equation*}
\big\|\|MZ_1\|_2\big\|_{\psi_2}\ \leq \ c\|M\|_F,
\end{equation*}
and
\begin{equation*}
\Big\|\langle v,Z_1\rangle^2 - \E[\langle v,Z_1\rangle^2]\Big\|_{\psi_1} \ \leq \ c\|v\|_2^2.
\end{equation*}
\end{lemma}
\proof The first statement is a slight reformulation of~\citep[Theorem 6.3.2]{Vershynin:2018}, while the second statement is a special case of~\citep[][Lemma 14]{Lopes:JMLR}.\qed
~\\

\begin{lemma}\label{lem:mitjagin} 
Suppose that Assumption~\ref{A:model} holds. Let $\{\delta_n\}\subset (0,1)$ be a sequence of numbers converging to 0 as $n\to\infty$, and let $\Theta(\delta_n)$ be a minimal $\delta_n$-net for $\Theta$ with respect to the metric $\rho$. Then, 
\begin{equation*}
\log\textup{card}(\Theta(\delta_n)) \ \lesssim \ \delta_n^{-1/\beta}.
\end{equation*}
\end{lemma}
\proof This result is a direct consequence of~(\cite{Kolmogorov:1993}~Theorem XVI), (see also \cite{Kolmogorov:1959}), and so we omit the details.
%  Let $\bar{\mathcal{E}}=\{Au\, |\,u\in\mathsf{B}_p(1)\}$, 
% and let $\bar{\mathcal{E}}(\delta_n)$ denote a corresponding $\delta_n$-net with respect to $\|\cdot\|_2$. It is clear that $\bar{\mathcal{E}}(\delta_n)\times \{\pm 1\}$ can be transformed into a $2\delta_n$-net (not necessarily minimal) for $\Theta$  by sending each point in $\bar{\mathcal{E}}(\delta_n)\times \{\pm 1\}$ to a nearest point in $\Theta$. Thus, it suffices show that
%\begin{equation}
%\log\textup{card}(\bar{\mathcal{E}}(\delta_n))\ \lesssim \ \delta_n^{-1/\beta}.
%\end{equation}
%To proceed, define the set
%$$\mathcal{W} \ = \  \Big\{w\in\R^p \,\Big|\, \ts\sum_{j=1}^p w_j^2/\sigma_j^2(A) \ \leq 1 \Big\},$$
%and let the columns of the matrix $U\in\R^{d\times p}$ contain the left singular vectors of $A$. It is straightforward to check that $\bar{\mathcal{E}}$ is contained in $U(\mathcal{W})$, which is isometric to $\mathcal{W}$. Hence, any $\delta_n$-net for $\mathcal{W}$ leads to a $\delta_n$-net for $\bar{\mathcal{E}}$. Finally, it is a classical fact~(\cite{Kolmogorov:1959},~\cite{Kolmogorov:1993}~Theorem XVI) that a $\delta_n$-net $\mathcal{W}(\delta_n)$ can be chosen for $\mathcal{W}$ so that
%$$\log\text{card}(\mathcal{W}(\delta_n)) \ \lesssim \ \delta_n^{-1/\beta},$$
%which completes the proof.
\qed

\end{document}